\documentclass[10pt]{article}

\usepackage[english]{babel}
\usepackage{amsmath}
\usepackage{amsthm}
\usepackage{bm}
\usepackage{authblk}
\usepackage[margin=1in]{geometry}
\usepackage{amsfonts}
\usepackage[immediate]{silence}
\usepackage{sectsty}
\usepackage{csquotes}
\usepackage[citestyle=numeric-comp, sorting=none, giveninits=true, maxbibnames=10]{biblatex}
\renewbibmacro{in:}{}
\DeclareNameAlias{default}{family-given}
\DeclareFieldFormat[article]{title}{#1}

\usepackage{hyperref}

\newtheorem{theorem}{Theorem}
\newtheorem{proposition}[theorem]{Proposition}
\newtheorem{lemma}[theorem]{Lemma}

\theoremstyle{remark}
\newtheorem*{remark}{Remark}

\theoremstyle{definition}
\newtheorem{definition}[theorem]{Definition}

\numberwithin{equation}{section}

\newcommand{\abs}[1]{\left\lvert#1\right\rvert}

\newcommand{\norm}[1]{\left\lVert#1\right\rVert}

\newcommand{\normz}[1]{\left\lVert#1\right\rVert_{\mathcal{F}_{\nu}^{0,1}}}

\newcommand{\norms}[1]{\left\lVert#1\right\rVert_{\dot{\mathcal{F}}_{\nu}^{s,1}}}

\newcommand{\normzz}[1]{\left\lVert#1\right\rVert_{\mathcal{F}^{0,1}}}

\newcommand{\normss}[1]{\left\lVert#1\right\rVert_{\dot{\mathcal{F}}^{s,1}}}

\newcommand{\normx}[2]{\left\lVert#1\right\rVert_{\dot{\mathcal{F}}_{\nu}^{#2,1}}}

\newcommand{\norml}[1]{\left\lVert#1\right\rVert_{l_{\nu}^{1}}}

\newcommand{\normls}[1]{\left\lVert#1\right\rVert_{l_{\nu}^{s}}}

\DeclareMathOperator{\supp}{supp}

\subsectionfont{\itshape}
\subsubsectionfont{\mdseries\itshape}

\begin{document}


\title{\textbf{Stability of a Two-Phase Stokes Problem with Surface Tension}}


\author{Jae Ho Choi}
\affil{Department of Mathematics, University of Pennsylvania, 
Philadelphia, PA 19104 \\ ORCID: https://orcid.org/0009-0003-8968-863X \\ E-mail address: choijh11@sas.upenn.edu}
\date{}





\maketitle
\newpage


\begin{abstract}
In this work, we study the well-posedness of a system of partial differential equations that model the dynamics of a two-dimensional Stokes bubble immersed in two-dimensional ambient Stokes fluid of the same viscosity that extends to infinity under the effect of surface tension. We assume that the two fluids are immiscible and incompressible and that there is no interfacial jump in the fluid velocity. For this PDE system, a circular fluid bubble is a steady-state solution. Given an initial contour for the fluid bubble which is sufficiently close to a circle, we show that there exists a unique, global-in-time solution. This unique solution decays to a circle exponentially fast, which means that circular fluid bubbles are stable steady-state solutions. We also obtain a result concerning the regularity of the unique solution, that although the initial perturbation around a circular contour is assumed to be of low regularity, any later perturbation becomes real analytic, hence smooth. 
\end{abstract}
\newpage


\setcounter{tocdepth}{1}
\tableofcontents
\newpage


\section{Introduction}
\subsection{Relevant Literature}
In this work, we study well-posedness of a PDE system for the dynamics of a two-dimensional fluid bubble immersed in two-dimensional ambient fluid of the same viscosity that extends to infinity under the effect of surface tension. The two fluids are immiscible and incompressible and there is no interfacial jump in the fluid velocity. The fluids are driven internally by the Stokes equation and interact with one another via surface tension around their interface. From now on, we will refer to this system as a two-phase Stokes problem with surface tension. Stokes problems are sometimes called quasi-stationary approximations of the Navier-Stokes problem because the Stokes equation, which is blind to inertial effects due to fluid motion, is being employed to describe the motion of slow yet non-stationary fluids.

The Navier-Stokes problem with surface tension has attracted mathematicians' attention since the 1980s, starting with the one-phase problem in which an isolated liquid is driven by capillary forces acting on its boundary. The one-phase problem was pioneered in a series of papers by Solonnikov \cite{solonnikov1987solvability, solonnikov1986unsteady, solonnikov1987transient, solonnikov1989unsteady, solonnikov1991solvability, ambrosio2003lectures, solonnikov2003q} and by Mogilevskii and Solonnikov \cite{mogilevskii1992solvability}, in which short-time existence for arbitrary data and long-time existence for small data were established in H\"older and anisotropic Sobolev-Slobodetskii spaces. Since then, well-posedness for the one-phase problem has been established in a multitude of settings, such as the case in which the fluid domain is either bounded, a perturbed infinite layer, or a perturbed half-space \cite{shibata2007free, shibata2008lp, shibata2011report}; and the case in which an infinite viscous incompressible fluid layer is bounded below and above by a solid surface and a free surface, respectively, experiencing surface tension and gravity \cite{allain1987small, beale1984large, beale1985large, tani1996small, tani1995large}.

The two-phase problem gained traction in the 1990s. The first well-posedness results were established by Denisova \cite{denisova1990priori, denisova1994problem} and Denisova and Solonnikov \cite{denisova1994solvability, denisova1995classical}. Since then, well-posedness for the two-phase problem has been established in a number of settings, such as the case in which the system is driven by thermo-capillary convection in bounded domains \cite{tanaka1995two}; and the case in which the free boundary is given as the graph of a function on a hyperplane \cite{pruess2010rayleigh, anger2010two, pruss2011analytic}, sometimes with gravity \cite{pruess2010rayleigh, pruss2011analytic}.

As for the quasi-stationary approximation of the Navier-Stokes problem, the first well-posedness results for one-phase Stokes flow were established by G\"unther and Prokert \cite{gunther1997existence} and Prokert \cite{prokert1999parabolic}. A handful of results concerning the regularity of solutions exist. Escher and Prokert \cite{escher2006analyticity} obtained joint spatial and temporal analyticity of the moving boundary for one-phase Stokes flow with surface tension. G\"unther and Prokert \cite{gunther1997existence} proved short-time existence and uniqueness of a solution for one-phase Stokes flow with a free boundary driven by surface tension in Sobolev spaces of sufficiently high order. Friedman and Reitich \cite{friedman2002quasi} proved joint analyticity of solution for three-dimensional one-phase Stokes flow.

In this work, instantaneous analyticity is established for the unique global-in-time solution of the two-phase
Stokes flow. One of the distinguishing aspects of this work is the analytical framework used to establish well-posedness of the two-phase Stokes problem with surface tension. A large number of aforementioned studies make use of the so-called direct mapping method, where the original physically motivated free boundary problem is transformed into an abstract PDE
problem on a fixed manifold.

\subsection{Connections to Muskat and Peskin Problems}
In oil fields, water is sometimes injected as an emulsifier to reduce the viscosity of crude oil to facilitate its extraction. As our two-phase Stokes problem is driven by surface tension and characterized by a low Reynolds number, it can serve as a rudimentary model to study the behavior of an oil droplet inside water. However, the Muskat model is a more refined and established model in this setting \cite{gancedo2019global,gancedo2023nonlinear}.

Having its roots in petrochemical engineering, the Muskat model is a PDE system describing the dynamics of incompressible fluids of different nature (e.g., oil and water) permeating porous media (e.g., tar sands) under gravity. The fluids’ motions are governed by a momentum equation called Darcy’s law.

Also closely related to our system is the Peskin model, which is a fluid-structure interaction (FSI) model describing the dynamics of a one-dimensional closed elastic string immersed in two-dimensional Stokes fluid. Originally, it emerged as a model for blood flow through heart valves \cite{cameron2021critical}. Being one of the simplest FSI models, it has since been used for other kinds of physical modelling and for building numerical algorithms.

Both the Muskat and Peskin problems have interesting connections to our two-phase Stokes problem with surface tension, which will be explained in depth in Sections \ref{sdlo} and \ref{ppp}.

\subsubsection{Spectral Decomposition of Linearized Operator} \label{sdlo}
Recently, there has been a flurry of mathematical activity on the Muskat model studying its well-posedness. During this process, a multitude of techniques have been devised and employed. Of particular interest to our Stokes problem is spectral decomposition of the linearized operator, which has been applied by Gancedo, Garc\'ia-Ju\'arez, Patel, and Strain to establish global regularity of a two-dimensional Muskat bubble which is unstable under gravity \cite{gancedo2019global}. This technique has also been employed to show global-in-time well-posedness of the Peskin model \cite{garcia2023peskin}. The main idea behind this technique is to linearize the dynamics equation of interest around a steady state solution, which separates the equation into a linear part, which in principle is Fourier analytically a well-behaved operator, and the remainder part, which is “small” in some appropriate sense that depends on a clever selection of the solution space. This linearization is valid only for a small neighborhood around the steady state. For example, the dynamics equation considered by Gancedo, Garc\'ia-Ju\'arez, Patel, and Strain is written in the form
\begin{align} \label{hilbert}
        \partial_{t}g + (-\Delta)^{1/2}g = \mathfrak{R},
\end{align}
where \(\mathfrak{R}\) denotes the part of the equation consisting of terms that are superlinear in \(g\). We note that the principal linear part, \((-\Delta)^{1/2}g\), is the Hilbert transform acting on the spatial derivative of \(g\). In the Fourier space, this equation becomes
\begin{align*}
    \partial_{t}\hat{g}(k) = -\abs{k}\hat{g}(k)+\hat{\mathfrak{R}}(k),
\end{align*}
which clearly reveals that the principal linear part is diagonalized. This explains why the technique is often called spectral decomposition of the linearized operator.

The family of Banach spaces used by Gancedo, Garc\'ia-Ju\'arez, Patel, and Strain that witness the remainder part \(\mathfrak{R}\) to be “small” are
\begin{align} \label{keyspaces}
    \dot{\mathcal{F}}_{\nu(t)}^{s,1} = \biggl\{f: \mathbb{T} \to \mathbb{R} \mid \sum_{k \neq 0} e^{\nu(t)\abs{k}}\abs{k}^{s}\abs{\hat{f}(k)} < \infty \biggr\},
\end{align}
where \(\nu(t) = \frac{t}{1+t}\nu_{0}\) for some \(\nu_{0}>0\). The first superscript, \(s \geq 0\), measures the regularity of functions in the space. The second superscript, \(1\), is simply to indicate that the \(l^{1}\) norm is taken with respect to the wave number \(k\).

Given a sufficiently small initial datum of low regularity describing the initial perturbation of the interface from a circle, which is a steady-state solution, we used spectral decomposition of the linearized operator to establish global-in-time existence and uniqueness of a two-dimensional bubble that satisfies the initial datum and our Stokes problem. The dynamics equation for our Stokes problem is written in the form (\ref{hilbert}), as in the Muskat model studied by Gancedo, Garc\'ia-Ju\'arez, Patel, and Strain. The Peskin model, which bears much similarity with ours, can also be written in that form \cite{cameron2021critical}.

The time-dependent exponential weight in the norm associated with (\ref{keyspaces}) leads to the remarkable property that even though the initial perturbation from a circle is of low regularity, it becomes instantaneously analytic.

\subsubsection{Parametrization} \label{ppp}
The mathematical formulation of the Peskin model is similar to that of ours. However, unlike our model, the Peskin model is driven by the elasticity of the string, which obeys the following general law of elasticity:
\begin{align} \label{elasticity}
    \partial_{\theta}\biggl(T(\abs{\partial_{\theta}\bm{X}})\cdot\frac{\partial_{\theta}\bm{X}}{\abs{\partial_{\theta}\bm{X}}}\biggr)\cdot\abs{\partial_{\theta}\bm{X}}^{-1}.
\end{align}
If we let \(T(\alpha)=\alpha\), then this law reduces to Hooke's law, which is commonly adopted for the analytical study of well-posedness for the Peskin problem. This difference in the mathematical nature of the driving force begets an important analytical consequence.

In the Peskin model, the closed elastic string is parametrized using the Lagrangian coordinate. As the elastic force is critically dependent on this parametrization of the string, it is impossible to choose an arbitrary parametrization to aid in the analysis without fundamentally altering the physical system. This is a major point of difference for the Peskin model from both the Muskat model and our PDE model. In the Muskat model, the normal velocity at the free boundary is well-defined, while the tangential velocity is ill-defined. As the dynamics of the boundary are completely determined by the normal velocity, one can take advantage of the degree of freedom “in the tangential direction” and choose a parametrization that yields nicely to one's analytical framework. In our PDE model, the sole force driving the system is surface tension, which depends exclusively on the geometry of the interface. Therefore, one can employ a convenient parametrization for the interface without affecting its actual dynamics.

To prove our results, we deploy a particular parametrization \cite{hou1994removing} of the fluids' interface that yields nicely to spectral decomposition of the linearized operator. This parametrization is unusual in the sense that the boundary of the fluid bubble is parametrized by the direction of its tangent vector and the length of the boundary. We adopt a certain frame in this parametrization that gives way to our analytical framework, in which the tangent vector is independent of the spatial variable and depends only on time. The same parametrization and change of frame had also been used for a Muskat problem \cite{gancedo2019global}. 

Originally, this particular frame emerged out of a strictly numerical context. Hou, Lowengrub, and Shelley devised it to improve numerical simulation of the motion of the free boundary driven by surface tension between two-dimensional, irrotational, incompressible fluids. Using their novel numerical scheme in which the tangent vector's lack of dependence on the spatial variable removed numerical stiffness, they computed flows that had been unobtainable, such as the motion of the Hele-Shaw interface moving under the competing effects of gravity and surface tension, and discovered new singularity formations, such as the roll-up and collision of vortex sheets with surface tension in two-dimensional Euler flow.

It is possible to cast our Stokes problem as a Peskin model whose force satisfies the general law of elasticity in (\ref{elasticity}) with \(T(\alpha)=1\). The most general setting in which well-posedness has been established for the Peskin problem is when \(T(\alpha)>0\) and \(T'(\alpha)>0\), which means that our Stokes problem corresponds to a degenerate case for which no well-posedness results are available. This suggests that none of the techniques that have been successful in establishing well-posedness of the Peskin model can be used for our Stokes problem, which sheds some light on the strength of our analytical framework.


\subsubsection{Problem Formulation} \label{introduction}
Let \(\Gamma\) be a time-dependent simple closed curve in \(\mathbb{R}^{2}\) that represents the interface between two immiscible fluids. Then the model is given by
\begin{align}
    \mu \Delta \bm{u} - \nabla p &= \bm{0} \quad \mbox{ on \(\mathbb{R}^{2} \setminus \Gamma\)}, \label{1.1} \\
    \nabla \cdot \bm{u} &= 0 \quad \mbox{ on \(\mathbb{R}^{2} \setminus \Gamma\)}, \\
    [\bm{u}] &= \bm{0}, \\
    \left[\Sigma(\bm{u}, p)\bm{n}\right] &= -\gamma\kappa\bm{n}, \label{1.4}
\end{align}
where \(\bm{u}\) and \(p\) denote the fluid velocity and the fluid pressure, respectively; \(\mu\) is the fluid viscosity, which is a constant within each fluid but may differ across the two fluids; \(\Sigma(\bm{u}, p)\) represents the stress tensor for a Newtonian fluid of viscosity \(\mu\); \(\bm{n}\) is the outward-pointing unit normal vector to the interface \(\Gamma\); \(\gamma\) is the surface tension coefficient which is a constant; \(\kappa\) is the signed curvature of the interface; and the notation \([\cdot]\) denotes the value of the variable \(\cdot\) on the boundary as approached from the interior fluid in the normal direction minus the value on the boundary as approached from the exterior fluid in the normal direction. We assume that the two fluids share the same viscosity \(\mu\), which we normalize to \(1\).

In words, this model says that the interior and exterior fluids are incompressible Stokes fluids with no interfacial jump in the fluid velocity and that they are driven by a stress imbalance along the interface given by \(-\gamma\kappa\bm{n}\). The observation that the interfacial force depends exclusively on the geometry of the interface via curvature \(\kappa\) is important, because it allows us to introduce a convenient parametrization for the interface without affecting the physical dynamics of the system.

In this model, there are two unknown variables to solve for: the two-dimensional fluid velocity \(\bm{u}\) and the scalar pressure \(p\). In this work, we study well-posedness of this model in terms of the fluid velocity by imposing on the fluid velocity a single-layer potential satisfying the specified model. Using this ansatz reduces the original problem to that of well-posedness for the PDE system for the interface dynamics. The latter is summarized in the main theorem of this work in Section \ref{maintheorem}. Throughout this work, we may suppress certain expressions' dependence on time \(t\) for readability.

\section{Preliminary Work}
\subsection{Key Function Spaces}
For any \(2\pi\)-periodic function \(f\), its Fourier coefficient is defined as
\begin{align*}
    \mathcal{F}(f)(k)=\frac{1}{2\pi}\int_{-\pi}^{\pi}f(\alpha)e^{-ik\alpha}d\alpha.
\end{align*}
We may sometimes write \(\hat{f}(k)\) to denote the Fourier coefficient of \(f\) with no intended difference in meaning.
We use families of Banach spaces \(\mathcal{F}_{\nu}^{0,1}\) and \(\dot{\mathcal{F}}_{\nu}^{s,1}\), \(s \geq 0\), equipped respectively with norms
\begin{align*}
    \normz{f}&=\sum_{k \in\mathbb{Z}}e^{\nu(t)\abs{k}}\abs{\hat{f}(k)}, \\
    \norms{f}&=\sum_{k\neq 0}e^{\nu(t)\abs{k}}\abs{k}^{s}\abs{\hat{f}(k)},
\end{align*}
where
\begin{align} \label{nu}
    \nu(t)=\frac{t}{1+t}\nu_{0}.
\end{align}
Observe that if \(\nu_{0}>0\), then \(0<\nu'(t) \leq \nu_{0}\). We also use Banach spaces \(\mathcal{F}^{0,1}\) and \(\dot{\mathcal{F}}^{s,1}\), \(s \geq 0\), equipped respectively with norms
\begin{align}
    \normzz{f}&=\sum_{k \in\mathbb{Z}}\abs{\hat{f}(k)}, \label{classicalwienernorm} \\
    \normss{f}&=\sum_{k\neq 0}\abs{k}^{s}\abs{\hat{f}(k)}. \nonumber
\end{align}
The space \(\mathcal{F}^{0,1}\) equipped with the norm in (\ref{classicalwienernorm}) is the classical Wiener algebra, i.e., the space of absolutely convergent Fourier series.
\begin{proposition} (Embeddings.) \label{p3}
For \(0<s_{1}\leq s_{2}\),
\begin{align*}
    \normx{f}{s_{1}} \leq \normx{f}{s_{2}}.
\end{align*}
\end{proposition}
\begin{proposition}(Estimates.) \label{p5}
Let \(n \geq 1\). Then
    \begin{align*}
        \normz{f_{1}f_{2}\cdots f_{n}}&\leq \Pi_{k=1}^{n}\normz{f_{k}}.
    \end{align*}
For \(s > 0\),
    \begin{align*}
        \norms{f_{1}f_{2}\cdots f_{n}} \leq b(n,s)\sum_{j=1}^{n}\norms{f_{j}}\Pi_{k=1,k \neq j}^{n}\normz{f_{k}},
    \end{align*}
where
    \begin{align*}
        b(n,s) =
        \begin{cases}
            1 & 0 \leq s \leq 1, \\
            n^{s-1} & s>1.
        \end{cases}
    \end{align*}
Moreover,
\begin{align*}
    \normx{f_{1}f_{2}}{0} \leq \normx{f_{1}}{0}\normz{f_{2}} + \normz{f_{1}}\normx{f_{2}}{0}.
\end{align*}
\end{proposition}
\begin{remark}
    The first two estimates in Proposition \ref{p5} hold with \(\mathcal{F}_{\nu}^{0,1}\) and \(\dot{\mathcal{F}}_{\nu}^{s,1}\) replaced by \(\mathcal{F}^{0,1}\) and \(\dot{\mathcal{F}}^{s,1}\), respectively. For proof of Proposition \ref{p5}, see Lemma \(5.1\) of~\cite{gancedo2019global}.
\end{remark}
For the following frequently used operator 
\begin{align} \label{m}
    \mathcal{M}(f)(\alpha) = \int_{0}^{\alpha}f(\eta)d\eta-\frac{\alpha}{2\pi}\int_{-\pi}^{\pi}f(\eta)d\eta,
\end{align}
we note that
\begin{align}
    \mathcal{F}(\mathcal{M}(f))(k) =
    \begin{cases}
        -\frac{i}{k}\hat{f}(k) & k \neq 0 \label{operatorMFourier} \\
        \sum_{j \neq 0}\frac{i}{j}\hat{f}(j) & k=0.
    \end{cases}
\end{align}
For \(N \geq 0\), we define high frequency cut-off operators \(\mathcal{J}_{N}\) and \(\mathcal{J}_{N}^{1}\) as
\begin{align}
    \mathcal{F}(\mathcal{J}_{N}f)(k) &= 1_{\abs{k} \leq N}\mathcal{F}(f)(k), \label{highfreqcutoffop2} \\
    \mathcal{F}(\mathcal{J}_{N}^{1}f)(k) &= 1_{\abs{k}\neq 1}1_{\abs{k} \leq N}\mathcal{F}(f)(k). \label{highfreqcutoffop}
\end{align}

\subsection{Boundary Integral Formulation}
For the fluid velocity, let us adopt the single-layer potential
\begin{align} \label{slpvelo}
    u_{j}(\bm{x}) = \frac{1}{4\pi}\int_{\Gamma}\sum_{i=1}^{2}(-\gamma\kappa(s)\bm{n}(s))_{i}G_{ij}(\bm{x} - \bm{y}(s))ds, \quad \bm{x} \in \mathbb{R}^{2},
\end{align}
where \(\bm{u}(\bm{x})=(u_{1}(\bm{x}),u_{2}(\bm{x}))\) and \(G=(G_{ij})\) given by
\begin{align*}
    G_{ij}(\bm{w}) = -\delta_{ij}\log\abs{\bm{w}} + \frac{w_{i}w_{j}}{\abs{\bm{w}}^{2}}
\end{align*}
is the Green's function for two-dimensional unbounded incompressible Stokes flow \cite{pozrikidis1992boundary}. The fluid velocity in our model does not suffer from the Stokes' paradox because the force density \(-\gamma \kappa \bm{n}\) along the interface integrates to \(0\). The single-layer potential ensures that the fluid velocity satisfies equations (\ref{1.1}) through (\ref{1.4}).

\subsection{Interface Parametrization}
We note that the interface's shape is determined entirely by its normal velocity; the tangential velocity can only alter the frame of parametrization. This means that the tangential velocity can be entered into the equations without affecting the interface's shape. We first parametrize the interface \(z(\alpha, t)\), where \(\alpha \in [-\pi,\pi)\). Let us define a tangential angle variable \(\theta\) by writing the tangent vector \(z_{\alpha}(\alpha,t)\) in complex variable notation
\begin{align} \label{tanparam}
    z_{\alpha}(\alpha, t) = \abs{z_{\alpha}(\alpha,t)}e^{i(\alpha + \theta(\alpha, t))}.
\end{align}
We can write
\begin{align} \label{velo}
    z_{t}(\alpha, t) = -U(\alpha, t)\bm{n}(\alpha, t) + T(\alpha, t)\bm{\tau}(\alpha, t),
\end{align}
which in complex variable notation becomes
\begin{align} \label{velocomplex}
    z_{t}(\alpha, t)=U(\alpha,t) ie^{i(\alpha+\theta(\alpha,t))}+T(\alpha,t) e^{i(\alpha+\theta(\alpha,t))},
\end{align}
keeping in mind that in complex variable notation
\begin{align*}
    \bm{\tau}(\alpha,t) &= e^{i(\alpha+\theta(\alpha,t))}, \\
    \bm{n}(\alpha,t) &= -ie^{i(\alpha+\theta(\alpha,t))}.
\end{align*}
After differentiating (\ref{tanparam}) with respect to \(t\) and then differentiating (\ref{velocomplex}) with respect to \(\alpha\), we equate their real and imaginary parts to derive evolution equations for the interface in terms of \(\theta\) and \(\abs{z_{\alpha}(\alpha,t)}\):
\begin{align}
    \abs{z_{\alpha}(\alpha,t)}_{t} &= -U(\alpha, t)-U(\alpha, t)\theta_{\alpha}(\alpha, t) + T_{\alpha}(\alpha, t) \label{intf1}, \\
    \theta_{t}(\alpha, t) &= \frac{1}{\abs{z_{\alpha}(\alpha,t)}}\biggl(U_{\alpha}(\alpha, t) + T(\alpha, t) + T(\alpha, t)\theta_{\alpha}(\alpha, t)\biggr). \label{intf2}
\end{align}
A particularly useful frame of parametrization can be chosen by requiring that the tangential speed \(T(\alpha,t)\) be of the form
\begin{equation}
    T(\alpha,t) = \int_{0}^{\alpha}(1 + \theta_{\eta}(\eta,t))U(\eta,t)d\eta - \frac{\alpha}{2\pi}\int_{-\pi}^{\pi}(1 + \theta_{\eta}(\eta,t))U(\eta,t)d\eta + T(0,t), \label{T}
\end{equation}
where \(T(0,t)\) is a number that depends on \(t\), which allows for a frame change. This frame of parametrization ensures that \(\abs{z_{\alpha}(\alpha,t)}\) is independent of \(\alpha\), i.e.,
\begin{equation*}
    \abs{z_{\alpha}(\alpha, t)} = \frac{1}{2\pi}\int_{-\pi}^{\pi}\abs{z_{\alpha}(\eta, t)}d\eta = \frac{L(t)}{2\pi},
\end{equation*}
where \(L(t)\) is the length of the interface at time \(t\). Using this tangential speed formula, (\ref{intf1}) and (\ref{intf2}) can be rewritten as
\begin{align}
    L_{t}(t) &= -\int_{-\pi}^{\pi}(1+\theta_{\alpha}(\alpha))U(\alpha)d\alpha \label{lengthevol} \\
    \theta_{t}(\alpha, t) &= \frac{2\pi}{L(t)}U_{\alpha}(\alpha) + \frac{2\pi}{L(t)}T(\alpha)(1 + \theta_{\alpha}(\alpha)). \label{thetaevol}
\end{align}
The use of this particular frame of parametrization for a fluid interface was pioneered by~\cite{hou1994removing} in the context of removing numerical stiffness in implementing interfacial flows with surface tension. From now on, we will refer to it as Hou-Lowengrub-Shelley (HLS) parametrization.

\subsection{The Interface Length \(L(t)\)}
We can derive an analytical expression for \(L(t)\) from the incompressibility of the internal fluid. In fact, this analytical expression and (\ref{lengthevol}) are equivalent provided that \(L(t)>0\) for all time \(t\). This fact follows from the proposition below, whose proof can be garnered from~\cite{gancedo2019global}.
\begin{proposition} \label{p2}
    Let \(V_{0}=\pi R^{2}\) be the initial volume of the internal fluid. For any \(t \geq 0\) such that \(L(t)>0\),
    \begin{align}
        \biggl(\frac{L(t)}{2\pi}\biggr)^{2}&=R^{2}\biggl(1+\frac{1}{2\pi}\mbox{Im}\int_{-\pi}^{\pi}\int_{0}^{\alpha}e^{i(\alpha-\eta)}\sum_{n\geq 1} \frac{i^{n}}{n!}(\theta(\alpha)-\theta(\eta))^{n}d\eta d\alpha\biggr)^{-1} \label{l}
    \end{align}
    implies
    \begin{align*}
        L_{t}(t)&=-\int_{-\pi}^{\pi}(1+\theta_{\alpha}(\alpha))U(\alpha)d\alpha.
    \end{align*}
\end{proposition}
\begin{remark}
    That \(V_{0}=\pi R^{2}\) is not to say that the internal fluid is initially a circle of radius \(R\).
\end{remark}
To derive (\ref{l}) from the incompressibility condition on the internal fluid, let \(\mathcal{D}\) be the region enclosed by the fluid boundary \(\Gamma\). Then the volume of the region \(\mathcal{D}\) is given by
\begin{align} \label{v1}
    V&=\int_{\mathcal{D}}dx\wedge dy = \frac{1}{2}\int_{-\pi}^{\pi}(-z_{2}(\alpha),z_{1}(\alpha))\cdot z_{\alpha}(\alpha)d\alpha,
\end{align}
where \(\wedge\) in (\ref{v1}) is the wedge product of differential forms. In complex variable notation,
\begin{align*}
        V=\frac{1}{2}\int_{-\pi}^{\pi}\mbox{Im}\biggl(\overline{z(\alpha)}z_{\alpha}(\alpha)\biggr)d\alpha=\frac{1}{2}\mbox{Im}\int_{-\pi}^{\pi}\overline{z(\alpha)}z_{\alpha}(\alpha)d\alpha.
\end{align*}
Using that
\begin{align*}
    z_{\alpha}(\alpha)&=\frac{L(t)}{2\pi}e^{i(\alpha+\theta(\alpha))} \\
    z(\alpha)&=z(0)+\int_{0}^{\alpha}z_{\eta}(\eta)d\eta,
\end{align*}
we can write
\begin{align}
    V&=\pi\biggl(\frac{L(t)}{2\pi}\biggr)^{2}\biggl(1+\frac{1}{2\pi}\mbox{Im}\int_{-\pi}^{\pi}\int_{0}^{\alpha}e^{i(\alpha-\eta)}\sum_{n=1}^{\infty}\frac{i^{n}}{n!}(\theta(\alpha)-\theta(\eta))^{n}d\eta d\alpha\biggr). \nonumber
    \end{align}
Since the internal fluid is incompressible, \(V_{0}=\pi R^{2}=V\), which implies
\begin{align*}
    \biggl(\frac{L(t)}{2\pi}\biggr)^{2} = R^{2}\biggl(1+\frac{1}{2\pi}\mbox{Im}\int_{-\pi}^{\pi}\int_{0}^{\alpha}e^{i(\alpha-\eta)}\sum_{n=1}^{\infty}\frac{i^{n}}{n!}(\theta(\alpha)-\theta(\eta))^{n}d\eta d\alpha\biggr)^{-1}.
\end{align*}
\begin{remark}
    The interior fluid's incompressibility combined with the isoperimetric inequality ensures that \(L(t)>0\) is satisfied for all \(t \geq 0\).
\end{remark}

\subsection{The Circular Interface under HLS Parametrization}
\label{steadystatesolutions}
The proposition below characterizes the circular interface under HLS parametrization.
\begin{proposition} \label{circularinterface}
    Let \(R>0\). The interface at time \(t\) is a circle of radius \(R\) if and only if
    \begin{align*}
        (\theta(\alpha,t),L(t))=(\hat{\theta}(0,t),2\pi R),
    \end{align*}
    where the parametrization is HLS.
\end{proposition}
\begin{proof}
First, we check that \((\theta(\alpha,t),L(t))=(\hat{\theta}(0,t),2\pi R)\) is a circle of radius \(R\) for a fixed \(t\). It suffices to show that the curve has a constant curvature \(\abs{\frac{d^{2}z}{ds^{2}}}\) of \(1/R\). Observe that
\begin{align*}
    \frac{d^{2}z}{ds^{2}} = \frac{d}{d\alpha}\biggl(\frac{dz}{d\alpha}\cdot\frac{d\alpha}{ds}\biggr)\frac{d\alpha}{ds}=\frac{d^{2}z}{d\alpha^{2}}\cdot\abs{z_{\alpha}(\alpha,t)}^{-2}=\frac{ie^{i(\alpha+\hat{\theta}(0,t))}}{R}.
\end{align*}
Since \(\hat{\theta}(0,t)\) is a real number,
\begin{align*}
    \abs{\frac{d^{2}z}{ds^{2}}}=\frac{1}{R},
\end{align*}
as needed. To prove the converse, suppose that the interface at time \(t\) is a circle of radius \(R\). Then \(L(t)=2\pi R\). That \(\abs{\frac{d^{2}z}{ds^{2}}}=\frac{1}{R}\) implies that \(\abs{1+\theta_{\alpha}(\alpha,t)}=1\). Due to the periodicity of \(\theta\), we have \(\theta_{\alpha}(\alpha, t)=0\), i.e., \(\theta(\alpha,t)\) depends only on time \(t\). Then \(\hat{\theta}(0,t) = \theta(\alpha,t)\), as needed.
\end{proof}

\section{Statement of the Main Theorem}
\label{maintheorem}
To study the simple two-dimensional model given by (\ref{1.1}) through (\ref{1.4}), we have adopted the single-layer potential form (\ref{slpvelo}) for the fluid velocity. To completely describe the dynamics of the fluid velocity, it is therefore sufficient to study the dynamics of the interface itself. To that end, we have taken the HLS parametrization of the interface to obtain a pair of dynamics equations, (\ref{lengthevol}) and (\ref{thetaevol}), for the interface. We have then reformulated the dynamics equation (\ref{lengthevol}) for the length of the interface into (\ref{l}).

The main theorem of our work is that the equations (\ref{l}) and (\ref{thetaevol}) for the dynamics of the interface have a unique solution that is global in time, provided that the initial datum is sufficiently small as measured by the norm of \(\dot{\mathcal{F}}^{1,1}\). The unique solution also decays exponentially in time in the norm of \(\dot{\mathcal{F}}_{\nu}^{1,1}\), where \(\nu\) is given in (\ref{nu}) and \(\nu_{0}>0\) is dependent on the initial datum. In view of Proposition \ref{circularinterface}, this implies that the initial interface decays exponentially to a circular shape.
\begin{theorem} \label{mainthm}
    Fix \(\gamma>0\). If the initial datum \(\theta^{0} \in \dot{\mathcal{F}}^{1,1}\) such that \(\abs{\mathcal{F}(\theta^{0})(0)}\) and \(\norm{\theta^{0}}_{\dot{\mathcal{F}}^{1,1}}\) are sufficiently small, then for any \(T \in (0,\infty)\) there exists a unique solution
    \begin{align*}
        \theta(\alpha,t) \in C([0,T];\dot{\mathcal{F}}_{\nu}^{1,1}) \cap L^{1}([0,T];\dot{\mathcal{F}}_{\nu}^{2,1})
    \end{align*}
    to the equations (\ref{l}) and (\ref{thetaevol}), where \(\nu\) is given in (\ref{nu}) and \(\nu_{0}>0\) is dependent on \(\theta^{0}\). The solution becomes instantaneously analytic. In particular, for any \(t \in [0,T]\)
    \begin{align*}
    \normx{\theta(t)}{1} + \biggl(\Lambda(\norm{\theta^{0}}_{\dot{\mathcal{F}}^{1,1}})-\nu_{0}\biggr)\int_{0}^{t}\normx{\theta(\tau)}{2}d\tau & \leq \norm{\theta^{0}}_{\dot{\mathcal{F}}^{1,1}},
    \end{align*}
    where \(\Lambda(\norm{\theta^{0}}_{\dot{\mathcal{F}}^{1,1}})\) is given in (\ref{lambdaexpression}). Moreover, \(\normx{\theta(t)}{1}\) decays exponentially in time.
\end{theorem}
\begin{remark}
    The assumption that the initial datum be ``sufficiently small" can be made explicit in the sense that for any \(\gamma>0\), one can explicitly derive an upper bound on the magnitudes of \(\abs{\mathcal{F}(\theta^{0})(0)}\) and \(\norm{\theta^{0}}_{\dot{\mathcal{F}}^{1,1}}\).
\end{remark}

\section{The Interfacial Fluid Velocity}
\subsection{Formulation in Complex Variable Notation}
We set out to rewrite (\ref{slpvelo}) in complex variable notation. The signed curvature \(\kappa\) that appears in the single-layer potential is defined by, in vector notation,
\begin{equation} \label{curvature_def}
    \bm{\tau}'(s) = -\kappa(s)\bm{n}(s), \nonumber
\end{equation}
where \(s\) denotes arclength parametrization. Letting \(\bm{\tau}=(\tau_{1},\tau_{2})\) and \(z=(z_{1},z_{2})\), we use the Jacobian between arclength parametrization and HLS parametrization to obtain
\begin{align*}
    \tau_{i}'(s) = \frac{d\tau_{i}}{ds}=\frac{d}{ds}\biggl(\frac{dz_{i}}{ds}\biggr) = \frac{d}{d\beta}\biggl(\frac{dz_{i}}{d\beta}\frac{d\beta}{ds}\biggr) \frac{d\beta}{ds} = \frac{d^{2}z_{i}}{d\beta^{2}}\abs{z_{\beta}(\beta,t)}^{-2},
\end{align*}
which yields, in vector notation,
\begin{align*}
    u_{j}(\bm{x})&=\frac{2\pi}{L(t)}\frac{\gamma}{4\pi}\int_{-\pi}^{\pi}z''(\beta)\cdot G_{\cdot j}(\bm{x} - z(\beta))d\beta,
\end{align*}
where
\begin{equation*}
    G_{\cdot j}(\bm{x} - z(\beta)) = (G_{1j}(\bm{x} - z(\beta)), G_{2j}(\bm{x} - z(\beta))).
\end{equation*}
Let \(\bm{x}=z(\alpha) \in \Gamma\). To rewrite the current expression for \(u_{j}(\bm{x}) = u_{j}(z(\alpha))\) in complex variable notation, we note that in complex variable notation,
\begin{align*}
    G_{\cdot j}(z(\alpha) - z(\beta)) &= G_{1j}(z(\alpha) - z(\beta)) +i G_{2j}(z(\alpha) - z(\beta)), \\
    z'(\beta) &= \frac{L(t)}{2\pi}e^{i(\beta+\theta(\beta))}.
\end{align*}
We apply integration by parts to obtain
\begin{align*}
    u_{j}(z(\alpha)) =& -\frac{2\pi}{L(t)}\frac{\gamma}{4\pi}\int_{-\pi}^{\pi}\mbox{Re}\biggl(\overline{z'(\beta)}\frac{d}{d\beta}\biggl(G_{\cdot j}(z(\alpha) - z(\beta))\biggr)\biggr)d\beta,
\end{align*}
where \(\overline{z'(\beta)}\) denotes the complex conjugate of \(z'(\beta)\) and
\begin{align*}
    \mbox{Re}\biggl(\overline{z'(\beta)}\frac{d}{d\beta}\biggl(G_{\cdot j}(z(\alpha)-z(\beta))\biggr)\biggr) = &\frac{L(t)}{2\pi}\biggl(\cos(\beta + \theta(\beta))\frac{d}{d\beta}\biggl(G_{1j}(z(\alpha) - z(\beta))\biggr) \\
    &+ \sin(\beta + \theta(\beta))\frac{d}{d\beta}\biggl(G_{2j}(z(\alpha) - z(\beta))\biggr)\biggr).
\end{align*}
Hence,
\begin{align*}
    u_{j}(z(\alpha)) = &-\frac{\gamma}{4\pi}\int_{-\pi}^{\pi}\cos(\beta + \theta(\beta))\frac{d}{d\beta}\biggl(G_{1j}(z(\alpha) - z(\beta))\biggr) \\
    &+ \sin(\beta + \theta(\beta))\frac{d}{d\beta}\biggl(G_{2j}(z(\alpha) - z(\beta))\biggr)d\beta.
\end{align*}
By changing the variable of integration from \(\beta\) to \(\beta'=\alpha-\beta\) and rewriting the sine and cosine in complex variable notation, we obtain
\begin{align}
    &u_{j}(z(\alpha)) \nonumber \\
    =& \frac{\gamma}{4\pi}\int_{-\pi}^{\pi}\frac{1}{2}\biggl(e^{i(\alpha-\beta + \theta(\alpha - \beta))}+e^{-i(\alpha - \beta + \theta(\alpha - \beta))}\biggr) \label{vjterm1} \\
    &\cdot \frac{d}{d\beta}\biggl(G_{1j}(z(\alpha) - z(\alpha - \beta))\biggr) \nonumber \\
    &+ \frac{1}{2i}\biggl(e^{i(\alpha-\beta+\theta(\alpha-\beta))}-e^{-i(\alpha-\beta+\theta(\alpha-\beta))}\biggr)\frac{d}{d\beta}\biggl(G_{2j}(z(\alpha) - z(\alpha - \beta))\biggr)d\beta. \label{vjterm2}
\end{align}
\subsection{The Normal Speed \(U\)}
To obtain the normal speed in complex variable notation, we take the dot product of (\ref{velo}) and \(-\bm{n}\) to get
\begin{align*}
    U = \bm{u}\cdot(-\bm{n}),
\end{align*}
which can be rewritten in complex variable notation as
\begin{equation*}
    U(\alpha) = \mbox{Re}\biggl((u_{1}(\alpha) - iu_{2}(\alpha))ie^{i(\alpha + \theta(\alpha))}\biggr).
\end{equation*}
To obtain an analytical expression for \(U(\alpha)\) in complex variable notation, we first simplify (\ref{vjterm1}) and (\ref{vjterm2}). We note that
\begin{align*}
    &G_{11}(z(\alpha) - z(\alpha - \beta)) = -\log\abs{z(\alpha) - z(\alpha-\beta)} + \frac{(z_{1}(\alpha) - z_{1}(\alpha-\beta))^{2}}{\abs{z(\alpha) - z(\alpha-\beta)}^{2}}, \\
    &G_{12}(z(\alpha) - z(\alpha-\beta)) = \frac{(z_{1}(\alpha) - z_{1}(\alpha-\beta))(z_{2}(\alpha) - z_{2}(\alpha - \beta))}{\abs{z(\alpha) - z(\alpha-\beta)}^{2}}, \\
    &G_{21}(z(\alpha) - z(\alpha-\beta)) = \frac{(z_{1}(\alpha) - z_{1}(\alpha-\beta))(z_{2}(\alpha) - z_{2}(\alpha - \beta))}{\abs{z(\alpha) - z(\alpha-\beta)}^{2}}, \\
    &G_{22}(z(\alpha) - z(\alpha-\beta)) = -\log\abs{z(\alpha) - z(\alpha - \beta)} + \frac{(z_{2}(\alpha) - z_{2}(\alpha - \beta))^{2}}{\abs{z(\alpha) - z(\alpha - \beta)}^{2}}.
\end{align*}
Letting
\begin{equation*}
    w(\alpha, \beta) = \int_{0}^{1}e^{i(\alpha + (s-1)\beta + \theta(\alpha + (s-1)\beta))}ds,
\end{equation*}
we can write
\begin{align*}
    z(\alpha) - z(\alpha - \beta) = \beta\int_{0}^{1}z_{\alpha}(\alpha + (s-1)\beta)ds = \frac{\beta L(t)}{2\pi}w(\alpha,\beta).
\end{align*}
Denoting the complex conjugate of \(w\) by \(\overline{w}\), we then obtain
\begin{align*}
    &\frac{\partial}{\partial\beta}\biggl(-\log\abs{z(\alpha) - z(\alpha-\beta)}\biggr) = -\frac{1}{\beta}-\frac{w_{\beta}}{2w}-\frac{\overline{w}_{\beta}}{2\overline{w}}.
\end{align*}
Moreover,
\begin{align*}
    &\frac{\partial}{\partial\beta}\biggl(\frac{(z_{1}(\alpha) - z_{1}(\alpha - \beta))^{2}}{\abs{z(\alpha) - z(\alpha - \beta)}^{2}}\biggr) \\
    =&\frac{1}{2}\cdot\biggl(\frac{1}{\overline{w}}+\frac{1}{w}\biggr)(w_{\beta}+\overline{w}_{\beta})-\frac{1}{4}\cdot\biggl(\frac{1}{\overline{w}} + \frac{1}{w}\biggr)^{2}(w_{\beta}\overline{w} + w\overline{w}_{\beta}).
\end{align*}
Similarly,
\begin{align*}
    &\frac{\partial}{\partial\beta}\biggl(\frac{(z_{2}(\alpha) - z_{2}(\alpha - \beta))^{2}}{\abs{z(\alpha) - z(\alpha-\beta)}^{2}}\biggr) \\
    =&-\frac{1}{2}\biggl(\frac{1}{\overline{w}}-\frac{1}{w}\biggr)(w_{\beta} - \overline{w}_{\beta})+\frac{1}{4}\biggl(\frac{1}{\overline{w}}-\frac{1}{w}\biggr)^{2}(w_{\beta}\overline{w}+w\overline{w}_{\beta}).
\end{align*}
Lastly,
\begin{align*}
    &\frac{\partial}{\partial\beta}\biggl(\frac{(z_{1}(\alpha) - z_{1}(\alpha - \beta))(z_{2}(\alpha) - z_{2}(\alpha - \beta))}{\abs{z(\alpha) - z(\alpha-\beta)}^{2}}\biggr) \\
    =&\frac{1}{2i}\biggl(\frac{w_{\beta}}{\overline{w}}-\frac{\overline{w}_{\beta}}{w}\biggr)-\frac{1}{4i}\biggl(\frac{w}{\overline{w}}-\frac{\overline{w}}{w}\biggr)\biggl(\frac{w_{\beta}}{w}+\frac{\overline{w}_{\beta}}{\overline{w}}\biggr).
\end{align*}
Hence,
\begin{align*}
    &\frac{\partial}{\partial\beta}\biggl(G_{11}(z(\alpha) - z(\alpha - \beta))\biggr) \\
    =&-\frac{1}{\beta}-\frac{w_{\beta}}{2w}-\frac{\overline{w}_{\beta}}{2\overline{w}}+\frac{1}{2}\biggl(\frac{1}{\overline{w}}+\frac{1}{w}\biggr)(w_{\beta}+\overline{w}_{\beta})-\frac{1}{4}\biggl(\frac{1}{\overline{w}}+\frac{1}{w}\biggr)^{2}(w_{\beta}\overline{w}+w\overline{w}_{\beta}),
\end{align*}
\begin{align*}
    &\frac{\partial}{\partial\beta}\biggl(G_{22}(z(\alpha) - z(\alpha - \beta))\biggr) \\
    =&-\frac{1}{\beta}-\frac{w_{\beta}}{2w}-\frac{\overline{w}_{\beta}}{2\overline{w}}-\frac{1}{2}\biggl(\frac{1}{\overline{w}}-\frac{1}{w}\biggr)(w_{\beta}-\overline{w}_{\beta})+\frac{1}{4}\biggl(\frac{1}{\overline{w}}-\frac{1}{w}\biggr)^{2}(w_{\beta}\overline{w}+w\overline{w}_{\beta}),
\end{align*}
and
\begin{align*}
    &\frac{\partial}{\partial\beta}\biggl(G_{12}(z(\alpha) - z(\alpha - \beta))\biggr)=\frac{\partial}{\partial\beta}\biggl(G_{21}(z(\alpha) - z(\alpha - \beta))\biggr) \\
    =&\frac{1}{2i}\biggl(\frac{w_{\beta}}{\overline{w}}-\frac{\overline{w}_{\beta}}{w}\biggr)-\frac{1}{4i}\biggl(\frac{w}{\overline{w}}-\frac{\overline{w}}{w}\biggr)\biggl(\frac{w_{\beta}}{w}+\frac{\overline{w}_{\beta}}{\overline{w}}\biggr).
\end{align*}
For notational convenience, let us write
\begin{align*}
    &w = C_{1} + L_{1} + N_{1}, \\
    &w^{-1} = C_{2} + L_{2} + N_{2}, \\
    &w_{\beta} = C_{\beta} + L_{\beta} + N_{\beta},
\end{align*}
where \(C_{1}\), \(L_{1}\), and \(N_{1}\) are the parts of \(w\) which are constant, linear, and superlinear in the variable \(\phi = \theta - \hat{\theta}(0)\), respectively; \(C_{2}\), \(L_{2}\), and \(N_{2}\) are the parts of \(w^{-1}\) which are constant, linear, and superlinear in \(\phi\), respectively; lastly, \(C_{\beta}\), \(L_{\beta}\), and \(N_{\beta}\) are the parts of \(w_{\beta}\) which are constant, linear, and superlinear in \(\phi\). Similarly, let us write
\begin{align*}
    &\frac{\partial}{\partial\beta}\biggl(G_{11}(z(\alpha) - z(\alpha - \beta))\biggr) = C_{11} + L_{11} + N_{11}, \\
    &\frac{\partial}{\partial\beta}\biggl(G_{12}(z(\alpha) - z(\alpha - \beta))\biggr) = \frac{\partial}{\partial\beta}\biggl(G_{21}(z(\alpha) - z(\alpha - \beta))\biggr) = C_{12} + L_{12} + N_{12}, \\
    &\frac{\partial}{\partial\beta}\biggl(G_{22}(z(\alpha) - z(\alpha - \beta))\biggr) = C_{22} + L_{22} + N_{22},
\end{align*}
where \(C_{11}\), \(L_{11}\), and \(N_{11}\) are the parts of \(\frac{\partial}{\partial\beta}\biggl(G_{11}(z(\alpha) - z(\alpha - \beta))\biggr)\) which are constant, linear, and superlinear in \(\phi\); \(C_{12}\), \(L_{12}\), and \(N_{12}\) are the parts of \(\frac{\partial}{\partial\beta}\biggl(G_{12}(z(\alpha) - z(\alpha - \beta))\biggr)\) which are constant, linear, and superlinear in \(\phi\); lastly, \(C_{22}\), \(L_{22}\), and \(N_{22}\) are the parts of \(\frac{\partial}{\partial\beta}\biggl(G_{22}(z(\alpha) - z(\alpha - \beta))\biggr)\) which are constant, linear, and superlinear in \(\phi\). Using these expressions, we can write
\begin{align}
    &u_{1}(z(\alpha))-iu_{2}(z(\alpha)) \label{conjugateu} \\
    =&\frac{\gamma}{4\pi}\int_{-\pi}^{\pi}\frac{e^{i\hat{\theta}(0)}e^{i(\alpha-\beta)}}{2}\biggl(1+i\phi(\alpha-\beta)+\sum_{n=2}^{\infty}\frac{(i\phi(\alpha-\beta))^{n}}{n!}\biggr) \nonumber \\
    &\cdot\biggl((C_{11}+L_{11}+N_{11})-(C_{22}+L_{22}+N_{22})-2i(C_{12}+L_{12}+N_{12})\biggr) \nonumber \\
    &+\frac{e^{-i\hat{\theta}(0)}e^{-i(\alpha-\beta)}}{2}\biggl(1-i\phi(\alpha-\beta)+\sum_{n=2}^{\infty}\frac{(-i\phi(\alpha-\beta))^{n}}{n!}\biggr) \nonumber \\
    &\cdot\biggl((C_{11}+L_{11}+N_{11})+(C_{22}+L_{22}+N_{22})\biggr)d\beta. \nonumber
\end{align}
Let
\begin{equation*}
    u_{1}(\alpha) - iu_{2}(\alpha) = \mathfrak{C}(\alpha) + \mathfrak{L}(\alpha) + \mathfrak{N}(\alpha),
\end{equation*}
where \(\mathfrak{C}\), \(\mathfrak{L}\), and \(\mathfrak{N}\) are the parts of \(u_{1}-iu_{2}\) which are constant, linear, and superlinear in \(\phi\), respectively. In particular, \(\mathfrak{C}(\alpha)=0\). Let \(U = U_{0} + U_{1} + U_{\geq 2}\), where \(U_{0}\), \(U_{1}\), and \(U_{\geq 2}\) are the parts of \(U\) which are constant, linear, and superlinear in \(\phi\), respectively. Then
\begin{equation} \label{U0}
    U_{0}(\alpha) = \mbox{Re}\biggl(ie^{i\alpha}e^{i\hat{\theta}(0)}\mathfrak{C}(\alpha)\biggr) = 0.
\end{equation}
To find expressions for \(U_{1}\) and \(U_{\geq 2}\), we first rewrite
\begin{align*}
    U(\alpha)&=\mbox{Re}\biggl(ie^{i\alpha}e^{i\hat{\theta}(0)}\mathfrak{L}(\alpha)\biggr)+\mbox{Re}\biggl(ie^{i\alpha}e^{i\hat{\theta}(0)}\biggl(\mathfrak{L}(\alpha)(e^{i\phi(\alpha)}-1)+\mathfrak{N}(\alpha)e^{i\phi(\alpha)}\biggr)\biggr).
\end{align*}
Then it is clear that
\begin{align}
    U_{1}(\alpha)&=\mbox{Re}\biggl(ie^{i\alpha}e^{i\hat{\theta}(0)}\mathfrak{L}(\alpha)\biggr), \label{U1} \\
    U_{\geq 2}(\alpha)&=\mbox{Re}\biggl(ie^{i\alpha}e^{i\hat{\theta}(0)}\biggl(\mathfrak{L}(\alpha)(e^{i\phi(\alpha)}-1) +\mathfrak{N}(\alpha)e^{i\phi(\alpha)}\biggr)\biggr). \label{U2}
\end{align}

\subsection{The Tangential Speed \(T\)}
\label{tangentialspeed}
Let us rewrite the right hand side of (\ref{thetaevol}) as
\begin{align*}
    \frac{2\pi}{L(t)}\biggl(U_{\alpha}(\alpha)+T(\alpha)(1+\phi_{\alpha}(\alpha))\biggr) = \mathcal{C}(\alpha)+\mathcal{L}(\alpha)+\mathcal{N}(\alpha),
\end{align*}
where \(\mathcal{C}\), \(\mathcal{L}\), and \(\mathcal{N}\) are the parts of the right hand side of the evolution equation for \(\theta\) which are constant, linear, and superlinear in the variable \(\phi=\theta-\hat{\theta}(0)\), respectively. We will completely determine the frame of parametrization by specifying the analytical expression for \(T(\alpha)\) such that \(\mathcal{C}=0\). To begin, let us rewrite the right hand side of (\ref{T}) as
\begin{align*}
    &\int_{0}^{\alpha}(1 + \phi_{\alpha}(\eta))U(\eta)d\eta - \frac{\alpha}{2\pi}\int_{-\pi}^{\pi}(1 + \phi_{\alpha}(\eta))U(\eta)d\eta + T(0) \\
    = &T_{0}(\alpha)+T_{1}(\alpha)+T_{\geq2}(\alpha),
\end{align*}
where \(T_{0}\), \(T_{1}\), and \(T_{\geq2}\) are the parts of \(T\) which are constant, linear, and superlinear in \(\phi\), respectively. We note that
\begin{align}
    T_{0}(\alpha) =& \int_{0}^{\alpha}U_{0}(\eta)d\eta - \frac{\alpha}{2\pi}\int_{-\pi}^{\pi}U_{0}(\eta)d\eta + T(0) = T(0), \label{T0} \\
    T_{1}(\alpha) =& \int_{0}^{\alpha}U_{1}(\eta)d\eta + \int_{0}^{\alpha}\phi_{\alpha}(\eta)U_{0}(\eta)d\eta \nonumber \\
    &- \frac{\alpha}{2\pi}\int_{-\pi}^{\pi}U_{1}(\eta)d\eta - \frac{\alpha}{2\pi}\int_{-\pi}^{\pi}\phi_{\alpha}(\eta)U_{0}(\eta)d\eta, \nonumber \\
    T_{\geq 2}(\alpha) =& \int_{0}^{\alpha}U_{\geq 2}(\eta)d\eta + \int_{0}^{\alpha}\phi_{\alpha}(\eta)U_{\geq 1}(\eta)d\eta \nonumber \\
    &- \frac{\alpha}{2\pi}\int_{-\pi}^{\pi}U_{\geq 2}(\eta)d\eta - \frac{\alpha}{2\pi}\int_{-\pi}^{\pi}\phi_{\alpha}(\eta)U_{\geq 1}(\eta)d\eta, \nonumber
\end{align}
where we define \(U_{\geq1} = U_{1} + U_{\geq2}\). Let \(T(0)=0\). Then using (\ref{U0}), we obtain
\begin{align*}
    \mathcal{C}(\alpha) = \frac{2\pi}{L(t)}\biggl((U_{0})_{\alpha}(\alpha)+T_{0}(\alpha)\biggr) = \frac{2\pi}{L(t)}T_{0}(\alpha) = \frac{2\pi}{L(t)}T(0) = 0.
\end{align*}

\section{Steady-State Solutions}
\label{steadystatesolutions2}
We know from Proposition \ref{circularinterface} that \((\theta(\alpha,t),L(t))=(\hat{\theta}(0,t),2\pi R)\) corresponds to a circle of radius \(R\). The circular interface becomes a solution to (\ref{lengthevol}) and (\ref{thetaevol}) if \(\hat{\theta}(0,t)\) is constant in time. In this case, the interface is stationary. The following proposition summarizes the existence of steady-state solutions to (\ref{lengthevol}) and (\ref{thetaevol}).
\begin{proposition}
    For any constant \(c\), the circle defined by
    \begin{align*}
        (\theta(\alpha,t),L(t)) = (c,2\pi R)
    \end{align*}
    is a time-independent solution of (\ref{lengthevol}) and (\ref{thetaevol}) in which \(T(\alpha,t)\) is given by (\ref{T}) and \(U(\alpha,t)\) is given by
    \begin{align*}
        U(\alpha,t)&=\mbox{Re}\biggl((u_{1}(\alpha,t)-iu_{2}(\alpha,t))ie^{i(\alpha+\theta(\alpha,t))}\biggr)
    \end{align*}
    with \(u_{1}(\alpha,t)-iu_{2}(\alpha,t)\) given by (\ref{conjugateu}).
\end{proposition}
\begin{proof}
    Let \((\theta(\alpha,t),L(t)) = (\hat{\theta}(0,t),2\pi R)\) be a circle of radius \(R\) such that \(\hat{\theta}(0,t) = c\) for some constant \(c\). Since \(U(\alpha,t)=T(\alpha,t)=0\), the right hand sides of (\ref{lengthevol}) and (\ref{thetaevol}) vanish. Since \((\theta_{t}(\alpha,t),L_{t}(t))=(0,0)\), (\ref{lengthevol}) and (\ref{thetaevol}) are indeed satisfied by \((\theta(\alpha,t),L(t))=(c,2\pi R)\).
\end{proof}

\section{The Principal Linear Operator for the \texorpdfstring{\(\theta\)}{t} Equation}
\label{principallinearoperatorforthetaequation}

Let us explicitly calculate the operator \(\mathcal{L}\). We note that
\begin{align*}
    \mathcal{L}(\alpha) = \frac{2\pi}{L(t)}\biggl((U_{1})_{\alpha}(\alpha)+T_{0}(\alpha)\phi_{\alpha}(\alpha)+T_{1}(\alpha)\biggr) = \frac{2\pi}{L(t)}\biggl((U_{1})_{\alpha}(\alpha)+T_{1}(\alpha)\biggr)
\end{align*}
by (\ref{T0}). By (\ref{U0}),
\begin{align*}
    T_{1}(\alpha) =& \int_{0}^{\alpha}U_{1}(\eta)d\eta + \int_{0}^{\alpha}\phi_{\alpha}(\eta)U_{0}(\eta)d\eta \\
    &- \frac{\alpha}{2\pi}\int_{-\pi}^{\pi}U_{1}(\eta)d\eta - \frac{\alpha}{2\pi}\int_{-\pi}^{\pi}\phi_{\alpha}(\eta)U_{0}(\eta)d\eta \\
    =&\int_{0}^{\alpha}U_{1}(\eta)d\eta-\frac{\alpha}{2\pi}\int_{-\pi}^{\pi}U_{1}(\eta)d\eta.
\end{align*}
Using (\ref{m}), we can write
\begin{align} \label{operatorL}
    \mathcal{L}(\alpha)=\frac{2\pi}{L(t)}\biggl((U_{1})_{\alpha}(\alpha) + \mathcal{M}(U_{1})(\alpha)\biggr).
\end{align}

\subsection{The Fourier Modes of \(\mathcal{L}\)}
In this Section, we verify that \(\mathcal{L}\) is the Hilbert transform acting on the spatial derivative of \(\theta\), up to the \(\pm 1\) Fourier modes, by checking that its Fourier multiplier is \(\abs{k}\) for \(\abs{k} > 1\). To that end, we will calculate \(\mathcal{F}(\mathcal{L})(k)\), the \(k\)th Fourier mode of \(\mathcal{L}(\alpha)\), for all \(k\in \mathbb{Z} \setminus \{0\}\). Using (\ref{operatorMFourier}), we obtain that for \(k \neq 0\),
\begin{align} \label{fourierlinear}
    \mathcal{F}(\mathcal{L})(k)=\frac{2\pi}{L(t)}\biggl(\mathcal{F}((U_{1})_{\alpha})(k)-\frac{i}{k}\mathcal{F}(U_{1})(k)\biggr).
\end{align}
First, we set out to find the expressions for \(U_{1}\) and \((U_{1})_{\alpha}\).
From (\ref{U1}), we obtain
\begin{align}
    U_{1}(\alpha)&=\mbox{Re}\biggl(ie^{i\alpha}e^{i\hat{\theta}(0)}\biggr)\mbox{Re}\mathfrak{L}(\alpha)-\mbox{Im}\biggl(ie^{i\alpha}e^{i\hat{\theta}(0)}\biggr)\mbox{Im}\mathfrak{L}(\alpha). \nonumber
\end{align}
We note that
\begin{align} 
    &U_{1}(\alpha)=\frac{\gamma}{4\pi}\int_{-\pi}^{\pi}\biggl( \label{u1exp} \\
    &\int_{0}^{1}e^{-i\beta s}\phi(\alpha+\beta(-1+s))ds\cdot\frac{-(-i+(i+\beta)e^{i\beta})(-1-2e^{i\beta}+e^{2i\beta})}{4(-1+e^{i\beta})^{2}} \nonumber \\
    +&\int_{0}^{1}e^{i\beta s}\phi(\alpha+\beta(-1+s))ds\cdot\frac{e^{-i\beta}(-1+2e^{i\beta}+e^{2i\beta})(\beta+i(-1+e^{i\beta}))}{4(-1+e^{i\beta})^{2}} \nonumber \\
    +&\int_{0}^{1}e^{-i\beta s}(-1+s)\phi(\alpha+\beta(-1+s))ds\cdot\frac{-(-1+e^{i\beta})\beta(-1-2e^{i\beta}+e^{2i\beta})}{4(-1+e^{i\beta})^{2}} \nonumber \\
    +&\int_{0}^{1}e^{i\beta s}(-1+s)\phi(\alpha+\beta(-1+s))ds\cdot\frac{-e^{-i\beta}(-1+e^{i\beta})\beta(-1+2e^{i\beta}+e^{2i\beta})}{4(-1+e^{i\beta})^{2}} \nonumber \\
    +&\int_{0}^{1}e^{-i\beta s}(-1+s)\phi'(\alpha+\beta(-1+s))ds\cdot\frac{-(-1+e^{i\beta})i\beta(-1-2e^{i\beta}+e^{2i\beta})}{4(-1+e^{i\beta})^{2}} \nonumber \\
    +&\int_{0}^{1}e^{i\beta s}(-1+s)\phi'(\alpha+\beta(-1+s))ds\cdot\frac{e^{-i\beta}(-1+e^{i\beta})i\beta(-1+2e^{i\beta}+e^{2i\beta})}{4(-1+e^{i\beta})^{2}} \nonumber \\
    +&\phi(\alpha-\beta)\cdot\frac{e^{-i\beta}(-1+e^{i\beta})(-i)(1+e^{i\beta}+e^{2i\beta}+e^{3i\beta})}{4(-1+e^{i\beta})^{2}}\biggr)d\beta. \nonumber
\end{align}
We obtain the expression for \((U_{1})_{\alpha}\) by differentiating (\ref{u1exp}) with respect to \(\alpha\). Now, taking the Fourier modes of \(U_{1}\) and \((U_{1})_{\alpha}\) and plugging them into (\ref{fourierlinear}), we obtain that for \(k \notin \{0, \pm 1\}\),
\begin{align*}
    &\mathcal{F}(\mathcal{L})(k)=\frac{2\pi}{L(t)}\frac{\gamma}{4\pi}\mathcal{F}(\phi)(k)\biggl( \\
    &\int_{-\pi}^{\pi}\frac{(i-(i+\beta)e^{i\beta})(-1-2e^{i\beta}+e^{2i\beta})}{4(-1+e^{i\beta})^{2}}\frac{e^{-i\beta}-e^{-i\beta k}}{\beta}d\beta\biggl(\frac{k}{k-1}+\frac{1}{k(1-k)}\biggr) \\
    +&\int_{-\pi}^{\pi}\frac{e^{-i\beta}(-1+2e^{i\beta}+e^{2i\beta})(\beta+i(-1+e^{i\beta}))}{4(-1+e^{i\beta})^{2}}\frac{e^{i\beta}-e^{-i\beta k}}{\beta}d\beta\biggl(\frac{k}{1+k}-\frac{1}{k(1+k)}\biggr) \\
    +&\int_{-\pi}^{\pi}\frac{-i(-1-2e^{i\beta}+e^{2i\beta})}{4(-1+e^{i\beta})}\frac{e^{-i\beta(1+k)}(e^{i\beta k} - e^{i\beta})}{\beta}d\beta\biggl(\frac{-k^{2}}{(-1+k)^{2}}+\frac{1}{(-1+k)^{2}}\biggr) \\
    +&\int_{-\pi}^{\pi}\frac{-i(-1-2e^{i\beta}+e^{2i\beta})e^{-i\beta k}}{4(-1+e^{i\beta})}d\beta\biggl(\frac{ik^{2}}{-1+k}-\frac{i}{-1+k}\biggr) \\
    +&\int_{-\pi}^{\pi}\frac{e^{-i\beta}i(-1+2e^{i\beta}+e^{2i\beta})}{4(-1+e^{i\beta})}\frac{e^{-i\beta k}(-1+e^{i\beta(1+k)})}{\beta}d\beta\biggl(\frac{-k^{2}}{(1+k)^{2}}+\frac{1}{(1+k)^{2}}\biggr)
\end{align*}
\begin{align*}
    +&\int_{-\pi}^{\pi}\frac{e^{-i\beta}(-1+2e^{i\beta}+e^{2i\beta})e^{-i\beta k}}{4(-1+e^{i\beta})}d\beta\biggl(\frac{-k^{2}}{1+k}+\frac{1}{1+k}\biggr) \\
    +&\int_{-\pi}^{\pi}\frac{1+2e^{i\beta}-e^{2i\beta}}{4(-1+e^{i\beta})}\frac{e^{-i\beta(1+k)}(e^{i\beta k}-e^{i\beta})}{\beta}d\beta\biggl(\frac{ik}{(-1+k)^{2}}-\frac{i}{k(-1+k)^{2}}\biggr) \\
    +&\int_{-\pi}^{\pi}\frac{(1+2e^{i\beta}-e^{2i\beta})e^{-i\beta k}}{4(-1+e^{i\beta})}d\beta\biggl(\frac{k}{-1+k}-\frac{1}{k(-1+k)}\biggr) \\
    +&\int_{-\pi}^{\pi}\frac{-e^{-i\beta}(-1+2e^{i\beta}+e^{2i\beta})}{4(-1+e^{i\beta})}\frac{e^{-i\beta k}(-1+e^{i\beta(1+k)})}{\beta}d\beta\biggl(\frac{ik}{(1+k)^{2}}-\frac{i}{k(1+k)^{2}}\biggr) \\
    +&\int_{-\pi}^{\pi}\frac{-e^{-i\beta}(-1+2e^{i\beta}+e^{2i\beta})e^{-i\beta k}}{4(-1+e^{i\beta})}d\beta\biggl(\frac{k}{1+k}-\frac{1}{k(1+k)}\biggr) \\
    +&\int_{-\pi}^{\pi}\frac{-e^{-i\beta}i(1+e^{i\beta}+e^{2i\beta}+e^{3i\beta})e^{-ik\beta}}{4(-1+e^{i\beta})}d
    \beta\biggl(ik-\frac{i}{k}\biggr)\biggr).
\end{align*}
For \(k \notin \{0, \pm1\}\), we define
\begin{align} 
    &J_{1}(k) = \label{J1} \\
    &\int_{-\pi}^{\pi}\frac{(i-(i+\beta)e^{i\beta})(-1-2e^{i\beta}+e^{2i\beta})}{4(-1+e^{i\beta})^{2}}\frac{e^{-i\beta}-e^{-i\beta k}}{\beta}d\beta\biggl(\frac{k}{k-1}+\frac{1}{k(1-k)}\biggr) \nonumber \\
    +&\int_{-\pi}^{\pi}\frac{e^{-i\beta}(-1+2e^{i\beta}+e^{2i\beta})(\beta+i(-1+e^{i\beta}))}{4(-1+e^{i\beta})^{2}}\frac{e^{i\beta}-e^{-i\beta k}}{\beta}d\beta \nonumber \\
    &\cdot\biggl(\frac{k}{1+k}-\frac{1}{k(1+k)}\biggr) \nonumber \\
    +&\int_{-\pi}^{\pi}\frac{-i(-1-2e^{i\beta}+e^{2i\beta})}{4(-1+e^{i\beta})}\frac{e^{-i\beta(1+k)}(e^{i\beta k} - e^{i\beta})}{\beta}d\beta \nonumber \\
    &\cdot\biggl(\frac{-k^{2}}{(-1+k)^{2}}+\frac{1}{(-1+k)^{2}}\biggr) \nonumber \\
    +&\int_{-\pi}^{\pi}\frac{-i(-1-2e^{i\beta}+e^{2i\beta})e^{-i\beta k}}{4(-1+e^{i\beta})}d\beta\biggl(\frac{ik^{2}}{-1+k}-\frac{i}{-1+k}\biggr) \nonumber \\
    +&\int_{-\pi}^{\pi}\frac{e^{-i\beta}i(-1+2e^{i\beta}+e^{2i\beta})}{4(-1+e^{i\beta})}\frac{e^{-i\beta k}(-1+e^{i\beta(1+k)})}{\beta}d\beta \nonumber \\
    &\cdot\biggl(\frac{-k^{2}}{(1+k)^{2}}+\frac{1}{(1+k)^{2}}\biggr) \nonumber \\
    +&\int_{-\pi}^{\pi}\frac{e^{-i\beta}(-1+2e^{i\beta}+e^{2i\beta})e^{-i\beta k}}{4(-1+e^{i\beta})}d\beta\biggl(\frac{-k^{2}}{1+k}+\frac{1}{1+k}\biggr) \nonumber \\
    +&\int_{-\pi}^{\pi}\frac{1+2e^{i\beta}-e^{2i\beta}}{4(-1+e^{i\beta})}\frac{e^{-i\beta(1+k)}(e^{i\beta k}-e^{i\beta})}{\beta}d\beta\biggl(\frac{ik}{(-1+k)^{2}}-\frac{i}{k(-1+k)^{2}}\biggr) \nonumber \\
    +&\int_{-\pi}^{\pi}\frac{(1+2e^{i\beta}-e^{2i\beta})e^{-i\beta k}}{4(-1+e^{i\beta})}d\beta\biggl(\frac{k}{-1+k}-\frac{1}{k(-1+k)}\biggr) \nonumber \\
    +&\int_{-\pi}^{\pi}\frac{-e^{-i\beta}(-1+2e^{i\beta}+e^{2i\beta})}{4(-1+e^{i\beta})}\frac{e^{-i\beta k}(-1+e^{i\beta(1+k)})}{\beta}d\beta \nonumber \\
    &\cdot\biggl(\frac{ik}{(1+k)^{2}}-\frac{i}{k(1+k)^{2}}\biggr) \nonumber \\
    +&\int_{-\pi}^{\pi}\frac{-e^{-i\beta}(-1+2e^{i\beta}+e^{2i\beta})e^{-i\beta k}}{4(-1+e^{i\beta})}d\beta\biggl(\frac{k}{1+k}-\frac{1}{k(1+k)}\biggr) \nonumber
\end{align}
and
\begin{align*}
    J_{2}(k)=\int_{-\pi}^{\pi}\frac{-e^{-i\beta}i(1+e^{i\beta}+e^{2i\beta}+e^{3i\beta})e^{-ik\beta}}{4(-1+e^{i\beta})}d
    \beta\biggl(ik-\frac{i}{k}\biggr).
\end{align*}
Then for \(\abs{k}>1\) we can write the \(k\)th Fourier mode of \(\mathcal{L}\) as
\begin{align} \label{ftl}
    \mathcal{F}(\mathcal{L})(k) = \frac{2\pi}{L(t)}\frac{\gamma}{4\pi}\mathcal{F}(\phi)(k)\biggl(J_{1}(k)+J_{2}(k)\biggr).
\end{align}
Since (\ref{operatorL}) is real, for \(k \in \mathbb{Z}^{+}\),
\begin{align} \label{positivekonly}
    \mathcal{F}(\mathcal{L})(-k)=\frac{1}{2\pi}\int_{-\pi}^{\pi}\mathcal{L}(\alpha)e^{ik\alpha}d\alpha=\overline{\frac{1}{2\pi}\int_{-\pi}^{\pi}\mathcal{L}(\alpha)e^{-ik\alpha}d\alpha}=\overline{\mathcal{F}(\mathcal{L})(k)}.
\end{align}
Hence, it suffices to compute \(\mathcal{F}(\mathcal{L})(k)\) only for \(k > 1\).

\subsubsection{Computing \(J_{2}(k)\)}
\label{computingJ2}
For \(J_{2}(k)\), it suffices to calculate the integral
\begin{align*}
    &ik\int_{-\pi}^{\pi}\frac{-e^{-i\beta}i(1+e^{i\beta}+e^{2i\beta}+e^{3i\beta})e^{-ik\beta}}{4(-1+e^{i\beta})}d\beta.
\end{align*}
Using that
\begin{align*}
    \frac{1}{-1+re^{i\beta}}=-\frac{1}{1-re^{i\beta}}=-\sum_{n=0}^{\infty}(re^{i\beta})^{n},
\end{align*}
we obtain
\begin{align*}
    &ik\int_{-\pi}^{\pi}\frac{-e^{-i\beta}i(1+e^{i\beta}+e^{2i\beta}+e^{3i\beta})e^{-ik\beta}}{4(-1+e^{i\beta})}d\beta \\
    =&-\frac{k}{4}\lim_{\epsilon \to 0^{+}}\lim_{r \to 1^{-}}\sum_{n=0}^{\infty}r^{n}\int_{\substack{[-\pi,\pi] \\ \setminus(-\epsilon, \epsilon)}}e^{-i\beta(k+1)}(1+e^{i\beta}+e^{2i\beta}+e^{3i\beta})e^{i\beta n}d\beta.
\end{align*}
To calculate the outer integral, we note that
\begin{align*}
    &\int_{-\pi}^{\pi}e^{-i\beta(k+1)}(1+e^{i\beta}+e^{2i\beta}+e^{3i\beta})e^{i\beta n}d\beta =
    \begin{cases}
    0 & \mbox{if } n \notin \{k+1,k,k-1,k-2\}, \\
    2\pi & \mbox{otherwise}.
    \end{cases}
\end{align*}
A similar calculation shows that
\begin{align*}
    &-\frac{k}{4}\lim_{\epsilon\to0^{+}}\lim_{r\to1^{-}}\sum_{n=0}^{\infty}r^{n}\int_{\epsilon}^{-\epsilon}e^{-i\beta(k+1)}(1+e^{i\beta}+e^{2i\beta}+e^{3i\beta})e^{i\beta n} d\beta = k\pi.
\end{align*}
Therefore,
\begin{align} \label{valueofJ2}
    J_{2}(k)=
    \begin{cases}
    \pi\biggl(k-\frac{1}{k}\biggr) &\mbox{if } k<-1, \\
    -\pi\biggl(k-\frac{1}{k}\biggr) &\mbox{if } k>1.
    \end{cases}
\end{align}

\subsubsection{Computing \(J_{1}(k)\)}
\label{computingJ1}
In view of (\ref{positivekonly}), we assume that \(k>1\). In this Section, we adopt the notational convention that any summation \(\sum\) in which the upper bound is strictly less than the lower bound is defined to be \(0\). For example, if \(k=2\), then (\ref{sumexample}) vanishes. To begin, we simplify the first two integrals in (\ref{J1}). The first integral can be written as
\begin{align*}
    &\int_{-\pi}^{\pi}\frac{(i-(i+\beta)e^{i\beta})(-1-2e^{i\beta}+e^{2i\beta})}{4(-1+e^{i\beta})^{2}}\frac{e^{-i\beta}-e^{-i\beta k}}{\beta}d\beta \\
    =&\int_{-\pi}^{\pi}\frac{-i(-1-2e^{i\beta}+e^{2i\beta})}{4(-1+e^{i\beta})}\frac{e^{-i\beta} - e^{-i\beta k}}{\beta}d\beta \\
    &+ \int_{-\pi}^{\pi}\frac{-\beta e^{i\beta}(-1-2e^{i\beta}+e^{2i\beta})}{4(-1+e^{i\beta})^{2}}\frac{e^{-i\beta}-e^{-i\beta k}}{\beta}d\beta,
\end{align*}
while the second integral can be written as
\begin{align*}
    &\int_{-\pi}^{\pi}\frac{e^{-i\beta}(-1+2e^{i\beta}+e^{2i\beta})(\beta+i(-1+e^{i\beta}))}{4(-1+e^{i\beta})^{2}}\frac{e^{i\beta}-e^{-i\beta k}}{\beta}d\beta \\
    =&\int_{-\pi}^{\pi}\frac{e^{-i\beta}(-1+2e^{i\beta}+e^{2i\beta})\beta}{4(-1+e^{i\beta})^{2}}\frac{e^{i\beta} - e^{-i\beta k}}{\beta}d\beta \\
    &+ \int_{-\pi}^{\pi}\frac{(-1+2e^{i\beta}+e^{2i\beta})i}{4(-1+e^{i\beta})}\frac{1-e^{-i\beta(k+1)}}{\beta}d\beta.
\end{align*}
For ease of notation, let us define
\begin{align*}
    g_{2}(k) &= \int_{-\pi}^{\pi}\frac{-\beta e^{i\beta}(-1-2e^{i\beta}+e^{2i\beta})}{4(-1+e^{i\beta})^{2}}\frac{e^{-i\beta}-e^{-i\beta k}}{\beta}d\beta, \\
    g_{3}(k) &= \int_{-\pi}^{\pi}\frac{e^{-i\beta}(-1+2e^{i\beta}+e^{2i\beta})\beta}{4(-1+e^{i\beta})^{2}}\frac{e^{i\beta} - e^{-i\beta k}}{\beta}d\beta, \\
    g_{5}(k) &= \int_{-\pi}^{\pi}\frac{-i(-1-2e^{i\beta}+e^{2i\beta})e^{-i\beta k}}{4(-1+e^{i\beta})}d\beta, \\
    g_{6}(k) &= \int_{-\pi}^{\pi}\frac{e^{-i\beta}(-1+2e^{i\beta}+e^{2i\beta})e^{-i\beta k}}{4(-1+e^{i\beta})}d\beta, \\
    g_{7}(k) &= \int_{-\pi}^{\pi}\frac{(1+2e^{i\beta}-e^{2i\beta})e^{-i\beta k}}{4(-1+e^{i\beta})}d\beta, \\
    g_{8}(k) &= \int_{-\pi}^{\pi}\frac{-e^{-i\beta}(-1+2e^{i\beta}+e^{2i\beta})e^{-i\beta k}}{4(-1+e^{i\beta})}d\beta.
\end{align*}
Then we can rewrite
\begin{align} 
    J_{1}(k)=&\frac{k+1}{k}g_{2}(k)+\frac{k-1}{k}g_{3}(k)+i(k+1)g_{5}(k)+(1-k)g_{6}(k) \label{J1simplified} \\
    &+\frac{k+1}{k}g_{7}(k)+\frac{k-1}{k}g_{8}(k). \nonumber
\end{align}
Let us first compute \(g_{2}(k)\). Using that
\begin{align*}
    1-e^{-i\beta(k-1)} &= i\beta(k-1)\int_{0}^{1}e^{-i\beta(k-1)s}ds \\
    \frac{1}{(-1+re^{i\beta})^{2}} &= \sum_{n=0}^{\infty}(1+n)(re^{i\beta})^{n},
\end{align*}
we obtain
\begin{align*}
    g_{2}(k)=&\frac{i(k-1)}{4} \\
    &\cdot\lim_{\epsilon \to 0^{+}}\lim_{r \to 1^{-}}\sum_{n=0}^{\infty}(1+n)r^{n}\int_{0}^{1}\int_{\substack{[-\pi,\pi] \\ \setminus (-\epsilon, \epsilon)}}(1+2e^{i\beta}-e^{2i\beta})\beta e^{-i\beta(k-1)s}e^{i\beta n}d\beta ds.
\end{align*}
To simplify this expression, we first calculate the integral
\begin{align*}
    &\int_{-\pi}^{\pi}(1+2e^{i\beta}-e^{2i\beta})\beta e^{-i\beta(k-1)s}e^{i\beta n}d\beta \\
    =&\frac{e^{-i\pi(n+s(1-k))}}{(n+s(1-k))^{2}}\biggl(-1-i\pi(n+s(1-k))\biggr) \\
    &+\frac{e^{i\pi(n+s(1-k))}}{(n+s(1-k))^{2}}\biggl(1-i\pi(n+s(1-k))\biggr) \\
    &+\frac{2e^{-i\pi(1+n+s(1-k))}}{(1+n+s(1-k))^{2}}\biggl(-1-i\pi(1+n+s(1-k))\biggr) \\
    &+\frac{2e^{i\pi(1+n+s(1-k))}}{(1+n+s(1-k))^{2}}\biggl(1-i\pi(1+n+s(1-k))\biggr) \\
    &+\frac{e^{-i\pi(2+n+s(1-k))}}{(2+n+s(1-k))^{2}}\biggl(1+i\pi(2+n+s(1-k))\biggr) \\
    &+\frac{e^{i\pi(2+n+s(1-k))}}{(2+n+s(1-k))^{2}}\biggl(-1+i\pi(2+n+s(1-k))\biggr).
\end{align*}
For \(t \in \{1, 2\}\), we note that
\begin{align*}
    &\int_{0}^{1}\frac{e^{i\pi(n+s(1-k))}}{(n+s(1-k))^{2}}(1-i\pi(n+s(1-k)))ds = \frac{1}{k-1}\int_{n-(k-1)}^{n}\frac{e^{i\pi s}(1-i\pi s)}{s^{2}}ds, \\
    &\int_{0}^{1}\frac{e^{i\pi(t+n+s(1-k))}}{(t+n+s(1-k))^{2}}(1-i\pi(t+n+s(1-k)))ds \\
    = &\frac{1}{k-1}\int_{t+n-(k-1)}^{t+n}\frac{e^{i\pi s}(1-i\pi s)}{s^{2}}ds.
\end{align*}
Then
\begin{align}
    &\frac{i(k-1)}{4}\sum_{n=0}^{\infty}(1+n)r^{n}\int_{0}^{1}\int_{-\pi}^{\pi}(1+2e^{i\beta}-e^{2i\beta})\beta e^{-i\beta(k-1)s}e^{i\beta n}d\beta ds \nonumber \\
    =&\frac{1}{4}\biggl(\sum_{n=0}^{k-1}(1+n)r^{n}\int_{-\pi}^{\pi}e^{i\beta n}(1-e^{-i\beta(k-1)})d\beta \nonumber \\
    &+ \sum_{n=0}^{k-2}(1+n)r^{n}\int_{-\pi}^{\pi}2e^{i\beta}e^{i\beta n}(1-e^{-i\beta(k-1)})d\beta \nonumber \\
    &+\sum_{n=0}^{k-3}(1+n)r^{n}\int_{-\pi}^{\pi}-e^{2i\beta}e^{i\beta n}(1-e^{-i\beta(k-1)})d\beta\biggr) \label{sumexample} \\
    +&\frac{1}{4}\biggl(\sum_{n=k}^{\infty}(1+n)r^{n}\int_{-\infty}^{\infty}i\biggl(\frac{e^{i\pi s}(1-i\pi s)}{s^{2}} - \frac{e^{-i\pi s}(1+i\pi s)}{s^{2}}\biggr) \label{t1} \\
    &\cdot 1_{[n-(k-1), n]}(s)ds \nonumber
\end{align}
\begin{align}
    +&\sum_{n=k-1}^{\infty}(1+n)r^{n}\cdot 2\int_{-\infty}^{\infty}i\biggl(\frac{e^{i\pi s}(1-i\pi s)}{s^{2}} - \frac{e^{-i\pi s}(1+i\pi s)}{s^{2}}\biggr) \label{t2} \\
    & \cdot 1_{[1+n-(k-1), 1+n]}(s)ds \nonumber \\
    +&\sum_{n=k-2}^{\infty}(1+n)r^{n}(-1) \nonumber \\
    &\cdot\int_{-\infty}^{\infty}i\biggl(\frac{e^{i\pi s}(1-i\pi s)}{s^{2}} - \frac{e^{-i\pi s}(1+i\pi s)}{s^{2}}\biggr)1_{[2+n-(k-1), 2+n]}(s)ds\biggr). \label{t3}
\end{align}
We will further simplify the terms in (\ref{t1}), (\ref{t2}), and (\ref{t3}). The term in (\ref{t1}) becomes
\begin{align*}
    &\int_{-\infty}^{\infty}i\biggl(\frac{e^{i\pi s}(1-i\pi s)}{s^{2}} - \frac{e^{-i\pi s}(1+i\pi s)}{s^{2}}\biggr)\sum_{n=k}^{\infty}(1+n)r^{n}1_{[n-(k-1), n]}(s)ds \\
    =& \sum_{j=1}^{k-2}\int_{j}^{1+j}i\biggl(\frac{e^{i\pi s}(1-i\pi s)}{s^{2}} - \frac{e^{-i\pi s}(1+i\pi s)}{s^{2}}\biggr)\sum_{n=k}^{k+j-1}(1+n)r^{n}ds \\
    &+\sum_{j=1}^{\infty}\int_{k+j-2}^{k+j-1}i\biggl(\frac{e^{i\pi s}(1-i\pi s)}{s^{2}} - \frac{e^{-i\pi s}(1+i\pi s)}{s^{2}}\biggr)\sum_{n=k+(j-1)}^{2k-2+(j-1)}(1+n)r^{n}ds.
\end{align*}
Next, the term in (\ref{t2}) becomes
\begin{align*}
    &2\int_{-\infty}^{\infty}i\biggl(\frac{e^{i\pi s}(1-i\pi s)}{s^{2}} - \frac{e^{-i\pi s}(1+i\pi s)}{s^{2}}\biggr)\sum_{n=k-1}^{\infty}(1+n)r^{n}1_{[1+n-(k-1), 1+n]}(s)ds \\
    =&\sum_{j=1}^{k-2}2\int_{j}^{1+j}i\biggl(\frac{e^{i\pi s}(1-i\pi s)}{s^{2}} - \frac{e^{-i\pi s}(1+i\pi s)}{s^{2}}\biggr)\sum_{n=k-1}^{k+j-2}(1+n)r^{n}ds \\
    &+2\sum_{j=1}^{\infty}\int_{k+j-2}^{k+j-1}i\biggl(\frac{e^{i\pi s}(1-i\pi s)}{s^{2}} - \frac{e^{-i\pi s}(1+i\pi s)}{s^{2}}\biggr)\sum_{n=k-1+(j-1)}^{2k-3+(j-1)}(1+n)r^{n}ds.
\end{align*}
Lastly, the term in (\ref{t3}) becomes
\begin{align*}
    -&\int_{-\infty}^{\infty}i\biggl(\frac{e^{i\pi s}(1-i\pi s)}{s^{2}} - \frac{e^{-i\pi s}(1+i\pi s)}{s^{2}}\biggr)\sum_{n=k-2}^{\infty}(1+n)r^{n}1_{[2+n-(k-1), 2+n]}(s)ds \\
    =-&\sum_{j=1}^{k-2}\int_{j}^{1+j}i\biggl(\frac{e^{i\pi s}(1-i\pi s)}{s^{2}} - \frac{e^{-i\pi s}(1+i\pi s)}{s^{2}}\biggr)\sum_{n=k-2}^{k+j-3}(1+n)r^{n}ds \\
    -&\sum_{j=1}^{\infty}\int_{k+j-2}^{k+j-1}i\biggl(\frac{e^{i\pi s}(1-i\pi s)}{s^{2}} - \frac{e^{-i\pi s}(1+i\pi s)}{s^{2}}\biggr)\sum_{n=k-2+(j-1)}^{2k-4+(j-1)}(1+n)r^{n}ds.
\end{align*}
Using that
\begin{align*}
    \frac{1}{4}\biggl(&\sum_{n=0}^{k-1}(1+n)\int_{-\pi}^{\pi}e^{i\beta n}(1-e^{-i\beta(k-1)})d\beta \\
    +&\sum_{n=0}^{k-2}(1+n)\int_{-\pi}^{\pi}2e^{i\beta}e^{i\beta n}(1-e^{-i\beta(k-1)})d\beta \\
    +&\sum_{n=0}^{k-3}(1+n)\int_{-\pi}^{\pi}-e^{2i\beta}e^{i\beta n}(1-e^{-i\beta(k-1)})d\beta\biggr) = \frac{\pi}{2}-\pi k,
\end{align*}
we obtain
\begin{align*}
    &\lim_{\epsilon \to 0^{+}}\lim_{r \to 1^{-}}\frac{i(k-1)}{4}\sum_{n=0}^{\infty}(1+n)r^{n}\int_{0}^{1}\int_{-\pi}^{\pi}(1+2e^{i\beta}-e^{2i\beta})\beta e^{-i\beta(k-1)s}e^{i\beta n}d\beta ds \\
    =&\frac{1}{4}\biggl((2k+1)\sum_{j=1}^{k-2}j\int_{j}^{1+j}i\biggl(\frac{e^{i\pi s}(1-i\pi s)}{s^{2}} - \frac{e^{-i\pi s}(1+i\pi s)}{s^{2}}\biggr)ds \\
    &+\sum_{j=1}^{k-2}j^{2}\int_{j}^{1+j}i\biggl(\frac{e^{i\pi s}(1-i\pi s)}{s^{2}} - \frac{e^{-i\pi s}(1+i\pi s)}{s^{2}}\biggr)ds \\
    &+(-1+k)(-2+3k)\sum_{j=1}^{\infty}\int_{k+j-2}^{k+j-1}i\biggl(\frac{e^{i\pi s}(1-i\pi s)}{s^{2}} - \frac{e^{-i\pi s}(1+i\pi s)}{s^{2}}\biggr)ds \\
    &+2(-1+k)\sum_{j=1}^{\infty}j\int_{k+j-2}^{k+j-1}i\biggl(\frac{e^{i\pi s}(1-i\pi s)}{s^{2}} - \frac{e^{-i\pi s}(1+i\pi s)}{s^{2}}\biggr)ds\biggr) \\
    &+\frac{\pi}{2} -\pi k.
\end{align*}
Using similar calculations, we obtain
\begin{align*}
    &\lim_{\epsilon \to 0^{+}}\lim_{r \to 1^{-}}\frac{i(k-1)}{4}\sum_{n=0}^{\infty}(1+n)r^{n}\int_{0}^{1}\int_{\epsilon}^{-\epsilon}(1+2e^{i\beta}-e^{2i\beta})\beta e^{-i\beta(k-1)s}e^{i\beta n}d\beta ds \\
    = &\frac{\pi}{2}(k-1).
\end{align*}
Therefore,
\begin{align*}
    g_{2}(k) = -&\frac{\pi}{2}k+\frac{1}{4}\biggl((2k+1)\sum_{j=1}^{k-2}j\int_{j}^{1+j}i\biggl(\frac{e^{i\pi s}(1-i\pi s)}{s^{2}} - \frac{e^{-i\pi s}(1+i\pi s)}{s^{2}}\biggr)ds \\
    +&\sum_{j=1}^{k-2}j^{2}\int_{j}^{1+j}i\biggl(\frac{e^{i\pi s}(1-i\pi s)}{s^{2}} - \frac{e^{-i\pi s}(1+i\pi s)}{s^{2}}\biggr)ds \\
    +&(-1+k)(-2+3k)\sum_{j=1}^{\infty}\int_{k+j-2}^{k+j-1}i\biggl(\frac{e^{i\pi s}(1-i\pi s)}{s^{2}} - \frac{e^{-i\pi s}(1+i\pi s)}{s^{2}}\biggr)ds \\
    +&2(-1+k)\sum_{j=1}^{\infty}j\int_{k+j-2}^{k+j-1}i\biggl(\frac{e^{i\pi s}(1-i\pi s)}{s^{2}} - \frac{e^{-i\pi s}(1+i\pi s)}{s^{2}}\biggr)ds\biggr).
\end{align*}
To simplify the expression for \(g_{2}(k)\), we note that for \(j \neq \{-1,0\}\),
\begin{align*}
    &\int_{j}^{j+1}i\biggl(\frac{e^{i\pi s}(1-i\pi s)}{s^{2}}-\frac{e^{-i\pi s}(1+i\pi s)}{s^{2}}\biggr)ds=0.
\end{align*}
Then \(g_{2}(k)=-\frac{\pi}{2}k\). Similarly, we can calculate that
\begin{align*}
    g_{3}(k)&=-\frac{\pi}{2}k, \quad g_{5}(k)= -\frac{i\pi}{2}, \quad g_{6}(k)=-\frac{\pi}{2}, \quad g_{7}(k)=-\frac{\pi}{2}, \quad g_{8}(k)=\frac{\pi}{2}.
\end{align*}
Plugging these values into (\ref{J1simplified}), we obtain
\begin{align*}
    J_{1}(k) = -\frac{\pi}{k}.
\end{align*}
Using (\ref{positivekonly}) and (\ref{valueofJ2}), we deduce that
\begin{align} \label{valueofJ1}
    J_{1}(k)=
    \begin{cases}
        -\frac{\pi}{k} & k > 1, \\
        \frac{\pi}{k} & k < -1.
    \end{cases}
\end{align}

\subsection{Summary}
Plugging the results of Sections \ref{computingJ2} and \ref{computingJ1} into (\ref{ftl}), we obtain that for \(k >1\),
\begin{align*}
    \mathcal{F}(\mathcal{L})(k) = -\frac{2\pi}{L(t)}\frac{\gamma}{4\pi}\mathcal{F}(\phi)(k)\pi k =-\frac{2\pi}{L(t)}\frac{\gamma}{4\pi}\mathcal{F}(\phi)(k)\pi \abs{k}.
\end{align*}
Since for \(k > 1\)
\begin{align*}
    \mathcal{F}(\mathcal{L})(-k) &= \overline{\mathcal{F}(\mathcal{L})(k)} = -\frac{2\pi}{L(t)}\frac{\gamma}{4\pi}\overline{\mathcal{F}(\phi)(k)}\pi k = -\frac{2\pi}{L(t)}\frac{\gamma}{4\pi}\mathcal{F}(\phi)(-k)\pi \abs{k},
\end{align*}
we conclude that for \(\abs{k} > 1\),
\begin{align}
    \mathcal{F}(\mathcal{L})(k) = -\frac{2\pi}{L(t)}\frac{\gamma}{4\pi}\mathcal{F}(\phi)(k)\pi \abs{k}. \nonumber
\end{align}
To calculate \(\mathcal{F}(\mathcal{L})(k)\) for \(\abs{k}=1\), we note that for \(k \in \mathbb{Z}\),
\begin{align*}
    \mathcal{F}((U_{1})_{\alpha})(k) = ik\mathcal{F}(U_{1})(k).
\end{align*}
Then (\ref{fourierlinear}) can be written as
\begin{align*}
    \mathcal{F}(\mathcal{L})(k) = \frac{2\pi}{L(t)}\biggl(ik\mathcal{F}(U_{1})(k) - \frac{i}{k}\mathcal{F}(U_{1})(k)\biggr),
\end{align*}
where \(k \neq 0\). Hence, \(\mathcal{F}(\mathcal{L})(\pm 1)=0\). This presents a technical challenge in dealing with the term in the evolution equation for \(\theta\) that induces an exponential decay in time of the initial perturbation of the interface. This challenge can be resolved by observing that the identity
\begin{align} \label{implicittheoremidentity}
    \int_{-\pi}^{\pi}z_{\alpha}(\alpha,t)d\alpha = 0
\end{align}
provides a means to control the \(\pm 1\) Fourier modes of \(\mathcal{L}\) using the other nonzero Fourier modes.

\section{Derivation of an \emph{a priori} Estimate}
\label{derivationofthekeynorminequality}
Before embarking on the derivation of a key \emph{a priori} estimate for \(\phi = \theta -\hat{\theta}(0)\), we first derive tight bounds around \(L(t)\), which remains valid as long as \(\norm{\phi(t)}_{\mathcal{F}^{0,1}}\) is sufficiently small for all \(t \geq 0\).
\begin{proposition} \label{tightcontrolL}
If \(\norm{\phi(t)}_{\mathcal{F}^{0,1}}\) is sufficiently small for all \(t \geq 0\), then
    \begin{align}
        \frac{R^{2}}{1+\frac{\pi}{2}(e^{2\normzz{\phi(t)}}-1)} \leq \biggl(\frac{L(t)}{2\pi}\biggr)^{2} \leq \frac{R^{2}}{1-\frac{\pi}{2}(e^{2\normzz{\phi(t)}}-1)}. \nonumber
    \end{align}
\end{proposition}
\begin{proof}
    Using integration by parts, we obtain
\begin{align*}
    &\mathcal{F}\biggl(\int_{0}^{\alpha}e^{-i\eta}(\phi(\alpha)-\phi(\eta))^{n}d\eta\biggr)(-1) =i\mathcal{F}((\phi(\pi)-\phi(\eta))^{n})(1).
\end{align*}
Then
\begin{align*}
    &\int_{-\pi}^{\pi}\int_{0}^{\alpha}e^{i(\alpha-\eta)}\sum_{n \geq 1}\frac{i^{n}}{n!}(\phi(\alpha)-\phi(\eta))^{n}d\eta d\alpha = 2\pi i\sum_{n \geq 1}\frac{i^{n}}{n!}\mathcal{F}((\phi(\pi)-\phi(\eta))^{n})(1).
\end{align*}
Hence,
\begin{align*}
    &\mbox{Im}\biggl(\int_{-\pi}^{\pi}\int_{0}^{\alpha}e^{i(\alpha-\eta)}\sum_{n \geq 1}\frac{i^{n}}{n!}(\phi(\alpha)-\phi(\eta))^{n}d\eta d\alpha\biggr) \\
    =&\pi\biggl(\sum_{n \geq 1}\frac{i^{n}}{n!}\mathcal{F}((\phi(\pi)-\phi(\eta))^{n})(1)+\sum_{n \geq 1}\frac{(-i)^{n}}{n!}\overline{\mathcal{F}((\phi(\pi)-\phi(\eta))^{n})(1)}\biggr).
\end{align*}
It follows that
\begin{align*}
    &\abs{\mbox{Im}\biggl(\int_{-\pi}^{\pi}\int_{0}^{\alpha}e^{i(\alpha-\eta)}\sum_{n \geq 1}\frac{i^{n}}{n!}(\phi(\alpha)-\phi(\eta))^{n}d\eta d\alpha\biggr)} \leq 2\pi\sum_{n\geq 1}\frac{\normzz{(\phi(\pi)-\phi(\cdot))^{n}}}{n!}.
\end{align*}
By Proposition \ref{p5},
\begin{align*}
    \normzz{(\phi(\pi)-\phi(\cdot))^{n}} \leq \normzz{\phi(\pi)-\phi(\cdot)}^{n}.
\end{align*}
Then
\begin{align*}
    \abs{\mbox{Im}\biggl(\int_{-\pi}^{\pi}\int_{0}^{\alpha}e^{i(\alpha-\eta)}\sum_{n \geq 1}\frac{i^{n}}{n!}(\phi(\alpha)-\phi(\eta))^{n}d\eta d\alpha\biggr)} &\leq 2\pi(e^{\phi(\pi)}e^{\normzz{\phi}}-1).
\end{align*}
Since
\begin{align*}
    \abs{\phi(\pi)} \leq \sum_{k \in \mathbb{Z}}\abs{\hat{\phi}(k)}=\normzz{\phi},
\end{align*}
we obtain
\begin{align}
    \abs{\mbox{Im}\biggl(\int_{-\pi}^{\pi}\int_{0}^{\alpha}e^{i(\alpha-\eta)}\sum_{n \geq 1}\frac{i^{n}}{n!}(\phi(\alpha)-\phi(\eta))^{n}d\eta d\alpha\biggr)} &\leq \pi^{2}(e^{2\normzz{\phi}}-1). \nonumber
\end{align}
This estimate shows that
\begin{align*}
    \frac{R^{2}}{1+\frac{\pi}{2}(e^{2\norm{\phi}_{\mathcal{F}^{0,1}}}-1)} \leq \biggl(\frac{L(t)}{2\pi}\biggr)^{2} \leq \frac{R^{2}}{1-\frac{\pi}{2}(e^{2\norm{\phi}_{\mathcal{F}^{0,1}}}-1)},
\end{align*}
as needed.
\end{proof}
Using Proposition \ref{tightcontrolL}, we can also prove the following useful estimate.
\begin{proposition} \label{p6}
    For sufficiently small \(\norm{\phi}_{\mathcal{F}^{0,1}}\),    \begin{align*}
        \abs{R\frac{2\pi}{L(t)}-1} \leq 1- \sqrt{1-\frac{\pi}{2}(e^{2\norm{\phi}_{\mathcal{F}^{0,1}}}-1)}.
    \end{align*}
\end{proposition}
\begin{proof}
    From Proposition \ref{tightcontrolL}, we obtain
    \begin{align*}
        \sqrt{1-\frac{\pi}{2}(e^{2\norm{\phi}_{\mathcal{F}^{0,1}}}-1)}-1 \leq \frac{2\pi R}{L(t)}-1 \leq \sqrt{1+\frac{\pi}{2}(e^{2\norm{\phi}_{\mathcal{F}^{0,1}}}-1)}-1.
    \end{align*}
    Then
    \begin{align*}
        \abs{\frac{2\pi R}{L(t)}-1} &\leq \max\biggl\{\abs{\sqrt{1-\frac{\pi}{2}(e^{2\norm{\phi}_{\mathcal{F}^{0,1}}}-1)}-1}, \abs{\sqrt{1+\frac{\pi}{2}(e^{2\norm{\phi}_{\mathcal{F}^{0,1}}}-1)}-1}\biggr\} \\
        &=1-\sqrt{1-\frac{\pi}{2}(e^{2\norm{\phi}_{\mathcal{F}^{0,1}}}-1)},
    \end{align*}
    as needed.
\end{proof}
Now, we derive a key \emph{a priori} estimate for \(\phi\). In Section \ref{principallinearoperatorforthetaequation}, we have shown that
\begin{align*}
    \mathcal{F}(\mathcal{L})(k)&=
    \begin{cases}
        0 & \mbox{if } \abs{k} = 1, \\
        \frac{2\pi}{L(t)}\frac{\gamma}{4\pi}\mathcal{F}(\phi)(k)(J_{1}(k)+J_{2}(k))& \mbox{if } \abs{k} > 1,
    \end{cases}
\end{align*}
where \(J_{1}\) and \(J_{2}\) are given in (\ref{valueofJ1}) and (\ref{valueofJ2}).
Let
\begin{align*}
    \widetilde{\mathcal{L}}(\alpha)&=\frac{L(t)}{2\pi}\mathcal{L}(\alpha), \quad \widetilde{\mathcal{N}}(\alpha)=\frac{L(t)}{2\pi}\mathcal{N}(\alpha).
\end{align*}
Then for \(\abs{k} \geq 1\),
\begin{align*}
    \frac{\partial}{\partial t}\mathcal{F}(\phi)(k)&=
    \begin{cases}
        \frac{2\pi}{L(t)}\mathcal{F}(\widetilde{\mathcal{N}})(k) & \mbox{ if \(\abs{k} = 1\)}, \\
        \frac{2\pi}{L(t)}\frac{\gamma}{4\pi}\mathcal{F}(\phi)(k)(J_{1}(k)+J_{2}(k))+\frac{2\pi}{L(t)}\mathcal{F}(\widetilde{\mathcal{N}})(k) & \mbox{ if \(\abs{k} > 1\)}.
    \end{cases}
\end{align*}
For convenience of notation, define \(J_{1}(k)=J_{2}(k)=0\) for \(\abs{k} = 1\) so that for \(k \in \mathbb{Z} \setminus \{0\}\),
\begin{align} \label{thetaevolFourier}
    \frac{\partial}{\partial t}\mathcal{F}(\phi)(k)=\frac{2\pi}{L(t)}\frac{\gamma}{4\pi}\mathcal{F}(\phi)(k)(J_{1}(k)+J_{2}(k))+\frac{2\pi}{L(t)}\mathcal{F}(\widetilde{\mathcal{N}})(k).
\end{align}
We observe that the principal linear term has a time-dependent coefficient. This dependence occurs, however, only via \(L(t)\). We choose an initial circular interface of radius \(R\) to perturb around and make the principal linear term independent of time by replacing \(L(t)\) with \(2\pi R\) and keeping an error term. That is, we rewrite (\ref{thetaevolFourier}) as
\begin{align} 
    \frac{\partial}{\partial t}\mathcal{F}(\phi)(k)=&\frac{1}{R}\frac{\gamma}{4\pi}\mathcal{F}(\phi)(k)(J_{1}(k)+J_{2}(k))+\frac{2\pi}{L(t)}\mathcal{F}(\widetilde{\mathcal{N}})(k) \label{phihatt} \\
    &+\frac{\gamma}{4\pi}\mathcal{F}(\phi)(k)(J_{1}(k)+J_{2}(k))\biggl(-\frac{1}{R}+\frac{2\pi}{L(t)}\biggr). \nonumber
\end{align}
We note that for \(k > 0\),
\begin{align} \label{absvalssame}
    \abs{\hat{\phi}(-k)}=\abs{\overline{\hat{\phi}(k)}}=\abs{\hat{\phi}(k)}
\end{align}
since \(\phi\) is real-valued. Then for \(s>0\),
\begin{align*}
    \norms{\phi}=\sum_{k\neq 0}e^{\nu(t)\abs{k}}\abs{k}^{s}\abs{\hat{\phi}(k)}=2\sum_{k \geq 1}e^{\nu(t)k}k^{s}\abs{\hat{\phi}(k)}.
\end{align*}
Differentiating this equation with respect to \(t\), we obtain
\begin{align*}
    &\frac{d}{dt}\norms{\phi} = 2\sum_{k\geq 1}e^{\nu(t)k}\nu'(t)k^{s+1}\abs{\hat{\phi}(k)}+2\sum_{k \geq 1}e^{\nu(t)k}k^{s}\frac{\hat{\phi}(k)\overline{\frac{\partial}{\partial t}\hat{\phi}(k)}+\overline{\hat{\phi}(k)}\frac{\partial}{\partial t}\hat{\phi}(k)}{2\abs{\hat{\phi}(k)}}.
\end{align*}
Let us simplify the second term. Using (\ref{phihatt}) and that \(J_{1}\) and \(J_{2}\) are real for \(k \geq 1\), we obtain
\begin{align*}
    &\hat{\phi}(k)\overline{\frac{\partial}{\partial t}\hat{\phi}(k)}+\overline{\hat{\phi}(k)}\frac{\partial}{\partial t}\hat{\phi}(k) \\
    =&\frac{1}{R}\frac{\gamma}{4\pi}(J_{1}+J_{2})(k)\abs{\hat{\phi}(k)}^{2}+\frac{2\pi}{L(t)}\overline{\mathcal{F}(\widetilde{\mathcal{N}})(k)}\hat{\phi}(k) \\
    &+\frac{\gamma}{4\pi}(J_{1}+J_{2})(k)\biggl(-\frac{1}{R}+\frac{2\pi}{L(t)}\biggr)\abs{\hat{\phi}(k)}^{2} \\
    &+\frac{1}{R}\frac{\gamma}{4\pi}(J_{1}+J_{2})(k)\abs{\hat{\phi}(k)}^{2}+\frac{2\pi}{L(t)}\mathcal{F}(\widetilde{\mathcal{N}})(k)\overline{\hat{\phi}(k)} \\
    &+\frac{\gamma}{4\pi}(J_{1}+J_{2})(k)\biggl(-\frac{1}{R}+\frac{2\pi}{L(t)}\biggr)\abs{\hat{\phi}(k)}^{2}.
\end{align*}
Then
\begin{align*}
    &2\sum_{k\geq 1}e^{\nu(t)k}k^{s}\frac{\hat{\phi}(k)\overline{\frac{\partial}{\partial t}\hat{\phi}(k)}+\overline{\hat{\phi}(k)}\frac{\partial}{\partial t}\hat{\phi}(k)}{2\abs{\hat{\phi}(k)}} \\
    =&\frac{2}{R}\frac{\gamma}{4\pi}\sum_{k \geq 1}e^{\nu(t)k}k^{s}(J_{1}+J_{2})(k)\abs{\hat{\phi}(k)} \\
    &+\frac{2\pi}{L(t)}\sum_{k\geq 1}e^{\nu(t)k}k^{s}\frac{\mathcal{F}(\widetilde{\mathcal{N}})(k)\overline{\hat{\theta}(k)}+\overline{\mathcal{F}(\widetilde{\mathcal{N}})(k)}\hat{\phi}(k)}{\abs{\hat{\phi}(k)}} \\
    &+2\frac{\gamma}{4\pi}\biggl(-\frac{1}{R}+\frac{2\pi}{L(t)}\biggr)\sum_{k\geq 1}e^{\nu(t)k}k^{s}(J_{1}+J_{2})(k)\abs{\hat{\phi}(k)}.
\end{align*}
Therefore,
\begin{align}
    &\frac{d}{dt}\norm{\phi}_{\dot{\mathcal{F}}_{\nu}^{s,1}} \nonumber \\
    =& 2\sum_{k \geq 1}e^{\nu(t)k}\nu'(t)k^{s+1}\abs{\hat{\phi}(k)}+\frac{2}{R}\frac{\gamma}{4\pi}\sum_{k \geq 1}e^{\nu(t)k}k^{s}(J_{1}+J_{2})(k)\abs{\hat{\phi}(k)} \label{terms12} \\
    &+\frac{2\pi}{L(t)}\sum_{k \geq 1}e^{\nu(t)k}k^{s}\frac{\mathcal{F}(\widetilde{\mathcal{N}})(k)\overline{\hat{\phi}(k)}+\overline{\mathcal{F}(\widetilde{\mathcal{N}})(k)}\hat{\phi}(k)}{\abs{\hat{\phi}(k)}} \label{term3} \\
    &+2\frac{\gamma}{4\pi}\biggl(-\frac{1}{R}+\frac{2\pi}{L(t)}\biggr)\sum_{k \geq 1}e^{\nu(t)k}k^{s}(J_{1}+J_{2})(k)\abs{\hat{\phi}(k)}. \label{term4}
\end{align}
First, let us estimate (\ref{term4}). Using Proposition \ref{p6}, we obtain
\begin{align}
    &\abs{2\frac{\gamma}{4\pi}\biggl(-\frac{1}{R}+\frac{2\pi}{L(t)}\biggr)\sum_{k \geq 1}e^{\nu(t)k}k^{s}(J_{1}+J_{2})(k)\abs{\hat{\phi}(k)}} \nonumber \\
    \leq&2\pi\frac{\gamma}{4\pi}\frac{1}{R}A\normzz{\phi}\sum_{k \geq 2}e^{\nu(t)k}k^{s+1}\abs{\hat{\phi}(k)}, \label{estimate1}
\end{align}
where
\begin{align*}
    A=A(\normzz{\phi})=\frac{1-\sqrt{1-\frac{\pi}{2}(e^{2\normzz{\phi}}-1)}}{\normzz{\phi}}.
\end{align*}
Next, let us estimate (\ref{terms12}) and (\ref{term3}). First, we note that
\begin{align}
    &2\sum_{k \geq 1}e^{\nu(t) k}\nu'(t)k^{s+1}\abs{\hat{\phi}(k)}+\frac{2}{R}\frac{\gamma}{4\pi}\sum_{k \geq 1}e^{\nu(t)k}k^{s}(J_{1}+J_{2})(k)\abs{\hat{\phi}(k)} \nonumber \\
    &+\frac{2\pi}{L(t)}\sum_{k \geq 1} e^{\nu(t)k}k^{s}\frac{\mathcal{F}(\widetilde{\mathcal{N}})(k)\overline{\hat{\phi}(k)}+\overline{\mathcal{F}(\widetilde{\mathcal{N}})(k)}\hat{\phi}(k)}{\abs{\hat{\phi}(k)}} \nonumber \\    \leq&\nu'(t)\norm{\phi}_{\dot{\mathcal{F}}_{\nu}^{s+1,1}}-\pi\frac{2}{R}\frac{\gamma}{4\pi}\sum_{k \geq 2}e^{\nu(t)k}k^{s+1}\abs{\hat{\phi}(k)}+\frac{2\pi}{L(t)}\norms{\widetilde{\mathcal{N}}}. \label{terms}
\end{align}
Plugging (\ref{estimate1}) and (\ref{terms}) into (\ref{terms12}), we obtain
\begin{align} 
    \frac{d}{dt}\norms{\phi}    \leq&\nu'(t)\normx{\phi}{s+1}-\pi\frac{2}{R}\frac{\gamma}{4\pi}\sum_{k \geq 2}e^{\nu(t)k}k^{s+1}\abs{\hat{\phi}(k)}+\frac{2\pi}{L(t)}\norms{\widetilde{\mathcal{N}}} \label{dissipativetermline} \\
    &+2\frac{\gamma}{4\pi}\frac{1}{R}A\normzz{\phi}\sum_{k \geq 2}e^{\nu(t)k}k^{s+1}\abs{\hat{\phi}(k)}. \nonumber
\end{align}
With the minus sign in the front, the second term in (\ref{dissipativetermline}) is associated with dissipation of the initial interfacial perturbation. It is clear that the \(\pm 1\) Fourier modes of \(\phi\) play no part in the dissipation. This presents a technical difficulty, because the norm of the function space that we intend to use involves all nonzero Fourier modes of \(\phi\). To resolve this issue, we note that (\ref{implicittheoremidentity}) and \(\hat{\phi}(0)=0\) imply
\begin{align*}
    0=\int_{-\pi}^{\pi}e^{i(\alpha+\hat{\phi}(1)e^{i\alpha}+\hat{\phi}(-1)e^{-i\alpha}+\sum_{\abs{k} > 1}\hat{\phi}(k)e^{ik\alpha})}d\alpha.
\end{align*}
This identity provides an implicit relation between the \(\pm 1\) Fourier modes and the other nonzero Fourier modes of \(\phi\), which allows us to control the former in terms of the latter. This observation is summarized in Proposition 4.1 of~\cite{gancedo2019global}. In particular, if \(\normzz{\phi} \in (0,\frac{1}{2}\log\frac{5}{4})\), then for all \(r \in (\normzz{\phi},\frac{1}{2}\log\frac{5}{4})\),
\begin{align*}
    \abs{\hat{\phi}(1)}+\abs{\hat{\phi}(-1)} \leq C_{I}(r)r\sum_{\abs{k} \geq 2}\abs{\hat{\phi}(k)}.
\end{align*}
Then
\begin{align*}
    -\sum_{k \geq 2}e^{\nu(t)k}k^{s+1}\abs{\hat{\phi}(k)} \leq -\frac{1}{2\biggl(C_{I}(\normzz{\phi})\normzz{\phi}+1\biggr)}\normx{\phi}{s+1}.
\end{align*}
Using this estimate and Propositions \ref{p3} and \ref{tightcontrolL}, we obtain
\begin{align} 
    &\frac{d}{dt}\norms{\phi} \label{aprioriestimateforNormsphi} \\
    \leq& \biggl(\nu'(t)-\frac{1}{2\biggl(C_{I}(\normzz{\phi})\normzz{\phi}+1\biggr)}\pi\frac{2}{R}\frac{\gamma}{4\pi}+\frac{\gamma}{4\pi}\frac{1}{R}A\normzz{\phi}\biggr)\normx{\phi}{s+1} \nonumber \\
    &+\frac{1}{R}\frac{1}{A_{1}}\norms{\widetilde{\mathcal{N}}}, \nonumber
\end{align}
where
\begin{align*}
    A_{1} &= \frac{1}{\sqrt{1+\frac{\pi}{2}(e^{2\normzz{\phi}}-1)}}.
\end{align*}

\section{Estimating \texorpdfstring{\(\widetilde{\mathcal{N}}\)}{N}}
In Section \ref{derivationofthekeynorminequality}, we derived an \emph{a priori} estimate containing the \(\dot{\mathcal{F}}_{\nu}^{s,1}\) norm of \(\widetilde{\mathcal{N}}\), where
\begin{align} \label{threeterms}
    \widetilde{\mathcal{N}}(\alpha) = (U_{\geq 2})_{\alpha}(\alpha) + T_{\geq 2}(\alpha)(1+\phi_{\alpha}(\alpha)) + T_{1}(\alpha)\phi_{\alpha}(\alpha).
\end{align}
We derive an estimate for each of the three terms separately. In Sections \ref{estimatingterm2} and \ref{estimaingterm3}, we will derive bounds for the second and third terms which depend on the \(\dot{\mathcal{F}}_{\nu}^{s,1}\) and \(\mathcal{F}_{\nu}^{0,1}\) norms of \(U_{1}\) and \(U_{\geq 2}\). In Sections \ref{estimationU1} and \ref{estimationU2}, respectively, we will derive bounds for these norms in terms of the corresponding norms of \(\phi\). We will estimate the first term in (\ref{threeterms}) in Section \ref{estimationU2deriv}.

\subsection{Estimating \texorpdfstring{\(T_{\geq 2}(\alpha)(1+\phi_{\alpha}(\alpha))\)}{Tg2}}
\label{estimatingterm2}
We prove the following estimates.
\begin{lemma} \label{estimatingterm2lem}
    For \(s \geq 1\),
    \begin{align*}
        &\norms{T_{\geq 2}(1+\phi_{\alpha})} \\
        \leq& \left(1+b(2,s)\normx{\phi}{1}\right) \\
        &\cdot\left(\norm{U_{\geq 2}}_{\dot{\mathcal{F}}_{\nu}^{s-1,1}} + b(2,s-1)\left(\norm{\phi}_{\dot{\mathcal{F}}_{\nu}^{s,1}}\norm{U_{\geq 1}}_{\mathcal{F}_{\nu}^{0,1}} + \norm{U_{\geq 1}}_{\dot{\mathcal{F}}_{\nu}^{s-1,1}}\norm{\phi}_{\dot{\mathcal{F}}_{\nu}^{1,1}}\right)\right) \\
        &+b(2,s)\normx{\phi}{s+1}\left(2\norm{U_{\geq 2}}_{\dot{\mathcal{F}}_{\nu}^{0,1}} + 2\left(\norm{\phi}_{\dot{\mathcal{F}}_{\nu}^{1,1}}\norm{U_{\geq 1}}_{\mathcal{F}_{\nu}^{0,1}} + \norm{\phi}_{\mathcal{F}_{\nu}^{1,1}}\norm{U_{\geq 1}}_{\dot{\mathcal{F}}_{\nu}^{0,1}}\right)\right).
    \end{align*}
    For \(0 \leq s <1\),
    \begin{align*}
        &\norms{T_{\geq 2}(1+\phi_{\alpha})} \\
        \leq&\left(1+b(2,s)\normx{\phi}{1}\right)\\
        &\cdot\left(\norm{U_{\geq 2}}_{\dot{\mathcal{F}}_{\nu}^{0,1}} + b(2,s)\biggl(\norm{\phi}_{\dot{\mathcal{F}}_{\nu}^{1,1}}\norm{U_{\geq 1}}_{\mathcal{F}_{\nu}^{0,1}} + \norm{\phi}_{\dot{\mathcal{F}}_{\nu}^{1,1}}\norm{U_{\geq 1}}_{\dot{\mathcal{F}}_{\nu}^{0,1}}\biggr)\right) \\
        &+b(2,s)\normx{\phi}{s+1}\left(2\norm{U_{\geq 2}}_{\dot{\mathcal{F}}_{\nu}^{0,1}} + 2\left(\norm{\phi}_{\dot{\mathcal{F}}_{\nu}^{1,1}}\norm{U_{\geq 1}}_{\mathcal{F}_{\nu}^{0,1}} + \norm{\phi}_{\mathcal{F}_{\nu}^{1,1}}\norm{U_{\geq 1}}_{\dot{\mathcal{F}}_{\nu}^{0,1}}\right)\right).
    \end{align*}
\end{lemma}
\begin{proof}
    Using Proposition \ref{p5}, we obtain that for \(s \geq 0\),
    \begin{align*}
        &\norms{T_{\geq 2}(1+\phi_{\alpha})} \\
        \leq& \norms{T_{\geq 2}}+b(2,s)\biggl(\norms{T_{\geq 2}}\normx{\phi}{1}+\normx{\phi}{s+1}\normz{T_{\geq 2}}\biggr).
    \end{align*}
    It remains to estimate \(\norm{T_{\geq 2}}_{\dot{\mathcal{F}}_{\nu}^{s,1}}\) and \(\norm{T_{\geq 2}}_{\mathcal{F}_{\nu}^{0,1}} \) in terms of \(U_{1}\) and \(U_{\geq 2}\). Since
    \begin{align*}
        T_{\geq 2}(\alpha)=&\mathcal{M}(U_{\geq 2})(\alpha)+\mathcal{M}(\phi_{\alpha}U_{\geq 1})(\alpha),
    \end{align*}
    using Proposition \ref{p5} we obtain that for \(s \geq 1\),
    \begin{align*}
        \norm{T_{\geq 2}}_{\dot{\mathcal{F}}_{\nu}^{s,1}} \leq&\norm{U_{\geq 2}}_{\dot{\mathcal{F}}_{\nu}^{s-1,1}} + b(2,s-1)(\norm{\phi}_{\dot{\mathcal{F}}_{\nu}^{s,1}}\norm{U_{\geq 1}}_{\mathcal{F}_{\nu}^{0,1}} + \norm{U_{\geq 1}}_{\dot{\mathcal{F}}_{\nu}^{s-1,1}}\norm{\phi}_{\dot{\mathcal{F}}_{\nu}^{1,1}}).
    \end{align*}
    Applying Proposition \ref{p5} to
    \begin{align*}
        \norm{T_{\geq 2}}_{\mathcal{F}_{\nu}^{0,1}} \leq&2\norm{U_{\geq 2}}_{\dot{\mathcal{F}}_{\nu}^{0,1}} + 2\norm{\phi_{\alpha}U_{\geq 1}}_{\dot{\mathcal{F}}_{\nu}^{0,1}},
    \end{align*}
    we obtain
    \begin{align*}
        \norm{T_{\geq 2}}_{\mathcal{F}_{\nu}^{0,1}} &\leq 2\norm{U_{\geq 2}}_{\dot{\mathcal{F}}_{\nu}^{0,1}} + 2(\norm{\phi}_{\dot{\mathcal{F}}_{\nu}^{1,1}}\norm{U_{\geq 1}}_{\mathcal{F}_{\nu}^{0,1}} + \norm{\phi}_{\mathcal{F}_{\nu}^{1,1}}\norm{U_{\geq 1}}_{\dot{\mathcal{F}}_{\nu}^{0,1}}).
    \end{align*}
    In the case of \(0 \leq s < 1\), applying Proposition \ref{p5} to
\begin{align*}
    \norm{T_{\geq 2}}_{\dot{\mathcal{F}}_{\nu}^{s,1}} \leq&\norm{U_{\geq 2}}_{\dot{\mathcal{F}}_{\nu}^{0,1}} + \norm{\phi_{\alpha}U_{\geq 1}}_{\dot{\mathcal{F}}_{\nu}^{0,1}},
\end{align*}
we obtain
\begin{align*}
    \norm{T_{\geq 2}}_{\dot{\mathcal{F}}_{\nu}^{s,1}} &\leq \norm{U_{\geq 2}}_{\dot{\mathcal{F}}_{\nu}^{0,1}} + b(2,s)\biggl(\norm{\phi}_{\dot{\mathcal{F}}_{\nu}^{1,1}}\norm{U_{\geq 1}}_{\mathcal{F}_{\nu}^{0,1}} + \norm{\phi}_{\dot{\mathcal{F}}_{\nu}^{1,1}}\norm{U_{\geq 1}}_{\dot{\mathcal{F}}_{\nu}^{0,1}}\biggr).
\end{align*}
\end{proof}

\subsection{Estimating \texorpdfstring{\(T_{1}(\alpha)\phi_{\alpha}(\alpha)\)}{t1}} \label{estimaingterm3}
We prove the following estimates.
\begin{lemma} \label{estimatingterm3lem}
    For \(s \geq 1\),
    \begin{align*}
        \norm{T_{1}\phi_{\alpha}}_{\dot{\mathcal{F}}_{\nu}^{s,1}} &\leq b(2,s)\normx{\phi}{1}\norm{U_{1}}_{\dot{\mathcal{F}}_{\nu}^{s-1,1}}+b(2,s)\normx{\phi}{s+1}2\norm{U_{1}}_{\dot{\mathcal{F}}_{\nu}^{0,1}}.
    \end{align*}
    For \(0 \leq s <1\),
    \begin{align*}
        \norm{T_{1}\phi_{\alpha}}_{\dot{\mathcal{F}}_{\nu}^{s,1}} &\leq b(2,s)\normx{\phi}{1}\norm{U_{1}}_{\dot{\mathcal{F}}_{\nu}^{0,1}}+b(2,s)\normx{\phi}{s+1}2\norm{U_{1}}_{\dot{\mathcal{F}}_{\nu}^{0,1}}.
    \end{align*}
\end{lemma}
\begin{proof}
    Using Proposition \ref{p5}, we obtain that for \(s \geq 0\),
\begin{align*}
    \norm{T_{1}\phi_{\alpha}}_{\dot{\mathcal{F}}_{\nu}^{s,1}} &\leq b(2,s)\biggl(\norm{T_{1}}_{\dot{\mathcal{F}}_{\nu}^{s,1}}\norm{\phi}_{\dot{\mathcal{F}}_{\nu}^{1,1}} + \norm{\phi}_{\dot{\mathcal{F}}_{\nu}^{s+1,1}}\norm{T_{1}}_{\mathcal{F}_{\nu}^{0,1}}\biggr).
\end{align*}
It remains to estimate \(\norm{T_{1}}_{\dot{\mathcal{F}}_{\nu}^{s,1}}\) and \(\norm{T_{1}}_{\mathcal{F}_{\nu}^{0,1}} \) in terms of \(U_{1}\) and \(U_{\geq 2}\). Recalling that \(T_{1}(\alpha) = \mathcal{M}(U_{1})(\alpha)\), we obtain that for \(s \geq 1\), \(\norm{T_{1}}_{\dot{\mathcal{F}}_{\nu}^{s,1}}=\norm{U_{1}}_{\dot{\mathcal{F}}_{\nu}^{s-1,1}}\). Moreover, \(\norm{T_{1}}_{\mathcal{F}_{\nu}^{0,1}} \leq 2\norm{U_{1}}_{\dot{\mathcal{F}}_{\nu}^{0,1}}\). In the case of \(0 \leq s < 1\), \(\norm{T_{1}}_{\dot{\mathcal{F}}_{\nu}^{s,1}} \leq \norm{U_{1}}_{\dot{\mathcal{F}}_{\nu}^{0,1}}\).
\end{proof}

\section{Estimating \texorpdfstring{\(U_{1}\)}{U1}}
\label{estimationU1}
To estimate the \(\dot{\mathcal{F}}_{\nu}^{s,1}\) and \(\mathcal{F}_{\nu}^{0,1}\) norms of \(U_{1}\), we first estimate the Fourier modes of \(U_{1}\).
\subsection{Estimating Fourier Modes of \texorpdfstring{\(U_{1}\)}{U1}}
\label{estimatingfouriermodesofU1}
For any norm \(\norm{\cdot}\), we can estimate (\ref{U1}) as
\begin{align}
    \norm{U_{1}} \leq \norm{ie^{i\alpha}e^{i\hat{\theta}(0)}\mathfrak{L}(\alpha)}. \label{U1withoutRe}
\end{align}
To estimate the \(\dot{\mathcal{F}}_{\nu}^{s,1}\) and \(\mathcal{F}_{\nu}^{0,1}\) norms of the right hand side of (\ref{U1withoutRe}), we write
\begin{align}
    ie^{i\alpha}e^{i\hat{\theta}(0)}\mathfrak{L}(\alpha) = \sum_{j=1}^{7}\frac{\gamma}{4\pi}\int_{-\pi}^{\pi}E_{j}(\alpha, \beta)d\beta, \nonumber
\end{align}
where
\begin{align}
    E_{1}(\alpha, \beta) = &\frac{-e^{i\beta}(-1+e^{i\beta})(i(-1+e^{i\beta})+\beta(1+e^{i\beta}))}{2(-1+e^{i\beta})^{2}} \nonumber \\
    &\cdot \int_{0}^{1}e^{-i\beta s}\phi(\alpha+\beta(-1+s))ds, \nonumber \\
    E_{2}(\alpha, \beta) = &\frac{i(-1-2i\beta+e^{2i\beta})}{2(-1+e^{i\beta})^{2}}\int_{0}^{1}e^{i\beta s}\phi(\alpha+\beta(-1+s))ds, \nonumber \\
    E_{3}(\alpha, \beta) = &\frac{-(-1+e^{i\beta})\beta e^{i\beta}(-1+e^{i\beta})}{2(-1+e^{i\beta})^{2}} \nonumber \\
    &\cdot\int_{0}^{1}e^{-i\beta s}(-1+s)\phi(\alpha+\beta(-1+s))ds, \nonumber \\
    E_{4}(\alpha, \beta) = &\frac{-(-1+e^{i\beta})\beta(1+e^{i\beta})}{2(-1+e^{i\beta})^{2}}\int_{0}^{1}e^{i\beta s}(-1+s)\phi(\alpha+\beta(-1+s))ds, \nonumber \\
    E_{5}(\alpha, \beta) = &\frac{-(-1+e^{i\beta})i\beta e^{i\beta}(-1+e^{i\beta})}{2(-1+e^{i\beta})^{2}} \nonumber \\
    &\cdot \int_{0}^{1}e^{-i\beta s}(-1+s)\phi'(\alpha+\beta(-1+s))ds, \nonumber \\
    E_{6}(\alpha, \beta) = &\frac{-(-1+e^{i\beta})i(-\beta)(1+e^{i\beta})}{2(-1+e^{i\beta})^{2}} \nonumber \\
    &\cdot \int_{0}^{1}e^{i\beta s}(-1+s)\phi'(\alpha+\beta(-1+s))ds, \nonumber \\
    E_{7}(\alpha, \beta) = &\frac{-(-1+e^{i\beta})i(-1+2e^{i\beta}+e^{2i\beta})}{2(-1+e^{i\beta})^{2}}\phi(\alpha-\beta). \nonumber
\end{align}
First, we calculate the Fourier modes of \(E_{1}(\alpha, \beta)\):
\begin{align*}
    \mathcal{F}(E_{1})(k,\beta) = \frac{-e^{i\beta}(i(-1+e^{i\beta})+\beta(1+e^{i\beta}))}{2(-1+e^{i\beta})}\cdot\int_{0}^{1}e^{-i\beta s}e^{ik\beta(-1+s)}ds \cdot \mathcal{F}(\phi)(k).
\end{align*}
Since
\begin{align}
    &\abs{\frac{\gamma}{4\pi}\int_{-\pi}^{\pi}\frac{-e^{i\beta}}{2(-1+e^{i\beta})}(i(-1+e^{i\beta})+\beta(1+e^{i\beta}))\int_{0}^{1}e^{-i\beta s}e^{ik\beta(-1+s)}dsd\beta} \nonumber \\
    \leq&\frac{\gamma}{4\pi}\biggl(2\pi + \frac{1}{2}\sqrt{1+\frac{\pi^{2}}{4}}\cdot\frac{1}{2}\pi^{2}\biggr), \nonumber
\end{align}
we obtain
\begin{align*}
    \abs{\frac{\gamma}{4\pi}\int_{-\pi}^{\pi}\mathcal{F}(E_{1})(k,\beta)d\beta} \leq \frac{\gamma}{4\pi}\biggl(2\pi+\frac{\pi^{2}}{4}\sqrt{1+\frac{\pi^{2}}{4}}\biggr)\abs{\mathcal{F}(\phi)(k)}. 
\end{align*}
Similarly, we obtain
\begin{align*}
    \abs{\frac{\gamma}{4\pi}\int_{-\pi}^{\pi}\mathcal{F}(E_{2})(k, \beta)d\beta} \leq& \biggl(\frac{\gamma}{4\pi}\biggl(\frac{1}{2}\sqrt{1+\frac{\pi^{2}}{4}}\cdot 2\pi+\pi^{2}\biggr) \\
    &+\frac{\gamma}{4\pi}\biggl(2\cdot \frac{\pi}{2}\cdot\frac{1}{2}\sqrt{1+\frac{\pi^{2}}{4}}\cdot 2\pi +\pi^{2}\biggr)\biggr) \abs{\mathcal{F}(\phi)(k)}, \\
    \abs{\frac{\gamma}{4\pi}\int_{-\pi}^{\pi}\mathcal{F}(E_{3})(k,\beta)d\beta} \leq& \frac{\gamma}{4\pi}\cdot\frac{\pi^{2}}{2}\abs{\mathcal{F}(\phi)(k)}, \\
    \abs{\frac{\gamma}{4\pi}\int_{-\pi}^{\pi}\mathcal{F}(E_{4})(k,\beta)d\beta} \leq& \frac{\gamma}{4\pi}\biggl(\frac{1}{2}\sqrt{1+\frac{\pi^{2}}{4}}\cdot\pi^{2}+2\pi\biggr)\abs{\mathcal{F}(\phi)(k)}, \\
    \abs{\frac{\gamma}{4\pi}\int_{-\pi}^{\pi}\mathcal{F}(E_{5})(k, \beta)d\beta} \leq& \frac{\gamma}{4\pi}\cdot\frac{\pi^{2}}{2}\abs{k}\cdot \abs{\mathcal{F}(\phi)(k)}, \\
    \abs{\frac{\gamma}{4\pi}\int_{-\pi}^{\pi}\mathcal{F}(E_{6})(k, \beta)d\beta} \leq& \frac{\gamma}{4\pi}\biggl(\frac{1}{2}\sqrt{1+\frac{\pi^{2}}{4}}\cdot \pi^{2} + 2\pi\biggr)\abs{k}\abs{\mathcal{F}(\phi)(k)}, \\
    \abs{\frac{\gamma}{4\pi}\int_{-\pi}^{\pi}\mathcal{F}(E_{7})(k,\beta)d\beta} \leq& \frac{\gamma}{4\pi}\biggl(\frac{1}{2}\sqrt{1+\frac{\pi^{2}}{4}}\cdot\frac{3}{2}\cdot 2\pi +\frac{1}{2}\cdot 4\cdot 5\biggr)\abs{\mathcal{F}(\phi)(k)}.
\end{align*}

\subsection{Estimating \texorpdfstring{\(\normz{U_{1}}\)}{U1}}
\label{estNormzU1}
In Section \ref{estimatingfouriermodesofU1}, we observed that
\begin{align*}
    \normz{U_{1}} \leq \sum_{j=1}^{7}\normz{\frac{\gamma}{4\pi}\int_{-\pi}^{\pi}E_{j}(\alpha, \beta)d\beta}.
\end{align*}
Since
\begin{align*}
    \normz{\frac{\gamma}{4\pi}\int_{-\pi}^{\pi}E_{1}(\alpha, \beta)d\beta} \leq&\frac{\gamma}{4\pi}\biggl(2\pi+\frac{\pi^{2}}{4}\sqrt{1+\frac{\pi^{2}}{4}}\biggr)\normz{\phi}, \\
    \normz{\frac{\gamma}{4\pi}\int_{-\pi}^{\pi}E_{2}(\alpha, \beta)d\beta} \leq&\biggl(\frac{\gamma}{4\pi}\biggl(\frac{1}{2}\sqrt{1+\frac{\pi^{2}}{4}}\cdot 2\pi+\pi^{2}\biggr) \\
    &+\frac{\gamma}{4\pi}\biggl(2\cdot \frac{\pi}{2}\cdot\frac{1}{2}\sqrt{1+\frac{\pi^{2}}{4}}\cdot 2\pi +\pi^{2}\biggr)\biggr)\normz{\phi}, \\
    \normz{\frac{\gamma}{4\pi}\int_{-\pi}^{\pi}E_{3}(\alpha, \beta)d\beta} \leq&\frac{\gamma}{4\pi}\cdot\frac{\pi^{2}}{2}\normz{\phi}, \\
    \normz{\frac{\gamma}{4\pi}\int_{-\pi}^{\pi}E_{4}(\alpha, \beta)d\beta} \leq&\frac{\gamma}{4\pi}\biggl(\frac{1}{2}\sqrt{1+\frac{\pi^{2}}{4}}\cdot\pi^{2}+2\pi\biggr)\normz{\phi}, \\
    \normz{\frac{\gamma}{4\pi}\int_{-\pi}^{\pi}E_{5}(\alpha, \beta)d\beta} \leq&\frac{\gamma}{4\pi}\cdot\frac{\pi^{2}}{2}\norm{\phi}_{\mathcal{F}_{\nu}^{1,1}}, \\
    \normz{\frac{\gamma}{4\pi}\int_{-\pi}^{\pi}E_{6}(\alpha, \beta)d\beta} \leq&\frac{\gamma}{4\pi}\biggl(\frac{1}{2}\sqrt{1+\frac{\pi^{2}}{4}}\cdot \pi^{2} + 2\pi\biggr)\norm{\phi}_{\mathcal{F}_{\nu}^{1,1}}, \\
    \normz{\frac{\gamma}{4\pi}\int_{-\pi}^{\pi}E_{7}(\alpha, \beta)d\beta} \leq&\frac{\gamma}{4\pi}\biggl(\frac{1}{2}\sqrt{1+\frac{\pi^{2}}{4}}\cdot\frac{3}{2}\cdot 2\pi +\frac{1}{2}\cdot 4\cdot 5\biggr)\normz{\phi},
\end{align*}
we obtain
\begin{align*}
    \normz{U_{1}} \leq H_{3}\normz{\phi}+H_{4}\norm{\phi}_{\mathcal{F}_{\nu}^{1,1}},
\end{align*}
where \(H_{3}\) and \(H_{4}\) are constants.

\subsection{Estimating \texorpdfstring{\(\norms{U_{1}}\)}{U1}}
In Section \ref{estimatingfouriermodesofU1}, we observed that
\begin{align*}
    \norms{U_{1}} \leq \sum_{j=1}^{7}\norms{\frac{\gamma}{4\pi}\int_{-\pi}^{\pi}E_{j}(\alpha, \beta)d\beta}.
\end{align*}
Since
\begin{align*}
    \norms{\frac{\gamma}{4\pi}\int_{-\pi}^{\pi}E_{1}(\alpha, \beta)d\beta} \leq&\frac{\gamma}{4\pi}\biggl(2\pi+\frac{\pi^{2}}{4}\sqrt{1+\frac{\pi^{2}}{4}}\biggr)\norms{\phi}, \\
    \norms{\frac{\gamma}{4\pi}\int_{-\pi}^{\pi}E_{2}(\alpha, \beta)d\beta} \leq& \biggl(\frac{\gamma}{4\pi}\biggl(\frac{1}{2}\sqrt{1+\frac{\pi^{2}}{4}}\cdot 2\pi+\pi^{2}\biggr) \\
    &+\frac{\gamma}{4\pi}\biggl(2\cdot \frac{\pi}{2}\cdot\frac{1}{2}\sqrt{1+\frac{\pi^{2}}{4}}\cdot 2\pi +\pi^{2}\biggr)\biggr)\norms{\phi}, \\
    \norms{\frac{\gamma}{4\pi}\int_{-\pi}^{\pi}E_{3}(\alpha, \beta)d\beta} \leq& \frac{\gamma}{4\pi}\cdot\frac{\pi^{2}}{2}\norms{\phi}, \\
    \norms{\frac{\gamma}{4\pi}\int_{-\pi}^{\pi}E_{4}(\alpha, \beta)d\beta} \leq& \frac{\gamma}{4\pi}\biggl(\frac{1}{2}\sqrt{1+\frac{\pi^{2}}{4}}\cdot\pi^{2}+2\pi\biggr) \norms{\phi}, \\
    \norms{\frac{\gamma}{4\pi}\int_{-\pi}^{\pi}E_{5}(\alpha, \beta)d\beta} \leq& \frac{\gamma}{4\pi}\cdot\frac{\pi^{2}}{2}\normx{\phi}{s+1}, \\
    \norms{\frac{\gamma}{4\pi}\int_{-\pi}^{\pi}E_{6}(\alpha, \beta)d\beta} \leq& \frac{\gamma}{4\pi}\biggl(\frac{1}{2}\sqrt{1+\frac{\pi^{2}}{4}}\cdot \pi^{2} + 2\pi\biggr)\normx{\phi}{s+1}, \\
    \norms{\frac{\gamma}{4\pi}\int_{-\pi}^{\pi}E_{7}(\alpha, \beta)d\beta} \leq& \frac{\gamma}{4\pi}\biggl(\frac{1}{2}\sqrt{1+\frac{\pi^{2}}{4}}\cdot\frac{3}{2}\cdot 2\pi +\frac{1}{2}\cdot 4\cdot 5\biggr)\norms{\phi},
\end{align*}
we obtain
\begin{align*}
    \norms{U_{1}} \leq H_{1}\norms{\phi} + H_{2}\normx{\phi}{s+1},
\end{align*}
where \(H_{1}\) and \(H_{2}\) are constants.

\section{Estimating \texorpdfstring{\(U_{\geq 2}\)}{U2}}
\label{estimationU2}
For any norm \(\norm{\cdot}\), we can estimate (\ref{U2}) as
\begin{align}
    \norm{U_{\geq 2}} &\leq \norm{ie^{i\alpha}e^{i\hat{\theta}(0)}\biggl(e^{i\phi(\alpha)}(\mathfrak{L}(\alpha) + \mathfrak{N}(\alpha))-\mathfrak{L}(\alpha)\biggr)}. \label{U2withoutRe}
\end{align}
To estimate the \(\dot{\mathcal{F}}_{\nu}^{s,1}\) and \(\mathcal{F}_{\nu}^{0,1}\) norms of the right hand side of (\ref{U2withoutRe}), we first note that
\begin{align}
    ie^{i\alpha}e^{i\hat{\theta}(0)}\biggl(e^{i\phi(\alpha)}(\mathfrak{L}(\alpha) + \mathfrak{N}(\alpha))-\mathfrak{L}(\alpha)\biggr) = \frac{\gamma}{4\pi}\int_{-\pi}^{\pi}B(\alpha, \beta)d\beta, \label{U2sum}
\end{align}
where \(B(\alpha, \beta) = \sum_{j=1}^{8}\widetilde{B_{j}}(\alpha, \beta) + \widetilde{B_{13}}(\alpha, \beta)\), in which
\begin{align*}
    \widetilde{B_{1}}(\alpha, \beta) = -&\sum_{\substack{j_{1}+j_{2}+n \geq 1 \\ j_{3}=1}} \frac{-i\beta e^{2i\beta}}{1-e^{i\beta}}\frac{(-1)^{j_{1}+j_{3}}i^{j_{1}+j_{2}+j_{3}}}{2j_{1}!j_{2}!j_{3}!}\phi(\alpha - \beta)^{j_{1}}\phi(\alpha)^{j_{2}} \\
    &\cdot\int_{0}^{1}e^{-i\beta s}\phi(\alpha + \beta(-1+s))^{j_{3}}(-1+s)ds \\
    &\cdot\biggl(\sum_{m=1}^{\infty}\frac{i\beta}{1-e^{i\beta}}\int_{0}^{1}e^{-i(s-1)\beta}\frac{(-i\phi(\alpha + (s-1)\beta))^{m}}{m!}ds\biggr)^{n} \\
    -&\sum_{\substack{j_{1}+j_{2}+j_{3} \geq 1 \\ n=1}} \frac{-i\beta e^{2i\beta}}{1-e^{i\beta}}\frac{(-1)^{j_{1}+j_{3}}i^{j_{1}+j_{2}+j_{3}}}{2j_{1}!j_{2}!j_{3}!}\phi(\alpha - \beta)^{j_{1}}\phi(\alpha)^{j_{2}} \\
    &\cdot\int_{0}^{1}e^{-i\beta s}\phi(\alpha + \beta(-1+s))^{j_{3}}(-1+s)ds \\
    &\cdot\biggl(\sum_{m=1}^{\infty}\frac{i\beta}{1-e^{i\beta}}\int_{0}^{1}e^{-i(s-1)\beta}\frac{(-i\phi(\alpha + (s-1)\beta))^{m}}{m!}ds\biggr)^{n} \\
    - &\frac{-i\beta e^{2i\beta}}{1-e^{i\beta}}\cdot\frac{1}{2}\int_{0}^{1}e^{-i\beta s}(-1+s)ds \\
    &\cdot\sum_{m=2}^{\infty}\frac{i\beta}{1-e^{i\beta}}\int_{0}^{1}e^{-i(s-1)\beta}\frac{(-i\phi(\alpha+(s-1)\beta))^{m}}{m!}ds \\
    \widetilde{B_{2}}(\alpha, \beta) =-&\frac{1}{2}\sum_{j_{1}+j_{2}+j_{3}+n \geq 1} \frac{-i\beta e^{2i\beta}}{1-e^{i\beta}}\frac{(-1)^{j_{1}+j_{3}}i^{j_{1}+j_{2}+j_{3}}}{j_{1}!j_{2}!j_{3}!}\phi(\alpha - \beta)^{j_{1}}\phi(\alpha)^{j_{2}} \\
    &\cdot\int_{0}^{1}e^{-i\beta s}\phi(\alpha + \beta(-1+s))^{j_{3}}(-1+s)\phi'(\alpha + \beta(-1+s))ds \\
    &\cdot\biggl(\sum_{m=1}^{\infty}\frac{i\beta}{1-e^{i\beta}}\int_{0}^{1}e^{-i(s-1)\beta}\frac{(-i\phi(\alpha+(s-1)\beta))^{m}}{m!}ds\biggr)^{n} \\
    \widetilde{B_{3}}(\alpha, \beta) = \frac{1}{2}&\sum_{\substack{j_{1}=j_{2}=0 \\ j_{3}+j_{4}+n\geq 2}} (n+1)\frac{-\beta^{2}e^{-i\beta}}{(1-e^{-i\beta})^{2}}\frac{i^{j_{1}+j_{2}+j_{3}+j_{4}}(-1)^{j_{3}}}{j_{1}!j_{2}!j_{3}!j_{4}!}\phi(\alpha)^{j_{1}}\phi(\alpha - \beta)^{j_{2}} \\
    &\cdot\int_{0}^{1}e^{-i\beta s}\phi(\alpha + \beta(-1+s))^{j_{3}}ds \\
    &\cdot \int_{0}^{1}e^{i\beta s}\phi(\alpha + \beta(-1+s))^{j_{4}}(-1+s)ds \\
    &\cdot\biggl(\sum_{m=1}^{\infty}\frac{-i^{m}}{m!}e^{-i\alpha}\frac{i\beta}{1-e^{-i\beta}}\int_{0}^{1}e^{i(\alpha + (s-1)\beta)}\phi(\alpha + (s-1)\beta)^{m}ds\biggr)^{n} \\
    +&\frac{1}{2}\sum_{\substack{j_{1}+j_{2}\geq 1 \\ j_{3}+j_{4}+n \geq 1}} (n+1)\frac{-\beta^{2}e^{-i\beta}}{(1-e^{-i\beta})^{2}}\frac{i^{j_{1}+j_{2}+j_{3}+j_{4}}(-1)^{j_{3}}}{j_{1}!j_{2}!j_{3}!j_{4}!}\phi(\alpha)^{j_{1}}\phi(\alpha - \beta)^{j_{2}} \\
    &\cdot\int_{0}^{1}e^{-i\beta s}\phi(\alpha + \beta(-1+s))^{j_{3}}ds \\
    &\cdot \int_{0}^{1}e^{i\beta s}\phi(\alpha + \beta(-1+s))^{j_{4}}(-1+s)ds \\
    &\cdot\biggl(\sum_{m=1}^{\infty}\frac{-i^{m}}{m!}e^{-i\alpha}\frac{i\beta}{1-e^{-i\beta}}\int_{0}^{1}e^{i(\alpha + (s-1)\beta)}\phi(\alpha + (s-1)\beta)^{m}ds\biggr)^{n}
\end{align*}
\begin{align*}
    +&\frac{1}{2}\cdot 2 \cdot \frac{-\beta^{2}e^{-i\beta}}{(1-e^{-i\beta})^{2}}\int_{0}^{1}e^{-i\beta s}ds\int_{0}^{1}e^{i\beta s}(-1+s)ds \sum_{m=2}^{\infty}\frac{-i^{m}}{m!} \\
    &\cdot e^{-i\alpha}\frac{i\beta}{1-e^{-i\beta}}\int_{0}^{1}e^{i(\alpha + (s-1)\beta)}\phi(\alpha + (s-1)\beta)^{m}ds \\
    \widetilde{B_{4}}(\alpha, \beta) = \frac{1}{2}&\sum_{j_{1}+j_{2}+j_{3}+j_{4}+n \geq 1} (n+1)\frac{-\beta^{2}e^{-i\beta}}{(1-e^{-i\beta})^{2}}\frac{i^{j_{1}+j_{2}+j_{3}+j_{4}}(-1)^{j_{3}}}{j_{1}!j_{2}!j_{3}!j_{4}!} \\
    &\cdot\phi(\alpha)^{j_{1}}\phi(\alpha - \beta)^{j_{2}}\int_{0}^{1}e^{-i\beta s}\phi(\alpha + \beta(-1+s))^{j_{3}}ds \\
    &\cdot\int_{0}^{1}e^{i\beta s}\phi(\alpha + \beta(-1+s))^{j_{4}}(-1+s)\phi'(\alpha + \beta(-1+s))ds \\
    &\cdot\biggl(\sum_{m=1}^{\infty}\frac{-i^{m}}{m!}e^{-i\alpha}\frac{i\beta}{1-e^{-i\beta}}\int_{0}^{1}e^{i(\alpha + (s-1)\beta)}\phi(\alpha + (s-1)\beta)^{m}ds\biggr)^{n} \\
    \widetilde{B_{5}}(\alpha, \beta) = \frac{1}{2}&\sum_{\substack{j_{1}=j_{2}=0 \\ j_{3}+n \geq 2}} \frac{i\beta}{1-e^{-i\beta}}\frac{i^{j_{1}+j_{2}+j_{3}}(-1)^{j_{3}}}{j_{1}!j_{2}!j_{3}!}\phi(\alpha)^{j_{1}}\phi(\alpha - \beta)^{j_{2}} \\
    &\cdot\int_{0}^{1}e^{-i\beta s}\phi(\alpha + \beta(-1+s))^{j_{3}}(-1+s)ds \\
    &\cdot\biggl(\sum_{m=1}^{\infty}\frac{-i\beta}{1-e^{-i\beta}}\int_{0}^{1}e^{i(s-1)\beta}\frac{(i\phi(\alpha + (s-1)\beta))^{m}}{m!}ds\biggr)^{n} \\
    + &\frac{1}{2}\sum_{\substack{j_{1}+j_{2}\geq 1 \\ j_{3}+n \geq 1}} \frac{i\beta}{1-e^{-i\beta}}\frac{i^{j_{1}+j_{2}+j_{3}}(-1)^{j_{3}}}{j_{1}!j_{2}!j_{3}!}\phi(\alpha)^{j_{1}}\phi(\alpha - \beta)^{j_{2}} \\
    &\cdot\int_{0}^{1}e^{-i\beta s}\phi(\alpha + \beta(-1+s))^{j_{3}}(-1+s)ds \\
    &\cdot\biggl(\sum_{m=1}^{\infty}\frac{-i\beta}{1-e^{-i\beta}}\int_{0}^{1}e^{i(s-1)\beta}\frac{(i\phi(\alpha + (s-1)\beta))^{m}}{m!}ds\biggr)^{n} \\
    + &\frac{1}{2}\cdot\frac{i\beta}{1-e^{-i\beta}}\int_{0}^{1}e^{-i\beta s}(-1+s)ds \\
    &\cdot\sum_{m=2}^{\infty}\frac{-i\beta}{1-e^{-i\beta}}\int_{0}^{1}e^{i(s-1)\beta}\frac{(i\phi(\alpha + (s-1)\beta))^{m}}{m!}ds \\
    \widetilde{B_{6}}(\alpha, \beta) = \frac{1}{2}&\sum_{\substack{j_{1}=j_{2}=0 \\ j_{3}+n \geq 2}} \frac{i\beta}{1-e^{-i\beta}}\frac{i^{j_{1}+j_{2}+j_{3}}(-1)^{j_{2}}}{j_{1}!j_{2}!j_{3}!}\phi(\alpha)^{j_{1}}\phi(\alpha - \beta)^{j_{2}} \\
    &\cdot\int_{0}^{1}e^{i\beta s}\phi(\alpha + \beta(-1+s))^{j_{3}}(-1+s)ds \\
    &\cdot\biggl(\sum_{m=1}^{\infty}\frac{-i\beta}{1-e^{-i\beta}}\int_{0}^{1}e^{i(s-1)\beta}\frac{(i\phi(\alpha + (s-1)\beta))^{m}}{m!}ds\biggr)^{n} \\
    +& \frac{1}{2}\sum_{\substack{j_{1}+j_{2}\geq 1 \\ j_{3}+n \geq 1}} \frac{i\beta}{1-e^{-i\beta}}\frac{i^{j_{1}+j_{2}+j_{3}}(-1)^{j_{2}}}{j_{1}!j_{2}!j_{3}!}\phi(\alpha)^{j_{1}}\phi(\alpha - \beta)^{j_{2}} \\
    &\cdot\int_{0}^{1}e^{i\beta s}\phi(\alpha + \beta(-1+s))^{j_{3}}(-1+s)ds \\
    &\cdot\biggl(\sum_{m=1}^{\infty}\frac{-i\beta}{1-e^{-i\beta}}\int_{0}^{1}e^{i(s-1)\beta}\frac{(i\phi(\alpha + (s-1)\beta))^{m}}{m!}ds\biggr)^{n}
\end{align*}
\begin{align*}
    +&\frac{1}{2}\cdot\frac{i\beta}{1-e^{-i\beta}}\int_{0}^{1}e^{i\beta s}(-1+s)ds \\
    &\cdot\sum_{m=2}^{\infty}\frac{-i\beta}{1-e^{-i\beta}}\int_{0}^{1}e^{i(s-1)\beta}\frac{(i\phi(\alpha + (s-1)\beta))^{m}}{m!}ds \\
    \widetilde{B_{7}}(\alpha, \beta) =& \frac{1}{2}\sum_{j_{1}+j_{2}+j_{3}+n \geq 1} \frac{i\beta}{1-e^{-i\beta}}\frac{i^{j_{1}+j_{2}+j_{3}}(-1)^{j_{3}}}{j_{1}!j_{2}!j_{3}!}\phi(\alpha)^{j_{1}}\phi(\alpha - \beta)^{j_{2}} \\
    &\cdot\int_{0}^{1}e^{-i\beta s}\phi(\alpha + \beta(-1+s))^{j_{3}}(-1+s)\phi'(\alpha + \beta(-1+s))ds \\
    &\cdot\biggl(\sum_{m=1}^{\infty}\frac{-i\beta}{1-e^{-i\beta}}\int_{0}^{1}e^{i(s-1)\beta}\frac{(i\phi(\alpha + (s-1)\beta))^{m}}{m!}ds\biggr)^{n} \\
    \widetilde{B_{8}}(\alpha, \beta) =& \frac{1}{2}\sum_{j_{1}+j_{2}+j_{3}+n \geq 1} \frac{i\beta}{1-e^{-i\beta}}\frac{i^{j_{1}+j_{2}+j_{3}}(-1)^{j_{2}}}{j_{1}!j_{2}!j_{3}!}\phi(\alpha)^{j_{1}}\phi(\alpha-\beta)^{j_{2}} \\
    &\cdot\int_{0}^{1}e^{i\beta s}\phi(\alpha + \beta(-1+s))^{j_{3}}(-1+s)\phi'(\alpha + \beta(-1+s))ds \\
    &\cdot\biggl(\sum_{m=1}^{\infty}\frac{-i\beta}{1-e^{-i\beta}}\int_{0}^{1}e^{i(s-1)\beta}\frac{(i\phi(\alpha + (s-1)\beta))^{m}}{m!}ds\biggr)^{n} \\
    \widetilde{B_{13}}(\alpha, \beta) =& \sum_{j_{1}+j_{2} \geq 2} \frac{i^{j_{1}+j_{2}}}{j_{1}!j_{2}!}\phi(\alpha)^{j_{1}}\phi(\alpha - \beta)^{j_{2}}\biggl((-1)^{j_{2}}\frac{e^{i\beta}(1+e^{i\beta})}{2(-1+e^{i\beta})}-\frac{1}{2}\biggr).
\end{align*}

\subsection{Estimating Fourier Modes of \texorpdfstring{\(U_{\geq 2}\)}{U2}}
\label{efmodesU2}
In our calculations, we adopt the notational convention that any product \(\prod\) in which the upper bound is strictly less than the lower bound is defined to be \(1\). To calculate the Fourier modes of \(U_{\geq 2}\), we frequently use the identity
\begin{align} \label{convformula}
    \mathcal{F}(g_{1}g_{2}\cdots g_{n})(k_{1}) = \sum_{k_{2}, \dots, k_{n} \in \mathbb{Z}}\biggl(\prod_{d=1}^{n-1}\mathcal{F}(g_{d})(k_{d}-k_{d+1})\biggr)\mathcal{F}(g_{n})(k_{n}).
\end{align}
The following estimate is used frequently.
\begin{lemma} \label{inbound}
For \(n \geq 0\), let
\begin{align*}
    &I_{n}(k_{1}, k_{j_{1}+1}, k_{j_{1}+j_{2}+1}, \dots, k_{j_{1}+j_{2}+n}, k_{j_{1}+j_{2}+n+1}, \beta) \\
    =& \prod_{d=1}^{n}\biggl(\frac{i\beta e^{i\beta}}{1-e^{i\beta}}\int_{0}^{1}e^{-is\beta}e^{i(s-1)\beta(k_{j_{1}+j_{2}+d}-k_{j_{1}+j_{2}+d+1})}ds\biggr) \cdot e^{-i\beta(k_{1}-k_{j_{1}+1})} \\
    &\cdot\int_{0}^{1}e^{-i\beta s}e^{i\beta(-1+s)k_{j_{1}+j_{2}+n+1}}(-1+s)ds.
\end{align*}
Then
    \begin{align*}
        \abs{\frac{\gamma}{4\pi}\int_{-\pi}^{\pi}I_{n}(k_{1}, k_{j_{1}+1}, k_{j_{1}+j_{2}+1}, \dots, k_{j_{1}+j_{2}+n}, k_{j_{1}+j_{2}+n+1}, \beta) \cdot \frac{-i\beta e^{2i\beta}}{1-e^{i\beta}} d\beta} \leq C_{n},
    \end{align*}
    where
    \begin{align*}
        C_{n} = \frac{\gamma}{4\pi}\biggl((n+1) \cdot \biggl(\frac{\pi}{2}\biggr)^{n}\cdot\frac{1}{2}\sqrt{1+\frac{\pi^{2}}{4}}\cdot \pi^{2} + 2\pi\biggr).
    \end{align*}
\end{lemma}
\begin{proof}
    We note that
    \begin{align*}
        \int_{0}^{1}e^{-is\beta}e^{i(s-1)\beta k}ds =
        \begin{cases}
            \frac{i(e^{-i\beta} - e^{-i\beta k})}{\beta(1-k)} & \mbox{if \(k \neq 1\),} \\
            e^{-i\beta} & \mbox{if \(k = 1\).}
        \end{cases}
    \end{align*}
    In the case of \(n \geq 1\), suppose that \(0 \leq l \leq n\) and \(l\) elements of \(\{k_{j_{1}+j_{2}+d}-k_{j_{1}+j_{2}+d+1}\}_{d=1}^{n}\) satisfy \(k_{j_{1}+j_{2}+d}-k_{j_{1}+j_{2}+d+1} = 1\). Reordering the subscripts such that \(k_{j_{1}+j_{2}+d} - k_{j_{1}+j_{2}+d+1} \neq 1\) for \(d = 1, \dots, n-l\), we obtain
    \begin{align*}
        I_{n} = &e^{-i\beta(k_{1}-k_{j_{1}+1})}\prod_{d=1}^{n-l}\frac{-(1-e^{-i\beta(-1+k_{j_{1}+j_{2}+d}-k_{j_{1}+j_{2}+d+1})})}{(1-e^{i\beta})(1-k_{j_{1}+j_{2}+d} + k_{j_{1}+j_{2}+d+1})}\biggl(\frac{i\beta}{1-e^{i\beta}}\biggr)^{l} \\
        &\cdot\int_{0}^{1}e^{-i\beta s}e^{i\beta(-1+s)k_{j_{1}+j_{2}+n+1}}(-1+s)ds.
    \end{align*}
    If \(k_{j_{1}+j_{2}+d} - k_{j_{1}+j_{2}+d+1} > 1\), then
    \begin{align*}
        \frac{-(1-e^{-i\beta(-1+k_{j_{1}+j_{2}+d}-k_{j_{1}+j_{2}+d+1})})}{1-e^{i\beta}} = e^{-i\beta}\sum_{r_{d}=0}^{-2+k_{j_{1}+j_{2}+d}-k_{j_{1}+j_{2}+d+1}}(e^{-i\beta})^{r_{d}}.
    \end{align*}
    If \(k_{j_{1}+j_{2}+d}-k_{j_{1}+j_{2}+d+1} < 1\), then
    \begin{align*}
        \frac{-(1-e^{-i\beta(-1+k_{j_{1}+j_{2}+d}-k_{j_{1}+j_{2}+d+1})})}{1-e^{i\beta}} = -\sum_{r_{d}=0}^{-(k_{j_{1}+j_{2}+d}-k_{j_{1}+j_{2}+d+1})}(e^{i\beta})^{r_{d}}.
    \end{align*}
    Without loss of generality, suppose that \(k_{j_{1}+j_{2}+d} - k_{j_{1}+j_{2}+d+1} < 1\) only for \(d=w, \dots, n-l\). Then
    \begin{align*}
        &\prod_{d=1}^{n-l}\frac{-(1-e^{-i\beta(-1+k_{j_{1}+j_{2}+d}-k_{j_{1}+j_{2}+d+1})})}{1-e^{i\beta}} \\
        =& e^{-i\beta}\sum_{r_{1}=0}^{-2+k_{j_{1}+j_{2}+1}-k_{j_{1}+j_{2}+2}}(e^{-i\beta})^{r_{1}}\cdot \cdots \cdot e^{-i\beta}\sum_{r_{w-1}=0}^{-2+k_{j_{1}+j_{2}+w-1}-k_{j_{1}+j_{2}+w}}(e^{-i\beta})^{r_{w-1}} \\
        &\cdot(-1)\sum_{r_{w}=0}^{-(k_{j_{1}+j_{2}+w}-k_{j_{1}+j_{2}+w+1})}(e^{i\beta})^{r_{w}} \cdot \cdots \cdot (-1)\sum_{r_{n-l} = 0}^{-(k_{j_{1}+j_{2}+n-l} - k_{j_{1}+j_{2}+n-l+1})}(e^{i\beta})^{r_{n-l}} \\
        =&\sum_{r_{1}=0}^{-2+k_{j_{1}+j_{2}+1}-k_{j_{1}+j_{2}+2}} \cdots \sum_{r_{w-1}=0}^{-2+k_{j_{1}+j_{2}+w-1}-k_{j_{1}+j_{2}+w}} \\
        &\sum_{r_{w}=0}^{-(k_{j_{1}+j_{2}+w}-k_{j_{1}+j_{2}+w+1})} \cdots \sum_{r_{n-l} = 0}^{-(k_{j_{1}+j_{2}+n-l}-k_{j_{1}+j_{2}+n-l+1})} \\
        &(e^{-i\beta})^{w-1}(-1)^{n-l-w+1}(e^{-i\beta})^{r_{1}+\cdots+r_{w-1}}(e^{i\beta})^{r_{w}+\cdots+r_{n-l}}.
    \end{align*}
    Hence,
    \begin{align*}
        I_{n} = &\prod_{d=1}^{n-l}\frac{1}{1-k_{j_{1}+j_{2}+d}+k_{j_{1}+j_{2}+d+1}} \cdot \\
        &\sum_{r_{1}=0}^{-2+k_{j_{1}+j_{2}+1}-k_{j_{1}+j_{2}+2}} \cdots \sum_{r_{w-1}=0}^{-2+k_{j_{1}+j_{2}+w-1}-k_{j_{1}+j_{2}+w}} \\
        &\sum_{r_{w}=0}^{-(k_{j_{1}+j_{2}+w}-k_{j_{1}+j_{2}+w+1})} \cdots \sum_{r_{n-l} = 0}^{-(k_{j_{1}+j_{2}+n-l}-k_{j_{1}+j_{2}+n-l+1})} \\
        &(-1)^{n-l-w+1}(e^{-i\beta})^{w-1}(e^{-i\beta})^{r_{1}+\cdots+r_{w-1}}(e^{i\beta})^{r_{w}+\cdots+r_{n-l}} \\
        &\cdot\biggl(\frac{i\beta}{1-e^{i\beta}}\biggr)^{l}\cdot\int_{0}^{1}e^{-i\beta s}e^{i\beta(-1+s)k_{j_{1}+j_{2}+n+1}}(-1+s)ds \cdot e^{-i\beta(k_{1}-k_{j_{1}+1})}.
    \end{align*}
    Let
    \begin{align*}
        J_{n} = &\frac{\gamma}{4\pi}\int_{-\pi}^{\pi}(e^{-i\beta})^{w-1+r_{1}+\cdots+r_{w-1}-(r_{w}+\cdots+r_{n-l})}\biggl(\frac{i\beta}{1-e^{i\beta}}\biggr)^{l}e^{-i\beta(k_{1}-k_{j_{1}+1})} \\
        &\cdot\int_{0}^{1}e^{-i\beta s}e^{i\beta(-1+s)k_{j_{1}+j_{2}+n+1}}(-1+s)ds\cdot \frac{-i\beta e^{2i\beta}}{1-e^{i\beta}} d\beta.
    \end{align*}
    Using that for all \(l \geq 0\),
    \begin{align*}
        \abs{\biggl(\frac{i\beta}{1-e^{i\beta}}\biggr)^{l}-1} \leq \abs{\beta}\cdot l\cdot\biggl(\frac{\pi}{2}\biggr)^{l-1}\cdot\frac{1}{2}\sqrt{1+\frac{\pi^{2}}{4}},
    \end{align*}
    we obtain
    \begin{align*}
        \abs{J_{n}} &\leq \frac{\gamma}{4\pi}\biggl(\int_{-\pi}^{\pi}\abs{\biggl(\frac{i\beta}{1-e^{i\beta}}\biggr)^{l+1}-1}d\beta + 2\pi\biggr) \\
        &\leq \frac{\gamma}{4\pi}\biggl((l+1) \cdot \biggl(\frac{\pi}{2}\biggr)^{l}\cdot\frac{1}{2}\sqrt{1+\frac{\pi^{2}}{4}}\cdot \pi^{2} + 2\pi\biggr) \\
        &\leq \frac{\gamma}{4\pi}\biggl((n+1) \cdot \biggl(\frac{\pi}{2}\biggr)^{n}\cdot\frac{1}{2}\sqrt{1+\frac{\pi^{2}}{4}}\cdot \pi^{2} + 2\pi\biggr) \\
        &=C_{n}.
    \end{align*}
    Thus,
    \begin{align*}
        &\abs{\frac{\gamma}{4\pi}\int_{-\pi}^{\pi}I_{n}(k_{1}, k_{j_{1}+1}, k_{j_{1}+j_{2}+1}, \dots, k_{j_{1}+j_{2}+n}, k_{j_{1}+j_{2}+n+1}, \beta) \cdot \frac{-i\beta e^{2i\beta}}{1-e^{i\beta}} d\beta} \\
        \leq& \prod_{d=1}^{n-l}\frac{1}{\abs{1-k_{j_{1}+j_{2}+d}+k_{j_{1}+j_{2}+d+1}}}\cdot\sum_{r_{1}=0}^{-2+k_{j_{1}+j_{2}+1}-k_{j_{1}+j_{2}+2}} \cdots \sum_{r_{w-1}=0}^{-2+k_{j_{1}+j_{2}+w-1}-k_{j_{1}+j_{2}+w}} \\
        &\cdot\sum_{r_{w}=0}^{-(k_{j_{1}+j_{2}+w}-k_{j_{1}+j_{2}+w+1})} \cdots \sum_{r_{n-l} = 0}^{-(k_{j_{1}+j_{2}+n-l}-k_{j_{1}+j_{2}+n-l+1})} \abs{J_{n}} \\
        \leq& C_{n}.
    \end{align*}
    In the case of \(n=0\),
    \begin{align*}
        &\abs{\frac{\gamma}{4\pi}\int_{-\pi}^{\pi}I_{0}(k_{1},k_{j_{1}+1}, k_{j_{1}+j_{2}+1}, \beta)\cdot\frac{-i\beta e^{2i\beta}}{1-e^{i\beta}}d\beta} \\
        \leq& \frac{\gamma}{4\pi}\biggl(\int_{-\pi}^{\pi}\abs{\beta}\cdot\frac{1}{2}\sqrt{1+\frac{\pi^{2}}{4}}d\beta +2\pi\biggr) \\
        =& C_{0}.
    \end{align*}
\end{proof}
Let us calculate the Fourier modes of \(\widetilde{B_{1}}(\alpha, \beta)\). We can write \(\widetilde{B_{1}} = \sum_{j=1}^{3}\widetilde{B_{1,j}}\), where
\begin{align*}
    \widetilde{B_{1,1}}(\alpha, \beta) = -&\sum_{\substack{j_{1}+j_{2}+n \geq 1 \\ j_{3}=1}} \frac{-i\beta e^{2i\beta}}{1-e^{i\beta}}\frac{(-1)^{j_{1}+j_{3}}i^{j_{1}+j_{2}+j_{3}}}{2j_{1}!j_{2}!j_{3}!}\phi(\alpha - \beta)^{j_{1}}\phi(\alpha)^{j_{2}} \\
    &\cdot\int_{0}^{1}e^{-i\beta s}\phi(\alpha + \beta(-1+s))^{j_{3}}(-1+s)ds \\
    &\cdot\biggl(\sum_{m=1}^{\infty}\frac{i\beta}{1-e^{i\beta}}\int_{0}^{1}e^{-i(s-1)\beta}\frac{(-i\phi(\alpha + (s-1)\beta))^{m}}{m!}ds\biggr)^{n}
\end{align*}
\begin{align*}
    \widetilde{B_{1,2}}(\alpha, \beta) = -&\sum_{\substack{j_{1}+j_{2}+j_{3} \geq 1 \\ n=1}} \frac{-i\beta e^{2i\beta}}{1-e^{i\beta}}\frac{(-1)^{j_{1}+j_{3}}i^{j_{1}+j_{2}+j_{3}}}{2j_{1}!j_{2}!j_{3}!}\phi(\alpha - \beta)^{j_{1}}\phi(\alpha)^{j_{2}} \\
    &\cdot\int_{0}^{1}e^{-i\beta s}\phi(\alpha + \beta(-1+s))^{j_{3}}(-1+s)ds \\
    &\cdot\biggl(\sum_{m=1}^{\infty}\frac{i\beta}{1-e^{i\beta}}\int_{0}^{1}e^{-i(s-1)\beta}\frac{(-i\phi(\alpha + (s-1)\beta))^{m}}{m!}ds\biggr)^{n} \\
    \widetilde{B_{1,3}}(\alpha, \beta) = - &\frac{-i\beta e^{2i\beta}}{1-e^{i\beta}}\cdot\frac{1}{2}\int_{0}^{1}e^{-i\beta s}(-1+s)ds \\
    &\cdot\sum_{m=2}^{\infty}\frac{i\beta}{1-e^{i\beta}}\int_{0}^{1}e^{-i(s-1)\beta}\frac{(-i\phi(\alpha+(s-1)\beta))^{m}}{m!}ds.
\end{align*}
Due to similarity in calculations, we only calculate the Fourier modes of \(\widetilde{B_{1,1}}(\alpha, \beta)\) for illustration. We note that
\begin{align*}
    &\mathcal{F}(\widetilde{B_{1,1}})(k_{1}, \beta) = -\sum_{j_{1}+j_{2}+n \geq 1} \frac{-i\beta e^{2i\beta}}{1-e^{i\beta}}\frac{(-1)^{j_{1}+1}i^{j_{1}+j_{2}+1}}{2j_{1}!j_{2}!} \\
    &\cdot\mathcal{F}\biggl(\phi(\alpha-\beta)^{j_{1}}\phi(\alpha)^{j_{2}}\int_{0}^{1}e^{-i\beta s}\phi(\alpha + \beta(-1+s))(-1+s)ds \\
    &\cdot\biggl(\sum_{m=1}^{\infty}\frac{i\beta}{1-e^{i\beta}}\int_{0}^{1}e^{-i(s-1)\beta}\frac{(-i\phi(\alpha+(s-1)\beta))^{m}}{m!}ds\biggr)^{n}\biggr)(k_{1}).
\end{align*}
Letting
\begin{align*}
    P(k) &= \sum_{m=1}^{\infty}\frac{(-i)^{m}}{m!}\mathcal{F}(\phi^{m})(k),
\end{align*}
we can write
\begin{align*}
    &\mathcal{F}\biggl(\phi(\alpha-\beta)^{j_{1}}\phi(\alpha)^{j_{2}}\int_{0}^{1}e^{-i\beta s}\phi(\alpha + \beta(-1+s))(-1+s)ds \\
    &\cdot\biggl(\sum_{m=1}^{\infty}\frac{i\beta}{1-e^{i\beta}}\int_{0}^{1}e^{-i(s-1)\beta}\frac{(-i\phi(\alpha+(s-1)\beta))^{m}}{m!}ds\biggr)^{n}\biggr)(k_{1}) \\
    =&\sum_{k_{2}, \dots, k_{j_{1}+j_{2}+n+1} \in \mathbb{Z}}\prod_{d=1}^{j_{1}+j_{2}}\mathcal{F}(\phi)(k_{d}-k_{d+1})\mathcal{F}(\phi)(k_{j_{1}+j_{2}+n+1}) \\
    &\cdot\prod_{d=1}^{n}P(k_{j_{1}+j_{2}+d}-k_{j_{1}+j_{2}+d+1}) \\
    &\cdot I_{n}(k_{1}, k_{j_{1}+1}, k_{j_{1}+j_{2}+1}, \dots, k_{j_{1}+j_{2}+n}, k_{j_{1}+j_{2}+n+1}, \beta).
\end{align*}
Then
\begin{align*}
    &\frac{\gamma}{4\pi}\int_{-\pi}^{\pi}\mathcal{F}(\widetilde{B_{1,1}})(k_{1}, \beta)d\beta \\
    = -&\sum_{j_{1}+j_{2}+n\geq 1}\frac{(-1)^{j_{1}+1}i^{j_{1}+j_{2}+1}}{2j_{1}!j_{2}!} \\
    &\cdot\sum_{k_{2}, \dots, k_{j_{1}+j_{2}+n+1} \in \mathbb{Z}}\prod_{d=1}^{j_{1}+j_{2}}\mathcal{F}(\phi)(k_{d}-k_{d+1})\mathcal{F}(\phi)(k_{j_{1}+j_{2}+n+1}) \\
    &\cdot\prod_{d=1}^{n}P(k_{j_{1}+j_{2}+d}-k_{j_{1}+j_{2}+d+1})\frac{\gamma}{4\pi}\int_{-\pi}^{\pi}\frac{-i\beta e^{2i\beta}}{1-e^{i\beta}}
\end{align*}
\begin{align*}
    &\cdot I_{n}(k_{1}, k_{j_{1}+1}, k_{j_{1}+j_{2}+1}, \dots, k_{j_{1}+j_{2}+n}, k_{j_{1}+j_{2}+n+1}, \beta)d\beta.
\end{align*}
By Lemma \ref{inbound},
\begin{align*}
    &\abs{\frac{\gamma}{4\pi}\int_{-\pi}^{\pi}\mathcal{F}(\widetilde{B_{1,1}})(k_{1}, \beta)d\beta} \\
    \leq& \sum_{j_{1}+j_{2}+n\geq 1}\frac{C_{n}}{2j_{1}!j_{2}!} \sum_{k_{2}, \dots, k_{j_{1}+j_{2}+n+1} \in \mathbb{Z}}\prod_{d=1}^{j_{1}+j_{2}}\abs{\mathcal{F}(\phi)(k_{d}-k_{d+1})}\abs{\mathcal{F}(\phi)(k_{j_{1}+j_{2}+n+1})} \\
    &\cdot\prod_{d=1}^{n}\abs{P(k_{j_{1}+j_{2}+d}-k_{j_{1}+j_{2}+d+1})}.
\end{align*}
Similarly,
\begin{align*}
    &\abs{\frac{\gamma}{4\pi}\int_{-\pi}^{\pi}\mathcal{F}(\widetilde{B_{1,2}})(k_{1}, \beta)d\beta} \\
    \leq& \sum_{j_{1}+j_{2}+j_{3} \geq 1}\frac{C_{1}}{2j_{1}!j_{2}!j_{3}!}\sum_{k_{2}, \dots, k_{j_{1}+j_{2}+2} \in \mathbb{Z}}\prod_{d=1}^{j_{1}+j_{2}}\abs{\mathcal{F}(\phi)(k_{d}-k_{d+1})}\abs{\mathcal{F}(\phi^{j_{3}})(k_{j_{1}+j_{2}+2})} \\
    &\cdot\abs{P(k_{j_{1}+j_{2}+1}-k_{j_{1}+j_{2}+2})}
\end{align*}
and
\begin{align*}
    &\abs{\frac{\gamma}{4\pi}\int_{-\pi}^{\pi}\mathcal{F}(\widetilde{B_{1,3}})(k_{1}, \beta)d\beta} \leq \frac{1}{2}\cdot\frac{\gamma}{4\pi}\biggl(\frac{\pi}{2}\sqrt{1+\frac{\pi^{2}}{4}}\cdot\pi^{2}+2\pi\biggr)\abs{\widetilde{P}(k_{1})},
\end{align*}
where
\begin{align*}
    \widetilde{P}(k) &= \sum_{m=2}^{\infty}\frac{(-i)^{m}}{m!}\mathcal{F}(\phi^{m})(k).
\end{align*}
Calculations associated with the rest of the terms in (\ref{U2sum}) can be done similarly.

\subsection{Estimating \texorpdfstring{\(\normz{U_{\geq 2}}\)}{U2}}
\label{estNormzU2}
Let us prove the following estimate for \(\normz{U_{\geq 2}}\).
\begin{lemma} \label{estimateforNormzofU2}
    For some monotone increasing functions \(D_{1}\) and \(D_{2}\) of \(\normz{\phi}\),
    \begin{align*}
        \normz{U_{\geq 2}} \leq D_{1}(\normz{\phi})\normz{\phi}^{2}+D_{2}(\normz{\phi})\normz{\phi}\normz{\phi'}.
    \end{align*}
\end{lemma}
Before commencing the proof of Lemma \ref{estimateforNormzofU2}, let us introduce the setup for the proof. For ease of notation, we define the \(l_{\nu}^{1}\) norm of a sequence \(a=a(k)\) defined on \(\mathbb{Z}\) by
\begin{align*}
    \norml{a} = \sum_{k \in \mathbb{Z}}e^{\nu(t)\abs{k}}\abs{a(k)}.
\end{align*}
The following estimate of the \(l_{\nu}^{1}\) norm of the convolution is frequently used.
\begin{proposition} \label{convest}
    If \(a_{1}, \dots, a_{n}\) are sequences on \(\mathbb{Z}\) whose \(l_{\nu}^{1}\) norms are finite, then
    \begin{align*}
        \norml{a_{1}*\cdots * a_{n}} \leq \prod_{j=1}^{n}\norml{a_{j}}.
    \end{align*}
\end{proposition}
\begin{proof}
    It suffices to prove the case when \(n=2\). We have
    \begin{align*}
        \norml{a*b}&=\sum_{k \in \mathbb{Z}}e^{\nu(t)\abs{k}}\abs{(a*b)(k)}  \\
        &\leq \sum_{j \in \mathbb{Z}}\sum_{k \in \mathbb{Z}}e^{\nu(t)\abs{k-j}}e^{\nu(t)\abs{j}}\abs{a(k-j)}\abs{b(j)} \\
        &=\sum_{j \in \mathbb{Z}}e^{\nu(t)\abs{j}}\abs{b(j)}\sum_{k \in \mathbb{Z}}e^{\nu(t)\abs{k-j}}\abs{a(k-j)} \\
        &=\norml{a}\norml{b},
    \end{align*}
    as needed.
\end{proof}
We note that
\begin{align}
    \norml{P} &\leq \sum_{m=1}^{\infty}\frac{\norm{\phi}_{\mathcal{F}_{\nu}^{0,1}}^{m}}{m!} = e^{\norm{\phi}_{\mathcal{F}_{\nu}^{0,1}}}-1, \label{normp} \\
    \norml{\widetilde{P}}&\leq\sum_{m=2}^{\infty}\frac{\norm{\phi}_{\mathcal{F}_{\nu}^{0,1}}^{m}}{m!} =e^{\norm{\phi}_{\mathcal{F}_{\nu}^{0,1}}}-\norm{\phi}_{\mathcal{F}_{\nu}^{0,1}}-1. \label{normptilde}
\end{align}
To begin the proof of Lemma \ref{estimateforNormzofU2}, we observe that
\begin{align*}
    \normz{U_{\geq 2}} \leq \sum_{j=1}^{8}\normz{\frac{\gamma}{4\pi}\int_{-\pi}^{\pi}\widetilde{B_{j}}(\alpha, \beta)d\beta} + \normz{\frac{\gamma}{4\pi}\int_{-\pi}^{\pi}\widetilde{B_{13}}(\alpha, \beta)d\beta}.
\end{align*}
This means that it suffices to estimate each of the \(\mathcal{F}_{\nu}^{0,1}\) norms on the right hand side. By Proposition \ref{convest} and (\ref{normp}), we obtain
\begin{align*}
    &\sum_{k_{1} \in \mathbb{Z}}e^{\nu(t)\abs{k_{1}}}\abs{\frac{\gamma}{4\pi}\int_{-\pi}^{\pi}\mathcal{F}(\widetilde{B_{1,1}})(k_{1}, \beta)d\beta} \\
    \leq&\sum_{j_{1}+j_{2}+n \geq 1} \frac{C_{n}}{2j_{1}!j_{2}!}\norm{\phi}_{\mathcal{F}_{\nu}^{0,1}}^{j_{1}+j_{2}+1}(e^{\norm{\phi}_{\mathcal{F}_{\nu}^{0,1}}}-1)^{n}.
\end{align*}
By Propositions \ref{p5} and \ref{convest} and (\ref{normp}), we obtain
\begin{align*}
    &\sum_{k_{1} \in \mathbb{Z}}e^{\nu(t)\abs{k_{1}}}\abs{\frac{\gamma}{4\pi}\int_{-\pi}^{\pi}\mathcal{F}(\widetilde{B_{1,2}})(k_{1}, \beta)d\beta} \\
    \leq& \sum_{j_{1}+j_{2}+j_{3} \geq 1}\frac{C_{1}}{2j_{1}j_{2}!j_{3}!}\norm{\phi}_{\mathcal{F}_{\nu}^{0,1}}^{j_{1}+j_{2}+j_{3}}(e^{\norm{\phi}_{\mathcal{F}_{\nu}^{0,1}}}-1).
\end{align*}
By (\ref{normptilde}), we obtain
\begin{align*}
    \sum_{k_{1} \in \mathbb{Z}}e^{\nu(t)\abs{k_{1}}}\abs{\frac{\gamma}{4\pi}\int_{-\pi}^{\pi}\mathcal{F}(\widetilde{B_{1,3}})(k_{1}, \beta)d\beta}\leq \frac{C_{2}}{2}(e^{\norm{\phi}_{\mathcal{F}_{\nu}^{0,1}}}-\norm{\phi}_{\mathcal{F}_{\nu}^{0,1}}-1).
\end{align*}
Calculations associated with the rest of the terms in (\ref{U2sum}) can be done similarly. This completes the proof of Lemma \ref{estimateforNormzofU2}.

\subsection{Estimating \texorpdfstring{\(\norms{U_{\geq 2}}\)}{U2}}
We prove a useful estimate for \(\norms{U_{\geq 2}}\), \(s>0\).
\begin{lemma} \label{estimateforNormsofU2}
    For \(s>0\),
    \begin{align*}
        \norms{U_{\geq 2}} \leq& F_{1}(\normz{\phi})\normz{\phi}\norms{\phi}+F_{2}(\normz{\phi})\normz{\phi}^{2}\norms{\phi}\\
        &+F_{3}(\normz{\phi})\normz{\phi'}\norms{\phi}+F_{4}(\normz{\phi})\normz{\phi}\norms{\phi'},
    \end{align*}
    where \(F_{1}\), \(F_{2}\), \(F_{3}\), and \(F_{4}\) are monotone increasing functions of \(\normz{\phi}\).
\end{lemma}
Before commencing the proof of Lemma \ref{estimateforNormsofU2}, let us introduce the setup for the proof. For ease of notation, we define the \(l_{\nu}^{s}\) norm of a sequence \(a=a(k)\) defined on \(\mathbb{Z}\) by
\begin{align*}
    \normls{a} = \sum_{k \neq 0}e^{\nu(t)\abs{k}}\abs{k}^{s}\abs{a_{k}}.
\end{align*}
The following estimate of the \(l_{\nu}^{s}\) norm of the convolution is frequently used.
\begin{proposition} \label{convests}
    Let \(s>0\). If \(a_{1}, \dots, a_{n}\) are sequences on \(\mathbb{Z}\) whose \(l_{\nu}^{s}\) norms are finite, then
    \begin{align*}
        \normls{a_{1}*\cdots * a_{n}} \leq b(n,s)\sum_{j=1}^{n}\normls{a_{j}}\prod_{\substack{k=1 \\ k \neq j}}^{n}\norml{a_{k}}.
    \end{align*}
\end{proposition}
\begin{proof}
    We note that for any \(k_{1}, \dots, k_{n} \in \mathbb{Z}\),
    \begin{align*}
        \abs{k_{1}}^{s} \leq b(n,s)(\abs{k_{1}-k_{2}}^{s}+\abs{k_{2}-k_{3}}^{s}+\cdots+\abs{k_{n-1}-k_{n}}^{s}+\abs{k_{n}}^{s}),
    \end{align*}
    which follows from the convexity of the function \(\abs{\cdot}^{s}\) for \(s \geq 1\) and the triangle inequality for \(0<s<1\). Using (\ref{convformula}), we then obtain
    \begin{align*}
        \normls{a_{1}*\cdots*a_{n}} =&\sum_{k_{1} \in \mathbb{Z}}e^{\nu(t)\abs{k_{1}}}\abs{k_{1}}^{s}\abs{(a_{1}*\cdots * a_{n})(k_{1})} \\
        \leq&\sum_{k_{1} \in \mathbb{Z}}\sum_{k_{2},\dots,k_{n} \in \mathbb{Z}}e^{\nu(t)\abs{k_{1}}}\abs{k_{1}}^{s}\prod_{d=1}^{n-1}\abs{a_{d}(k_{d}-k_{d+1})}\abs{a_{n}(k_{n})} \\
        \leq& \sum_{j=2}^{n}\sum_{k_{1}\in \mathbb{Z}}\sum_{k_{2},\dots,k_{n}\in \mathbb{Z}}e^{\nu(t)\abs{k_{1}}}b(n,s)\abs{k_{j-1}-k_{j}}^{s} \\
        &\cdot\prod_{d=1}^{n-1}\abs{a_{d}(k_{d}-k_{d+1})}\abs{a_{n}(k_{n})} \\
        &+\sum_{k_{1}\in \mathbb{Z}}\sum_{k_{2},\dots,k_{n} \in \mathbb{Z}}e^{\nu(t)\abs{k_{1}}}b(n,s)\abs{k_{n}}^{s}\prod_{d=1}^{n-1}\abs{a_{d}(k_{d}-k_{d+1})}\abs{a_{n}(k_{n})}.
    \end{align*}
    For \(j \in \{2,\dots,n\}\), we have
    \begin{align*}
        &\sum_{k_{1}\in \mathbb{Z}}\sum_{k_{2},\dots,k_{n}\in \mathbb{Z}}e^{\nu(t)\abs{k_{1}}}b(n,s)\abs{k_{j-1}-k_{j}}^{s}\prod_{d=1}^{n-1}\abs{a_{d}(k_{d}-k_{d+1})}\abs{a_{n}(k_{n})} \\
        \leq&b(n,s)\sum_{k_{n} \in \mathbb{Z}}\abs{a_{n}(k_{n})}e^{\nu(t)\abs{k_{n}}}\sum_{k_{n-1}\in \mathbb{Z}}\abs{a_{n-1}(k_{n-1}-k_{n})}e^{\nu(t)\abs{k_{n-1}-k_{n}}} \\
        &\cdot\cdots\cdot\sum_{k_{j-1} \in \mathbb{Z}}\abs{k_{j-1}-k_{j}}^{s}\abs{a_{j-1}(k_{j-1}-k_{j})}e^{\nu(t)\abs{k_{j-1}-k_{j}}} \\
        &\cdot\cdots\cdot\sum_{k_{2}\in \mathbb{Z}}\abs{a_{2}(k_{2}-k_{3})}e^{\nu(t)\abs{k_{2}-k_{3}}}\sum_{k_{1}\in\mathbb{Z}}\abs{a_{1}(k_{1}-k_{2})}e^{\nu(t)\abs{k_{1}-k_{2}}}.
    \end{align*}
    Making the changes of variables
    \begin{align*}
        k_{1}'&=k_{1}-k_{2} \\
        k_{2}'&=k_{2}-k_{3} \\
        &\vdots \\
        k_{n-1}'&=k_{n-1}-k_{n}
    \end{align*}
    in that order, we obtain
    \begin{align*}
        &\sum_{k_{1}\in \mathbb{Z}}\sum_{k_{2},\dots,k_{n}\in \mathbb{Z}}e^{\nu(t)\abs{k_{1}}}b(n,s)\abs{k_{j-1}-k_{j}}^{s}\prod_{d=1}^{n-1}\abs{a_{d}(k_{d}-k_{d+1})}\abs{a_{n}(k_{n})} \\
        \leq&b(n,s)\normls{a_{j-1}}\prod_{\substack{k=1 \\ k\neq j}}^{n}\norml{a_{k}}.
    \end{align*}
    Similarly,
    \begin{align*}
        &\sum_{k_{1}\in \mathbb{Z}}\sum_{k_{2},\dots,k_{n} \in \mathbb{Z}}e^{\nu(t)\abs{k_{1}}}b(n,s)\abs{k_{n}}^{s}\prod_{d=1}^{n-1}\abs{a_{d}(k_{d}-k_{d+1})}\abs{a_{n}(k_{n})} \\
        \leq&b(n,s)\normls{a_{n}}\prod_{k=1}^{n-1}\norml{a_{k}}.
    \end{align*}
    This completes the proof.
\end{proof}
We note that
\begin{align}
    \normls{P} \leq&\biggl(\sum_{m=1}^{\infty}\frac{b(m,s)}{(m-1)!}\normz{\phi}^{m-1}\biggr)\norms{\phi}, \label{normlsp} \\
    \normls{\widetilde{P}}\leq&\left(\sum_{m=2}^{\infty}\frac{b(m,s)}{(m-1)!}\normz{\phi}^{m-1}\right)\norms{\phi}. \label{normlsptilde}
\end{align}
To begin the proof of Lemma \ref{estimateforNormsofU2}, we observe that
\begin{align*}
    \norms{U_{\geq 2}} \leq \sum_{j=1}^{8}\norms{\frac{\gamma}{4\pi}\int_{-\pi}^{\pi}\widetilde{B_{j}}(\alpha,\beta)d\beta}+\norms{\frac{\gamma}{4\pi}\int_{-\pi}^{\pi}\widetilde{B_{13}}(\alpha,\beta)d\beta}.
\end{align*}
This means that it suffices to estimate each of the \(\dot{\mathcal{F}}_{\nu}^{s,1}\) norms on the right hand side. By Proposition \ref{convests} and (\ref{normlsp}), we obtain
\begin{align*}
    &\sum_{k_{1} \neq 0}e^{\nu(t)\abs{k_{1}}}\abs{k_{1}}^{s}\abs{\frac{\gamma}{4\pi}\int_{-\pi}^{\pi}\mathcal{F}(\widetilde{B_{1,1}})(k_{1}, \beta)d\beta} \\
    \leq&\sum_{j_{1}+j_{2}+n\geq 1}\frac{C_{n}}{2j_{1}!j_{2}!}b(j_{1}+j_{2}+n+1,s) \\    \cdot&\biggl(\norm{\phi}_{\dot{\mathcal{F}}_{\nu}^{s,1}}\norm{\phi}_{\mathcal{F}_{\nu}^{0,1}}^{j_{1}+j_{2}}(e^{\norm{\phi}_{\mathcal{F}_{\nu}^{0,1}}}-1)^{n}(j_{1}+j_{2}+1) \\
    &+\biggl(\sum_{m=1}^{\infty}\frac{b(m,s)}{(m-1)!}\normz{\phi}^{m-1}\biggr)\norms{\phi}\norm{\phi}_{\mathcal{F}_{\nu}^{0,1}}^{j_{1}+j_{2}+1}(e^{\norm{\phi}_{\mathcal{F}_{\nu}^{0,1}}}-1)^{n-1}\cdot n\biggr).
\end{align*}
By Proposition \ref{convests} and (\ref{normlsp}), we obtain
\begin{align*}
    &\sum_{k_{1} \neq 0}e^{\nu(t)\abs{k_{1}}}\abs{k_{1}}^{s}\abs{\frac{\gamma}{4\pi}\int_{-\pi}^{\pi}\mathcal{F}(\widetilde{B_{1,2}})(k_{1}, \beta)d\beta} \\
    \leq&\sum_{j_{1}+j_{2}+j_{3} \geq 1}\frac{C_{1}}{2j_{1}!j_{2}!j_{3}!}b(j_{1}+j_{2}+2,s) \\
    \cdot& \biggl(\norms{\phi}\normz{\phi}^{j_{1}+j_{2}+j_{3}-1}(e^{\normz{\phi}}-1)\cdot(j_{1}+j_{2}) \\
    &+b(j_{3},s)\norms{\phi}\normz{\phi}\cdot j_{3}\cdot\normz{\phi}^{j_{1}+j_{2}}(e^{\normz{\phi}}-1) \\
    &+\biggl(\sum_{m=1}^{\infty}\frac{b(m,s)}{(m-1)!}\normz{\phi}^{m-1}\biggr)\norms{\phi}\normz{\phi}^{j_{1}+j_{2}+j_{3}}\biggr).
\end{align*}
By (\ref{normlsptilde}), we obtain
\begin{align*}
    &\sum_{k_{1} \neq 0} e^{\nu(t)\abs{k_{1}}}\abs{k_{1}}^{s}\abs{\frac{\gamma}{4\pi}\int_{-\pi}^{\pi}\mathcal{F}(\widetilde{B_{1,3}})(k_{1}, \beta)d\beta} \\
    \leq&\frac{C_{1}}{2}\left(\sum_{m=2}^{\infty}\frac{b(m,s)}{(m-1)!}\normz{\phi}^{m-1}\right)\norms{\phi}.
\end{align*}
Calculations associated with the rest of the terms in (\ref{U2sum}) can be done similarly. This completes the proof of Lemma \ref{estimateforNormsofU2}.

\section{Estimating \texorpdfstring{\((U_{\geq 2})_{\alpha}\)}{U2deriv}}
\label{estimationU2deriv}
In Section \ref{estimationU2}, we derived that
\begin{align*}
    U_{\geq 2}(\alpha) = \mbox{Re}\biggl(\frac{\gamma}{4\pi}\int_{-\pi}^{\pi}B(\alpha, \beta)d\beta\biggr),
\end{align*}
where
\begin{align} \label{B}
    B(\alpha, \beta) = \sum_{j=1}^{8}\widetilde{B_{j}}(\alpha, \beta) + \widetilde{B_{13}}(\alpha, \beta).
\end{align}
To estimate the \(\dot{\mathcal{F}}_{\nu}^{s,1}\) norm of \((U_{\geq 2})_{\alpha}\), we differentiate the right hand side with respect to \(\alpha\). For example, recalling that
\begin{align*}
    \widetilde{B_{1}}(\alpha, \beta) = \widetilde{B_{1,1}}(\alpha, \beta)+\widetilde{B_{1,2}}(\alpha, \beta)+\widetilde{B_{1,3}}(\alpha, \beta),
\end{align*}
we note that
\begin{align*}
    (\widetilde{B_{1,1}})_{\alpha}(\alpha, \beta) &= \sum_{j=1}^{4}B_{1,1}^{j}(\alpha, \beta),
\end{align*}
where
\begin{align}
    B_{1,1}^{1}(\alpha, \beta) = -&\sum_{j_{1}+j_{2}+n \geq 1}\frac{-i\beta e^{2i\beta}}{1-e^{i\beta}}\cdot\frac{(-1)^{j_{1}+1}i^{j_{1}+j_{2}+1}}{2j_{1}!j_{2}!}j_{1}\phi(\alpha-\beta)^{j_{1}-1} \nonumber \\
    &\cdot\phi_{\alpha}(\alpha-\beta)\phi(\alpha)^{j_{2}}\int_{0}^{1}e^{-i\beta s}\phi(\alpha+\beta(-1+s))(-1+s)ds \nonumber \\
    &\cdot \biggl(\sum_{m=1}^{\infty}\frac{i\beta}{1-e^{i\beta}}\int_{0}^{1}e^{-i(s-1)\beta}\frac{(-i\phi(\alpha+(s-1)\beta))^{m}}{m!}ds\biggr)^{n}, \nonumber
\end{align}
\begin{align}
    B_{1,1}^{2}(\alpha, \beta) = -&\sum_{j_{1}+j_{2}+n \geq 1}\frac{-i\beta e^{2i\beta}}{1-e^{i\beta}}\cdot\frac{(-1)^{j_{1}+1}i^{j_{1}+j_{2}+1}}{2j_{1}!j_{2}!}\phi(\alpha-\beta)^{j_{1}}j_{2}\phi(\alpha)^{j_{2}-1}\phi_{\alpha}(\alpha) \nonumber \\
    &\cdot\int_{0}^{1}e^{-i\beta s}\phi(\alpha+\beta(-1+s))(-1+s)ds \nonumber \\
    &\cdot \biggl(\sum_{m=1}^{\infty}\frac{i\beta}{1-e^{i\beta}}\int_{0}^{1}e^{-i(s-1)\beta}\frac{(-i\phi(\alpha+(s-1)\beta))^{m}}{m!}ds\biggr)^{n} \nonumber \\
    B_{1,1}^{3}(\alpha, \beta) = -&\sum_{j_{1}+j_{2}+n \geq 1}\frac{-i\beta e^{2i\beta}}{1-e^{i\beta}}\cdot\frac{(-1)^{j_{1}+1}i^{j_{1}+j_{2}+1}}{2j_{1}!j_{2}!}\phi(\alpha-\beta)^{j_{1}}\phi(\alpha)^{j_{2}} \nonumber \\
    &\cdot\int_{0}^{1}e^{-i\beta s}\phi_{\alpha}(\alpha+\beta(-1+s))(-1+s)ds \nonumber \\
    &\cdot \biggl(\sum_{m=1}^{\infty}\frac{i\beta}{1-e^{i\beta}}\int_{0}^{1}e^{-i(s-1)\beta}\frac{(-i\phi(\alpha+(s-1)\beta))^{m}}{m!}ds\biggr)^{n}, \nonumber \\
    B_{1,1}^{4}(\alpha, \beta) = -&\sum_{j_{1}+j_{2}+n \geq 1}\frac{-i\beta e^{2i\beta}}{1-e^{i\beta}}\cdot\frac{(-1)^{j_{1}+1}i^{j_{1}+j_{2}+1}}{2j_{1}!j_{2}!}\phi(\alpha-\beta)^{j_{1}}\phi(\alpha)^{j_{2}} \nonumber \\
    &\cdot\int_{0}^{1}e^{-i\beta s}\phi(\alpha+\beta(-1+s))(-1+s)ds \nonumber \\
    &\cdot n\biggl(\sum_{m=1}^{\infty}\frac{i\beta}{1-e^{i\beta}}\int_{0}^{1}e^{-i(s-1)\beta}\frac{(-i\phi(\alpha+(s-1)\beta))^{m}}{m!}ds\biggr)^{n-1} \nonumber \\
    &\cdot\biggl(\sum_{m=1}^{\infty}\frac{i\beta}{1-e^{i\beta}} \nonumber \\
    &\cdot\int_{0}^{1}e^{-i(s-1)\beta}\frac{m(-i\phi(\alpha+(s-1)\beta))^{m-1}(-i)\phi_{\alpha}(\alpha+(s-1)\beta)}{m!}ds\biggr). \nonumber
\end{align}
After differentiation, certain terms contain the second derivative of \(\phi\). To such terms, we apply integration by parts to re-express them in terms of lower-order derivatives of \(\phi\).

\subsection{Estimating Fourier Modes of \texorpdfstring{\((U_{\geq 2})_{\alpha}\)}{U2deriv}}
\label{estimatingfouriermodesofU2deriv}
We use arguments as in Section \ref{efmodesU2} to estimate the Fourier modes of \((U_{\geq 2})_{\alpha}\). For example,
\begin{align*}
    &\abs{\frac{\gamma}{4\pi}\int_{-\pi}^{\pi}\mathcal{F}(B_{1,1}^{1})(k_{1},\beta)d\beta} \leq \\
    &\sum_{j_{1}+j_{2}+n \geq 1}\frac{j_{1}C_{n}}{2j_{1}!j_{2}!}\sum_{k_{2}, \dots, k_{j_{1}+j_{2}+n+1} \in \mathbb{Z}}\prod_{d=1}^{j_{1}-1}\abs{\mathcal{F}(\phi)(k_{d}-k_{d+1})}\cdot\abs{\mathcal{F}(\phi')(k_{j_{1}}-k_{j_{1}+1})} \\
    &\cdot\prod_{d=j_{1}+1}^{j_{1}+j_{2}}\abs{\mathcal{F}(\phi)(k_{d}-k_{d+1})}\cdot\prod_{d=j_{1}+j_{2}+1}^{j_{1}+j_{2}+n}\abs{P(k_{d}-k_{d+1})}\cdot\abs{\mathcal{F}(\phi)(k_{j_{1}+j_{2}+n+1})}.
\end{align*}
Calculations associated with the rest of the terms in the derivative of (\ref{B}) can be done similarly.

\subsection{Estimating \texorpdfstring{\(\norms{(U_{\geq 2})_{\alpha}}\)}{U2deriv}}
We prove a useful estimate for \(\norms{(U_{\geq 2})_{\alpha}}\), \(s>0\).
\begin{lemma} \label{estimateforNormsofU2deriv}
    For \(s>0\),
    \begin{align*}
        \norms{(U_{\geq 2})_{\alpha}} \leq& R_{1}(\normz{\phi})\normz{\phi'}\norms{\phi}+R_{2}(\normz{\phi})\normz{\phi}\norms{\phi'}\\
        &+R_{3}(\normz{\phi})\normz{\phi}\normz{\phi'}\norms{\phi} \\
        &+R_{4}(\normz{\phi})\normz{\phi'}^{2}\norms{\phi}+R_{5}(\normz{\phi})\normz{\phi'}\norms{\phi'},
    \end{align*}
    where \(R_{1}\), \(R_{2}\), \(R_{3}\), \(R_{4}\), and \(R_{5}\) are monotone increasing functions of \(\normz{\phi}\).
\end{lemma}
We use estimates from Section \ref{estimatingfouriermodesofU2deriv} to prove Lemma \ref{estimateforNormsofU2deriv}. For example,
\begin{align*}
    &\norms{\frac{\gamma}{4\pi}\int_{-\pi}^{\pi}B_{1,1}^{1}(\alpha, \beta)d\beta} \\
    \leq&\sum_{j_{1}+j_{2}+n \geq 1}\frac{j_{1}C_{n}}{2j_{1}!j_{2}!}\normls{\abs{\mathcal{F}(\phi)} * \cdots * \abs{\mathcal{F}(\phi)} * \abs{P} * \cdots * \abs{P} * \abs{\mathcal{F}(\phi')}} \\
    \leq&\sum_{j_{1}+j_{2}+n \geq 1}\frac{j_{1}C_{n}}{2j_{1}!j_{2}!}b(j_{1}+j_{2}+n+1,s) \\
    \cdot&\biggl(\norms{\phi}\normz{\phi}^{j_{1}+j_{2}-1}\normz{\phi'}(e^{\normz{\phi}}-1)^{n}(j_{1}+j_{2}) \\
    &+\norms{\phi'}\normz{\phi}^{j_{1}+j_{2}}(e^{\normz{\phi}}-1)^{n} \\
    &+\biggl(\sum_{m=1}^{\infty}\frac{b(m,s)}{(m-1)!}\normz{\phi}^{m-1}\biggr)\norms{\phi}\normz{\phi}^{j_{1}+j_{2}}\normz{\phi'}(e^{\normz{\phi}}-1)^{n-1}\cdot n\biggr).
\end{align*}
Calculations associated with the rest of the terms in \((U_{\geq 2})_{\alpha}\) can be done similarly. This completes the proof of Lemma \ref{estimateforNormsofU2deriv}.

\section{Proof of the Main Theorem}
\subsection{Proof of the Main \emph{a priori} Estimate}
\label{mainaprioriestimate}
To complete the estimate for the \(\dot{\mathcal{F}}_{\nu}^{s,1}\) norm of \(\widetilde{\mathcal{N}}\), we let \(s=1\). Recalling (\ref{threeterms}), we use Lemmas \ref{estimatingterm2lem} and \ref{estimatingterm3lem}, the estimates of the \(\mathcal{F}_{\nu}^{0,1}\) norm of \(U_{1}\) and \(U_{\geq 2}\) in Sections \ref{estNormzU1} and \ref{estNormzU2}, respectively, Lemma \ref{estimateforNormsofU2deriv}, and Proposition \ref{p3} to obtain an estimate for the \(\dot{\mathcal{F}}_{\nu}^{1,1}\) norm of \(\widetilde{\mathcal{N}}\). We substitute this estimate of the \(\dot{\mathcal{F}}_{\nu}^{1,1}\) norm of \(\widetilde{\mathcal{N}}\) into (\ref{aprioriestimateforNormsphi}) and use the fact that \(C_{I}\), \(A\), \(A_{1}^{-1}\), \(R_{1}\), \(R_{2}\), \(R_{3}\), \(R_{4}\), \(R_{5}\), \(D_{1}\), and \(D_{2}\) are all monotone increasing and Proposition \ref{p3} to obtain
\begin{align} \label{differentialinequality}
    \frac{d}{dt}\normx{\phi}{1} &\leq -\biggl(\Lambda(\normx{\phi}{1})-\nu'(t)\biggr)\normx{\phi}{2},
\end{align}
where
\begin{align*}
    \Lambda(\normx{\phi}{1}) =& \frac{1}{2\biggl(C_{I}(\normx{\phi}{1})\normx{\phi}{1}+1\biggr)}\pi\frac{2}{R}\frac{\gamma}{4\pi}-\frac{\gamma}{4\pi}\frac{1}{R}A(\normx{\phi}{1})\normx{\phi}{1} \\
    &-\frac{1}{R}\frac{1}{A_{1}(\normx{\phi}{1})}\biggl(R_{1}(\normx{\phi}{1})\normx{\phi}{1}+R_{2}(\normx{\phi}{1})\normx{\phi}{1}\\
    &+R_{3}(\normx{\phi}{1})\normx{\phi}{1}^{2}+R_{4}(\normx{\phi}{1})\normx{\phi}{1}^{2} \\
    &+R_{5}(\normx{\phi}{1})\normx{\phi}{1} \\
    &+3\left(H_{3}\normx{\phi}{1}+H_{4}\normx{\phi}{1}\right) \\
    &+3\left(D_{1}(\normx{\phi}{1})\normx{\phi}{1}^{2}+D_{2}(\normx{\phi}{1})\normx{\phi}{1}^{2}\right)\left(1+2\normx{\phi}{1}\right) \\
    &+\left(D_{1}(\normx{\phi}{1})\normx{\phi}{1}+D_{2}(\normx{\phi}{1})\normx{\phi}{1}\right)\left(1+2\normx{\phi}{1}\right)
\end{align*}
\begin{align*}
    &+6\normx{\phi}{1}\left(H_{3}\normx{\phi}{1}+H_{4}\normx{\phi}{1}\right) \\
    &+2\left(H_{3}\normx{\phi}{1}+H_{4}\normx{\phi}{1}\right)\biggr).
\end{align*}
Integrating (\ref{differentialinequality}) with respect to time, we obtain
\begin{align}
    \normx{\phi(t)}{1}+\int_{0}^{t}(\Lambda(\normx{\phi(\tau)}{1})-\nu'(\tau))\normx{\phi(\tau)}{2}d\tau \leq \norm{\phi_{0}}_{\dot{\mathcal{F}}^{1,1}}. \nonumber
\end{align}
Choose the initial datum such that \(\Lambda(\norm{\phi_{0}}_{\dot{\mathcal{F}}^{1,1}})>0\) and \(\nu_{0} \in (0,\Lambda(\norm{\phi_{0}}_{\dot{\mathcal{F}}^{1,1}}))\). Then \(\Lambda(\norm{\phi_{0}}_{\dot{\mathcal{F}}^{1,1}})-\nu'(0) >0\). It follows that for all \(t \in [0,\infty)\),
\begin{align*}
    \normx{\phi(t)}{1} &\leq \norm{\phi_{0}}_{\dot{\mathcal{F}}^{1,1}}-\int_{0}^{t}\biggl(\Lambda(\normx{\phi(\tau)}{1})-\nu'(\tau)\biggr)\normx{\phi(\tau)}{2}d\tau \\
    &\leq \norm{\phi_{0}}_{\dot{\mathcal{F}}^{1,1}}-\int_{0}^{t}\biggl(\Lambda(\norm{\phi_{0}}_{\dot{\mathcal{F}}^{1,1}})-\nu_{0}\biggr)\normx{\phi(\tau)}{2}d\tau.
\end{align*}

\subsection{Boundedness of \texorpdfstring{\(\mathcal{F}(\theta)(0)\)}{Ftheta0}}
Using the \emph{a priori} estimate for \(\phi = \theta - \hat{\theta}(0)\) in Section \ref{mainaprioriestimate}, we now show that \(\hat{\theta}(0)\) is bounded in time. We take the zeroth Fourier mode of (\ref{thetaevol}) and then integrate with respect to time to obtain
\begin{align*}
    \hat{\theta}(0)-\hat{\theta}_{0}(0) =&\int_{0}^{t}\frac{2\pi}{L(\tau)}\mathcal{F}\biggl(T_{1}(\alpha)(1+\theta_{\alpha}(\alpha))\biggr)(0)d\tau \\
    &+\int_{0}^{t}\frac{2\pi}{L(t)}\mathcal{F}\biggl(T_{\geq 2}(\alpha)(1+\theta_{\alpha}(\alpha))\biggr)(0)d\tau.
\end{align*}
Combining the estimate
\begin{align*}
    \abs{\mathcal{F}\biggl(T_{1}(\alpha)(1+\theta_{\alpha}(\alpha))\biggr)(0)} &\leq \norm{T_{1}}_{\mathcal{F}^{0,1}}(1+\norm{\phi}_{\mathcal{F}^{1,1}})
\end{align*}
with
\begin{align*}
    \norm{T_{1}}_{\mathcal{F}^{0,1}} \leq& 2\normx{U_{1}}{0} \leq 2\norm{U_{1}}_{\mathcal{F}_{\nu}^{0,1}} \leq2(H_{3}+H_{4})\normx{\phi}{1}
\end{align*}
\begin{align*}
    \norm{T_{\geq 2}}_{\mathcal{F}^{0,1}} \leq& 2\biggl(D_{1}(\normx{\phi}{1})\normx{\phi}{1}^{2}+D_{2}(\normx{\phi}{1})\normx{\phi}{1}^{2} \\
    &+2\normx{\phi}{1}\biggl(H_{3}\normx{\phi}{1}+H_{4}\normx{\phi}{1}\biggr) \\
    &+2\normx{\phi}{1}\biggl(D_{1}(\normx{\phi}{1})\normx{\phi}{1}^{2}+D_{2}(\normx{\phi}{1})\normx{\phi}{1}^{2}\biggr)\biggr)
\end{align*}
and
\begin{align*}
    \frac{2\pi}{L(t)} &\leq \frac{\sqrt{1+\frac{\pi}{2}(e^{2\norm{\phi_{0}}_{\dot{\mathcal{F}}^{1,1}}}-1)}}{R},
\end{align*}
we obtain
\begin{align*}
    \abs{\mathcal{F}(\theta)(0)} \leq&Y(\norm{\phi_{0}}_{\dot{\mathcal{F}}^{1,1}}),
\end{align*}
where
\begin{align} 
    &Y(\norm{\phi_{0}}_{\dot{\mathcal{F}}^{1,1}}) \label{uniformboundfunction} \\
    =&\abs{\mathcal{F}(\theta_{0})(0)} \nonumber \\
    &+\frac{\sqrt{1+\frac{\pi}{2}(e^{2\norm{\phi_{0}}_{\dot{\mathcal{F}}^{1,1}}}-1)}}{R}2(H_{3}+H_{4})\cdot\frac{\norm{\phi_{0}}_{\dot{\mathcal{F}}^{1,1}}}{\Lambda(\norm{\phi_{0}}_{\dot{\mathcal{F}}^{1,1}})-\nu_{0}} \nonumber \\
    &+ \frac{\sqrt{1+\frac{\pi}{2}(e^{2\norm{\phi_{0}}_{\dot{\mathcal{F}}^{1,1}}}-1)}}{R}\cdot 2(H_{3}+H_{4})\norm{\phi_{0}}_{\dot{\mathcal{F}}^{1,1}}\cdot\frac{\norm{\phi_{0}}_{\dot{\mathcal{F}}^{1,1}}}{\Lambda(\norm{\phi_{0}}_{\dot{\mathcal{F}}^{1,1}})-\nu_{0}} \nonumber \\
    &+\frac{\sqrt{1+\frac{\pi}{2}(e^{2\norm{\phi_{0}}_{\dot{\mathcal{F}}^{1,1}}}-1)}}{R} \nonumber \\
    &\cdot 2\biggl(D_{1}(\norm{\phi_{0}}_{\dot{\mathcal{F}}^{1,1}})\norm{\phi_{0}}_{\dot{\mathcal{F}}^{1,1}}+D_{2}(\norm{\phi_{0}}_{\dot{\mathcal{F}}^{1,1}})\norm{\phi_{0}}_{\dot{\mathcal{F}}^{1,1}} \nonumber \\
    &+2\biggl(H_{3}\norm{\phi_{0}}_{\dot{\mathcal{F}}^{1,1}}+H_{4}\norm{\phi_{0}}_{\dot{\mathcal{F}}^{1,1}}\biggr) \nonumber \\
    &+2\biggl(D_{1}(\norm{\phi_{0}}_{\dot{\mathcal{F}}^{1,1}})\norm{\phi_{0}}_{\dot{\mathcal{F}}^{1,1}}^{2}+D_{2}(\norm{\phi_{0}}_{\dot{\mathcal{F}}^{1,1}})\norm{\phi_{0}}_{\dot{\mathcal{F}}^{1,1}}^{2}\biggr)\biggr) \nonumber \\
    &\cdot\frac{\norm{\phi_{0}}_{\dot{\mathcal{F}}^{1,1}}}{\Lambda(\norm{\phi_{0}}_{\dot{\mathcal{F}}^{1,1}})-\nu_{0}} \nonumber \\
    &+\frac{\sqrt{1+\frac{\pi}{2}(e^{2\norm{\phi_{0}}_{\dot{\mathcal{F}}^{1,1}}}-1)}}{R} \nonumber \\
    &\cdot 2\biggl(D_{1}(\norm{\phi_{0}}_{\dot{\mathcal{F}}^{1,1}})\norm{\phi_{0}}_{\dot{\mathcal{F}}^{1,1}}^{2}+D_{2}(\norm{\phi_{0}}_{\dot{\mathcal{F}}^{1,1}})\norm{\phi_{0}}_{\dot{\mathcal{F}}^{1,1}}^{2} \nonumber \\
    &+2\norm{\phi_{0}}_{\dot{\mathcal{F}}^{1,1}}\biggl(H_{3}\norm{\phi_{0}}_{\dot{\mathcal{F}}^{1,1}}+H_{4}\norm{\phi_{0}}_{\dot{\mathcal{F}}^{1,1}}\biggr) \nonumber \\
    &+2\norm{\phi_{0}}_{\dot{\mathcal{F}}^{1,1}}\biggl(D_{1}(\norm{\phi_{0}}_{\dot{\mathcal{F}}^{1,1}})\norm{\phi_{0}}_{\dot{\mathcal{F}}^{1,1}}^{2}+D_{2}(\norm{\phi_{0}}_{\dot{\mathcal{F}}^{1,1}})\norm{\phi_{0}}_{\dot{\mathcal{F}}^{1,1}}^{2}\biggr)\biggr) \nonumber \\
    &\cdot\frac{\norm{\phi_{0}}_{\dot{\mathcal{F}}^{1,1}}}{\Lambda(\norm{\phi_{0}}_{\dot{\mathcal{F}}^{1,1}})-\nu_{0}}. \nonumber
\end{align}
Hence, \(\mathcal{F}(\theta)(0)\) is bounded in time.

\subsection{Regularization Argument}
In this Section, we outline an argument to prove Theorem \ref{mainthm}. We first regularize the original evolution equations for the interface for each \(N \in \mathbb{N}\). The sequence of solutions to these regularized equations yields a solution to the original evolution equations for the interface. To obtain a solution to each regularized equation, we employ Picard's theorem in the Banach space setting as stated in Theorems 3.1 and 3.3 of~\cite{majda2002vorticity}. To obtain a candidate for a solution to the original evolution equations for the interface, we use the Aubin-Lions lemma as stated in Lemma 7.5 of~\cite{gancedo2019global}.

\subsubsection{Regularized Equations for Interface Dynamics}
We recall that, under HLS parametrization, the dynamics of the interface are governed by
\begin{align}
    &\theta_{t}(\alpha) = \frac{2\pi}{L(t)}(U_{\alpha}(\theta)(\alpha)+T(\theta)(\alpha)(1+\theta_{\alpha}(\alpha))), \label{originalsystemstart}\\
    &L(t) =2\pi R \biggl(1+\frac{1}{2\pi}\mbox{Im}\int_{-\pi}^{\pi}\int_{0}^{\alpha}e^{i(\alpha-\eta)}\sum_{n \geq 1}\frac{i^{n}}{n!}(\theta(\alpha)-\theta(\eta))^{n}d\eta d\alpha\biggr)^{-\frac{1}{2}}.
\end{align}
For each \(N \in \mathbb{N}\), we define the regularized ordinary differential equation for the interface
\begin{align*}
    \frac{d\theta_{N}}{dt} = (\mathcal{J}_{N}^{1} \circ G_{N})(\theta_{N}),
\end{align*}
where \(\mathcal{J}_{N}^{1}\) is the high frequency cut-off operator introduced in (\ref{highfreqcutoffop}) and
\begin{align*}
    G_{N}(\theta_{N}) = &R^{-1}\biggl(1+\frac{1}{2\pi}\mbox{Im}\int_{-\pi}^{\pi}\int_{0}^{\alpha}e^{i(\alpha-\eta)}\sum_{n \geq 1}\frac{i^{n}}{n!}(\theta_{N}(\alpha)-\theta_{N}(\eta))^{n}d\eta d\alpha\biggr)^{\frac{1}{2}} \\
    &\cdot\biggl((U_{\alpha})_{N}(\theta_{N})+T_{N}(\theta_{N})\biggl(1+(\theta_{N})_{\alpha}\biggr)\biggr),
\end{align*}
in which
\begin{align}
    (U_{\alpha})_{N}(\theta_{N})(\alpha) &= (\mathcal{J}_{N}\circ\mbox{Re})\biggl(W(\theta_{N})(\alpha)\biggr) \nonumber \\
    U_{N}(\theta_{N})(\alpha) &=(\mathcal{J}_{N}\circ\mbox{Re})\biggl(V(\theta_{N})(\alpha)\biggr) \nonumber \\
    T_{N}(\theta_{N})(\alpha) &= \mathcal{M}\biggl(\biggl(1+(\theta_{N})_{\alpha}(\alpha)\biggr)U_{N}(\theta_{N})(\alpha)\biggr), \nonumber
\end{align}
\(\mathcal{J}_{N}\) is the high frequency cut-off operator defined in (\ref{highfreqcutoffop2}), and
\begin{align*}
    V(\theta_{N})(\alpha) = &\sum_{j=1}^{7}\frac{\gamma}{4\pi}\int_{-\pi}^{\pi}E_{j}(\theta_{N})(\alpha, \beta)d\beta \\
    &+\sum_{j=1}^{8}\frac{\gamma}{4\pi}\int_{-\pi}^{\pi}\widetilde{B_{j}}(\theta_{N})(\alpha,\beta)d\beta+\frac{\gamma}{4\pi}\int_{-\pi}^{\pi}\widetilde{B_{13}}(\theta_{N})(\alpha,\beta)d\beta, \\
    W(\theta_{N})(\alpha) =& V_{\alpha}(\theta_{N})(\alpha) \\
    =&\sum_{j=1}^{7}\frac{\gamma}{4\pi}\int_{-\pi}^{\pi}(E_{j})_{\alpha}(\theta_{N})(\alpha, \beta)d\beta \\
    &+\sum_{j=1}^{8}\frac{\gamma}{4\pi}\int_{-\pi}^{\pi}(\widetilde{B_{j}})_{\alpha}(\theta_{N})(\alpha,\beta)d\beta+\frac{\gamma}{4\pi}\int_{-\pi}^{\pi}(\widetilde{B_{13}})_{\alpha}(\theta_{N})(\alpha,\beta)d\beta.
\end{align*}

\subsubsection{Applying Picard's Theorem}
\label{applyingpicardtheorem}
We need to specify an appropriate Banach space for Picard's theorem. Due to HLS parametrization,
\begin{align}
    &\int_{-\pi}^{\pi}e^{i(\alpha+\phi(\alpha,t))}d\alpha=0, \label{originalsystemend}
\end{align}
which constrains the \(\pm 1\) Fourier modes of \(\phi(\alpha,t)\) to be completely determined by the rest of its nonzero Fourier modes. For this reason, we seek from the outset a solution whose \(\pm 1\) Fourier modes remain zero in time. For each \(N \in \mathbb{N}\), let
\begin{align*}
    H_{N}^{m} = \biggl\{f \in H^{m}([-\pi,\pi)) : \supp(\hat{f}) \subseteq [-N,N], \mbox{ \(\hat{f}(\pm 1) = 0\), } \mbox{Im\((f) = 0\)}\biggr\},
\end{align*}
which is a Banach space. For \(M>0\), let
\begin{align*}
    O^{M} = \{f \in H_{N}^{m} : \norm{f}_{H^{m}} < M\}.
\end{align*}
For sufficiently small \(M>0\), the conditions for Picard's theorem are met with \(B=H_{N}^{m}\), \(O=O^{M}\), and \(F=\mathcal{J}_{N}^{1} \circ G_{N}\). By Picard's theorem, we obtain that, for any \(\theta_{N,0} \in O^{M}\), there exists a time \(T_{N}>0\) such that the ordinary differential equation
\begin{align}
    \frac{d\theta_{N}}{dt} &= (\mathcal{J}_{N}^{1} \circ G_{N})(\theta_{N}), \label{ode} \\
    \theta_{N}(0) &= \theta_{N,0} \in O^{M} \nonumber
\end{align}
has a unique local solution \(\theta_{N} \in C^{1}([0,T_{N});O^{M})\).

\subsubsection{Derivation of an \emph{a priori} Estimate}
\label{keyaprioriestimateregularization}
For each \(N \in \mathbb{N}\), define \(\phi_{N}(\alpha,t) = \theta_{N}(\alpha,t) - \hat{\theta}_{N}(0,t)\).
Then we may write
\begin{align} \label{truncevol}
    \frac{d\theta_{N}}{dt} = \mathcal{L}_{N}(\theta_{N})+\mathcal{N}_{N}(\theta_{N}),
\end{align}
where \(\mathcal{L}_{N}(\theta_{N})\) and \(\mathcal{N}_{N}(\theta_{N})\) are the parts of the right hand side of (\ref{ode}) which are linear and superlinear in the variable \(\theta_{N}\), respectively.
The estimates presented in Sections \ref{estimationU1}, \ref{estimationU2}, and \ref{estimationU2deriv} and Lemmas \ref{estimatingterm2lem} and \ref{estimatingterm3lem} can be used to bound the right hand side of (\ref{truncevol}).
Ultimately, we obtain
\begin{align*}
    \frac{d}{dt}\normx{\phi_{N}}{1} &\leq -\biggl(\widetilde{\Lambda}(\normx{\phi_{N}}{1})-\nu'(t)\biggr)\normx{\phi_{N}}{2},
\end{align*}
where
\begin{align} \label{lambdaexpression}
    \widetilde{\Lambda}(\normx{\phi_{N}}{1}) =& \pi\frac{2}{R}\frac{\gamma}{4\pi}-\frac{\gamma}{4\pi}\frac{1}{R}A(\normx{\phi_{N}}{1})\normx{\phi_{N}}{1} \\
    &-\frac{1}{R}\frac{1}{A_{1}(\normx{\phi_{N}}{1})} \nonumber \\
    &\cdot\biggl(R_{1}(\normx{\phi_{N}}{1})\normx{\phi_{N}}{1}+R_{2}(\normx{\phi_{N}}{1})\normx{\phi_{N}}{1} \nonumber \\
    &+R_{3}(\normx{\phi_{N}}{1})\normx{\phi_{N}}{1}^{2}+R_{4}(\normx{\phi_{N}}{1})\normx{\phi_{N}}{1}^{2} \nonumber \\
    &+R_{5}(\normx{\phi_{N}}{1})\normx{\phi_{N}}{1} \nonumber \\
    &+3\left(H_{3}\normx{\phi_{N}}{1}+H_{4}\normx{\phi_{N}}{1}\right) \nonumber \\
    &+3\left(D_{1}(\normx{\phi_{N}}{1})\normx{\phi_{N}}{1}^{2}+D_{2}(\normx{\phi_{N}}{1})\normx{\phi_{N}}{1}^{2}\right) \nonumber \\
    &\cdot\left(1+2\normx{\phi_{N}}{1}\right) \nonumber \\
    &+\left(D_{1}(\normx{\phi_{N}}{1})\normx{\phi_{N}}{1}+D_{2}(\normx{\phi_{N}}{1})\normx{\phi_{N}}{1}\right) \nonumber \\
    &\cdot\left(1+2\normx{\phi_{N}}{1}\right) \nonumber \\
    &+6\normx{\phi_{N}}{1}\left(H_{3}\normx{\phi_{N}}{1}+H_{4}\normx{\phi_{N}}{1}\right) \nonumber \\
    &+2\left(H_{3}\normx{\phi_{N}}{1}+H_{4}\normx{\phi_{N}}{1}\right)\biggr). \nonumber
\end{align}
Let \(\theta^{0} \in \dot{\mathcal{F}}^{1,1}\), \(\mbox{Im}(\theta^{0}) = 0\), such that \(\widetilde{\Lambda}(\norm{\theta^{0}}_{\dot{\mathcal{F}}^{1,1}})>0\).
To ensure that \(\widetilde{\Lambda}(\norm{\theta^{0}}_{\dot{\mathcal{F}}^{1,1}})\) is well-defined, we choose \(\norm{\theta^{0}}_{\dot{\mathcal{F}}^{1,1}}\) small enough so that all of the geometric series contained in this expression converge. We further require that
\begin{align*}
    \abs{\mathcal{F}(\theta^{0})(0)}+Y\biggl(\norm{\theta^{0}}_{\dot{\mathcal{F}}^{1,1}}\biggr) + 2^{1/2}\norm{\theta^{0}}_{\dot{\mathcal{F}}^{1,1}}<M,
\end{align*}
where the function \(Y\) is defined in (\ref{uniformboundfunction}). This condition ensures that the solution to the original evolution equations for the interface is global in time. For each \(N \in \mathbb{N}\), let \(\theta_{N,0} = \mathcal{J}_{N}^{1}\theta^{0} \in H_{N}^{m}\). Letting
\begin{align*}
    \phi^{0} &= \theta^{0} - \hat{\theta^{0}}(0) \\
    \phi_{N,0} &= \theta_{N,0} - \hat{\theta_{N,0}}(0),
\end{align*}
we choose \(\nu_{0}\) such that \(0<\nu_{0}<\Lambda(\norm{\phi^{0}}_{\dot{\mathcal{F}}^{1,1}})<\Lambda(\norm{\phi_{N,0}}_{\dot{\mathcal{F}}^{1,1}})\).
It follows that for all \(t \in [0,T_{N})\),
\begin{align*}
    \normx{\phi_{N}(t)}{1} + \biggl(\Lambda(\norm{\theta^{0}}_{\dot{\mathcal{F}}^{1,1}})-\nu_{0}\biggr)\int_{0}^{t}\normx{\phi_{N}(\tau)}{2}d\tau & \leq \norm{\theta^{0}}_{\dot{\mathcal{F}}^{1,1}}.
\end{align*}
Then, for all \(t \in [0,T_{N})\),
\begin{align*}
    \frac{d}{dt}\normx{\phi_{N}}{1} \leq&-\biggl(\Lambda(\norm{\phi_{N,0}}_{\dot{\mathcal{F}}^{1,1}})-\nu_{0}\biggr)\normx{\phi_{N}}{1},
\end{align*}
which implies that \(\normx{\phi_{N}}{1}\) decays exponentially on \([0,T_{N})\).

\subsubsection{A Remark on the Solution Being Global in Time}
The global-in-time nature of the solution to the original evolution equations for the interface is inherited from that of the solutions to the regularized equations. The latter is a consequence of the continuation property of Picard's theorem in the Banach space setting and that the zeroth Fourier mode of \(\theta_{N}\) is bounded in time.

\subsubsection{Applying Aubin-Lions' Lemma}
To apply Aubin-Lions' lemma, we set \(X_{0} = \dot{\mathcal{F}}_{\nu}^{2,1}\), \(X = \dot{\mathcal{F}}_{\nu}^{1,1}\), \(X_{1} = \dot{\mathcal{F}}_{\nu}^{0,1}\), \(p=\infty\), and let \(G=\{\theta_{N} : N \in \mathbb{N}\}\). By Aubin-Lions' lemma, \(G\) is relatively compact in \(L^{2}([0,T];\dot{\mathcal{F}}_{\nu}^{1,1})\). This means that there exists a subsequence convergent in \(L^{2}([0,T];\dot{\mathcal{F}}_{\nu}^{1,1})\). That is, there exists \(\theta \in L^{2}([0,T];\dot{\mathcal{F}}_{\nu}^{1,1})\) such that \(\theta_{N} \to \theta\) in \(L^{2}([0,T];\dot{\mathcal{F}}_{\nu}^{1,1})\) as \(N \to \infty\), where for notational convenience we continue to use \(\theta_{N}\) to denote the subsequence. Our application of Aubin-Lions' lemma provides a candidate for a solution to the original problem, but it remains silent on the dynamics of \(\mathcal{F}(\theta(t))(0)\). Part of our task is to specify its dynamics. We first articulate the sense in which \(\theta\) is to become a solution.
\begin{definition} \label{weaksolution}
    We say that \(\theta \in L^{\infty}([0,T];\dot{\mathcal{F}}_{\nu}^{1,1}) \cap L^{1}([0,T];\dot{\mathcal{F}}_{\nu}^{2,1})\) is a weak solution of (\ref{originalsystemstart}) through (\ref{originalsystemend}) if for any \(\psi \in C_{0}^{\infty}([-\pi,\pi) \times [0,T])\),
    \begin{align*}
        &\int_{-\pi}^{\pi}\theta(\alpha,T)\psi(\alpha,T)d\alpha -\int_{-\pi}^{\pi}\theta(\alpha,0)\psi(\alpha,0)d\alpha-\int_{-\pi}^{\pi}\int_{0}^{T}\theta(\alpha,t)\psi_{t}(\alpha,t)dt d\alpha \\
        =&\int_{-\pi}^{\pi}\int_{0}^{T}R^{-1}\biggl(1+\frac{1}{2\pi}\mbox{Im}\int_{-\pi}^{\pi}\int_{0}^{\alpha}e^{i(\alpha-\eta)}\sum_{n \geq 1}\frac{i^{n}}{n!}(\theta(\alpha,t)-\theta(\eta,t))^{n}d\eta d\alpha\biggr)^{1/2}\\
        &\cdot\biggl(U_{\alpha}(\theta)(\alpha,t)+T(\theta)(\alpha,t)(1+\theta_{\alpha}(\alpha,t))\biggr)\psi(\alpha,t)dt d\alpha.
    \end{align*}
\end{definition}
Let \(\phi(t) = \theta(t) - \mathcal{F}(\theta(t))(0)\). We specify the dynamics of \(\mathcal{F}(\theta)(0)\) by requiring that
\begin{align}
    \frac{d}{dt}\mathcal{F}(\theta)(0) = \frac{2\pi}{L(\phi)}\biggl(U_{\alpha}(\phi)+T(\phi)(1+\phi_{\alpha})\biggr) - \frac{d}{dt}\phi \label{dynamicseqnforzeromode}
\end{align}
with the initial condition \(\mathcal{F}(\theta(0))(0) = \mathcal{F}(\theta^{0})(0)\). The dynamics equation (\ref{dynamicseqnforzeromode}) for \(\mathcal{F}(\theta)(0)\) implies that \(\theta\) is indeed a solution to the original problem in the sense of Definition \ref{weaksolution}.

\subsubsection{Inheritance of the \emph{a priori} Estimate}
At the end of Section \ref{keyaprioriestimateregularization}, we obtained that for all \(t \in [0,T]\),
\begin{align*}
    \normx{\phi_{N}(t)}{1} + \biggl(\Lambda(\norm{\theta^{0}}_{\dot{\mathcal{F}}^{1,1}})-\nu_{0}\biggr)\int_{0}^{t}\normx{\phi_{N}(\tau)}{2}d\tau & \leq \norm{\theta^{0}}_{\dot{\mathcal{F}}^{1,1}}.
\end{align*}
By Fatou's lemma, for any \(t \in [0,T]\),
\begin{align*}
    \int_{0}^{t}\liminf_{N \to \infty}\norm{\phi_{N}(\tau)}_{\dot{\mathcal{F}}_{\nu}^{2,1}}d\tau \leq \liminf_{N \to \infty}\int_{0}^{t}\norm{\phi_{N}(\tau)}_{\dot{\mathcal{F}}_{\nu}^{2,1}}d\tau.
\end{align*}
Using that
\begin{align*}
    \liminf_{N \to \infty}\normx{\phi_{N}(t)}{1} = \normx{\phi(t)}{1}, \\
    \liminf_{N \to \infty}\normx{\phi_{N}(t)}{2} = \normx{\phi(t)}{2},
\end{align*}
we then obtain that for all \(t \in [0,T]\)
\begin{align*}
    &\normx{\phi(t)}{1}+\biggl(\Lambda(\norm{\theta^{0}}_{\dot{\mathcal{F}}^{1,1}})-\nu_{0}\biggr)\int_{0}^{t}\normx{\phi(\tau)}{2}d\tau \\
    \leq&\liminf_{N \to \infty}\normx{\phi_{N}(t)}{1}+\biggl(\Lambda(\norm{\theta^{0}}_{\dot{\mathcal{F}}^{1,1}})-\nu_{0}\biggr)\liminf_{N\to \infty}\int_{0}^{t}\normx{\phi_{N}(\tau)}{2}d\tau \\
    \leq&\liminf_{N \to \infty}\biggl(\normx{\phi_{N}(t)}{1}+\biggl(\Lambda(\norm{\theta^{0}}_{\dot{\mathcal{F}}^{1,1}})-\nu_{0}\biggr)\int_{0}^{t}\normx{\phi_{N}(\tau)}{2}d\tau\biggr) \\
    \leq&\norm{\theta^{0}}_{\dot{\mathcal{F}}^{1,1}}.
\end{align*}
As a consequence,
\begin{align*}
    \theta \in L^{\infty}([0,T];\dot{\mathcal{F}}_{\nu}^{1,1})\cap L^{1}([0,T];\dot{\mathcal{F}}_{\nu}^{2,1})
\end{align*}
and \(\normx{\phi(t)}{1}\) decays exponentially on \([0,T]\).

\subsubsection{Continuity in Time and Instantaneous Analyticity}
That \(\theta \in C([0,T];\dot{\mathcal{F}}_{\nu}^{1,1})\) follows from arguments similar to those outlined in~\cite{gancedo2020surface}. That \(\theta\) is instantaneously analytic is due to the exponential weight in the \(\normx{\cdot}{1}\) norm of \(\theta(\cdot, t)\) for any \(t>0\).

\section{Uniqueness}
Let \(\theta_{1}\) and \(\theta_{2}\) be two solutions to the original problem with the same initial datum, where the \(\pm 1\) Fourier modes remain zero in time. We may write
\begin{align*}
    \norm{\theta_{1}-\theta_{2}}_{\dot{\mathcal{F}}^{1,1}} = \sum_{k \neq 0}\abs{k}\abs{\mathcal{F}(\theta_{1}-\theta_{2})(k)} = 2\sum_{k>0}\abs{k}\abs{\mathcal{F}(\theta_{1}-\theta_{2})(k)}.
\end{align*}
Differentiating this equation with respect to time and then using estimates analogous to those contained in earlier Sections, we can obtain that, for sufficiently small \(\norm{\theta^{0}}_{\dot{\mathcal{F}}^{1,1}}\),
\begin{align*}
    &\frac{d}{dt}\norm{\phi_{1}-\phi_{2}}_{\dot{\mathcal{F}}^{1,1}} \leq \mathcal{E}\norm{\phi_{1}-\phi_{2}}_{\mathcal{F}^{1,1}},
\end{align*}
where \(\mathcal{E}\) is a coefficient that depends on \(\norm{\phi_{1}}_{\dot{\mathcal{F}}^{1,1}}\), \(\norm{\phi_{2}}_{\dot{\mathcal{F}}^{1,1}}\), \(\norm{\phi_{1}}_{\dot{\mathcal{F}}^{2,1}}\), and \(\norm{\phi_{2}}_{\dot{\mathcal{F}}^{2,1}}\) such that \(\mathcal{E}\) is integrable in time. Since the two solutions share the same initial datum, \(\phi_{1}=\phi_{2}\) by Gr\"onwall's inequality. Moreover, \(\mathcal{F}(\theta_{1})(0) = \mathcal{F}(\theta_{2})(0)\) since the dynamics of \(\mathcal{F}(\theta_{1})(0)\) and \(\mathcal{F}(\theta_{2})(0)\) are determined completely by \(\phi_{1}\) and \(\phi_{2}\), respectively, and they share the same initial condition \(\mathcal{F}(\theta^{0})(0)\).

\section{Acknowledgements}
This work was supported in part by the NSF under Grant DMS-2042144 (USA, awarded to Yoichiro Mori) and Grant DMS-2055271 (USA, awarded to Robert M. Strain).

\nocite{*}
\printbibliography

\end{document}